\documentclass[a4paper,12pt]{article}
\usepackage{amsfonts,amsmath,amssymb,amsthm}
\usepackage{txfonts}
\usepackage[utf8]{inputenc}
\usepackage{amsmath}
\usepackage{amssymb}
\usepackage{mathtools}
\usepackage{bbm}

\RequirePackage[colorlinks,citecolor=red,urlcolor=blue, linkcolor=blue]{hyperref}
\usepackage[margin=0.8in]{geometry}

\usepackage[dvips]{graphics}
\usepackage{epsfig,rotating}
\usepackage{tabularx}
\usepackage{bm}
\usepackage{bbm}
\usepackage{float}
\usepackage{mdframed}
\usepackage{lipsum}
\usepackage{enumitem}
\usepackage{tikz}
\usepackage{natbib}

\newtheorem{assumption}{Assumption}

\newtheorem{manualtheoreminner}{Theorem}
\newenvironment{manualtheorem}[1]{%
  \renewcommand\themanualtheoreminner{#1}%
  \manualtheoreminner
}{\endmanualtheoreminner}

\newtheorem{manuallemmainner}{Lemma}
\newenvironment{manuallemma}[1]{%
  \renewcommand\themanuallemmainner{#1}%
  \manuallemmainner
}{\endmanualtheoreminner}

\newtheorem{assump}{Assumption}
\newtheorem{remark}{Remark}%

\usepackage[nottoc,numbib]{tocbibind} 
\usepackage[toc,page]{appendix}

\newtheorem{thm}{Theorem}
\newtheorem{lemma}{Lemma}
\newtheorem{prop}{Proposition}
\newtheorem{corollary}{Corollary}

\newcommand{\beaa}{\begin{eqnarray*}}
\newcommand{\eeaa}{\end{eqnarray*}}
\newcommand{\bea}{\begin{eqnarray}}
\newcommand{\eea}{\end{eqnarray}}

\newcommand{\la}{\left\{}
\newcommand{\ra}{\right\}}
\newcommand{\lb}{\left(}
\newcommand{\rb}{\right)}
\newcommand{\mb}{\mathbb}
\newcommand{\ve}{\varepsilon}
\newcommand{\mc}{\mathcal}

\newcommand{\argmin}{\operatornamewithlimits{argmin}}

\begin{document}
\title{Time-uniform concentration bounds for iterative algorithms}
\author{Tuan Pham \\ \small{Department of Statistics and Data Science, University of Texas, Austin} \\  \small{\href{mailto:tuan.pham@utexas.edu}{tuan.pham@utexas.edu} }
\and Alessandro Rinaldo \\ \small{Department of Statistics and Data Science, University of Texas, Austin} \\  \small{\href{mailto:alessandro.rinaldo@austin.utexas.edu}{alessandro.rinaldo@austin.utexas.edu}}
\and Purnamrita Sarkar \\ \small{Department of Statistics and Data Science, University of Texas, Austin} \\  \small{\href{mailto:purna.sarkar@austin.utexas.edu}{purna.sarkar@austin.utexas.edu}}}
\date{}

\maketitle

\begin{abstract}
We develop a new framework for deriving   time-uniform concentration bounds for the output of stochastic sequential algorithms  satisfying certain recursive inequalities akin to those defining the almost-supermartingale processes introduced by \cite{robbins1971convergence}. Our approach is of wide applicability, and can be deployed in settings in which exponential supermartingale processes, required by prevailing methodologies for anytime-valid concentration inequalities, are not readily available. 
Our results can be viewed as quantitative versions of the classical Robbins–Siegmund Lemma.
 We demonstrate the effectiveness of our method by providing new and optimal time-uniform concentration bounds for Oja’s algorithm for streaming PCA, stochastic gradient descent, and stochastic approximations.
\end{abstract}

\tableofcontents

\section{Introduction}


Modern statistical and machine learning modeling and techniques are often applied to datasets whose size and complexity call for scalable iterative stochastic algorithms, which are naturally deployed in a sequential and often data-driven manner. Concretely, suppose that we are interested in estimating a parameter $\bm{\theta}^*$ using an iterative algorithm that outputs a sequence of estimators $\{ \widehat{\bm{\theta}}_t \}_{t=0,1,2,\ldots}$ computed in a sequential manner - that is, $\hat{\bm{\theta}}_t$ is an update of $\widehat{\bm{\theta}}_{t-1}$ after acquiring or processing a new observation. At each iterate $t$, the accuracy of the algorithm is measured by the (unobservable) quantity
\[
\mathrm{Loss}\!\left(\hat{\bm{\theta}}_t, \bm{\theta_*} \right), 
\]
for an appropriate loss function $\mathrm{Loss}(\cdot,\cdot)$. The standard approach in the literature to address the convergence properties of $\{ \widehat{\bm{\theta}}_t \}_{t=0,1,2,\ldots}$ is to establish high-probability and concentration bounds for the loss that hold at a fixed and large deterministic number of iterations. Specifically, for any given value of the iterate $T$ (say larger than some value $T_0$) and a probability parameter $\delta \in (0,1)$, the corresponding {\it fixed-time} concentration bound takes the form
\begin{equation}
    \label{eq:fixed-time}
\quad \mathbb{P} \!\left( \mathrm{Loss}\!\left(\hat{\bm{\theta}}_T, \bm{\theta_*} \right)  \geq r_{\mathrm{fixed-time}}(T,\delta)\right) \leq \delta, \quad \forall T \geq T_0, 
\end{equation}
with $r_{\mathrm{fixed-time}}(T,\delta)$ a monotonically vanishing function of $T$ that depends on $\delta$ and properties of the data-generating distribution (e.g., the dimension of the support). 

While useful in providing a certificate of convergence,  fixed-time concentration bounds however offer only a point-wise (in time) and thus weak guarantee. Instead, as illustrated in the recent and growing literature on any-time bounds and confidence sequences, it is more desirable and statistically safer  to establish {\it time-uniform} or {\it any-time} concentration bounds of the form 
   \begin{equation}
    \label{eq:any-time}
\quad \mathbb{P} \!\left( \mathrm{Loss}\!\left(\hat{\bm{\theta}}_t, \bm{\theta_*} \right)  \geq r_{\mathrm{any-time}}(t,\delta), \forall t \geq T_0 \right) \leq \delta, 
    \end{equation}
where the, deterministic or random, time boundary function $t \mapsto r_{\mathrm{any-time}}(t,\delta)$ is also vanishing in $t$.  Any-time concentration bounds deliver signifcant  algorithmic and statistical advantages over fixed-time ones. 

First, they yield a stronger and arguably more natural notion of algorithmic convergence, as they ensure that the sequence of estimators $\{ \widehat{\bm{\theta}}_t \}_{t=0,1,2,\ldots}$ will remain close to the target parameter simultaneously over all large iterates $t$, with high probability. While this is presumably how practitioners would typically regard the convergence behavior of a stochastic algorithm, such interpretation is not warranted by the fixed-time guarantee \eqref{eq:fixed-time}. 

Secondly, from a statistical viewpoint, any-time bounds enjoy remarkable features over traditional inference methods that make them well suited to handle data adaptive decisions, e.g. when to terminate the data collection process, stop the algorithm, fine tune and adjust model parameters on the fly based on the values observed so far. Specifically, they are valid under arbitrary stopping rules { \cite{waudby2024estimating,waudby2024time}, i.e. they remain valid at any data-dependent stopping time. What in more, they also enjoy {\it post-hoc validity } at arbitrary random (not just stopping) times. This is a crucial property in practical scenarios where data collection -- in the case of sequential algorithms, the estimation task itself -- may be interrupted or retrospectively adjusted. For example, practitioners may face unforeseen budget cuts requiring them to terminate an experiment at an unplanned or random time. Despite this change from the initial plan,  any-time concentration bounds remain valid. 
Finally, any-time bounds are adaptive to random sample sizes: one may stop early when a convergence criterion is met or continue beyond the planned sample size, all while maintaining correct coverage and without requiring any corrections. 
To reiterate, any-time bounds provide broad protection against data-adaptive decision.  We emphasize that ``data-peeking'' practices of this sort are widespread and, arguably, a fairly  natural way to carry out data analytic tasks and sequential estimation, to formulate and test scientific hypothesis and to design or adjust experiments. However, they are incompatible with and, in fact, often invalidate conventional statistical methodologies that are not designed to be data-adaptive.

Any-time concentration bounds are certainly not new in concept:  they are rooted in statistical sequential analysis \citep{robbins1970statistical,lai1976confidence} and have been deployed extensively in the study of bandit algorithms. However, due to  recent breakthroughs in the construction of time-uniform, martingale-based concentration inequalities due to \cite{howard2021time,waudby2024time,howard2020time,waudby2024estimating,10.1109/TIT.2023.3250099}, the past few years have witnessed a  
flurry of new results and applications on a multitude of problems, from  mean estimation, A/B testing, reinforcement learning, and bandit optimization \cite{howard2021time,waudby2024time,waudby2024estimating,howard2020time}; see also Section \ref{related results} for more details.  

While powerful and broadly applicable, this current prevailing approach relies in a fundamental way on the availability of martingale and exponential supermartingale processes (and maximal inequalities thereof), which combines naturally  with the use the Cramer-Chernoff method for deriving high probability bounds.
For the puspose of studying sequential algorithms, these martingale processes are constructed based on the sequence of values $\{\mathrm{Loss}\!\left(\hat{\bm{\theta}}_t, \bm{\theta_*} \right) \}_{t=0,1,\ldots,}$ or related quantities. 
While this is possible in many problems, there are important scenarios in which constructing an exponential supermartingale in this manner may be difficult or suboptimal. 
A primary example is Oja's algorithm for online PCA, which is used to compute principal eigenvectors of an unknown covariance matrix $\bm{\Sigma}$ in a sequential manner \cite{oja1982simplified,oja1985stochastic} based on, say, i.i.d. samples $\{ \bm{X}_t\}_{t \geq 1}$, revealed one at a time. In detail, algorithm is initialized with a random unit vector \( \hat{\bm{v}}_0 \sim \text{Uni} \big( \mathbb{S}^{p-1} \big) \), and, for a given a sequence of step sizes \( \{ \eta_t \}_{t \geq 1} \), implements  the online update rule 
\begin{equation}
    \label{eq:oja.intro}
\hat{\bm{v}}_{t} := \big( \bm{I}_p + \eta_{t}\bm{X}_{t} \bm{X}_{t}^{\top} \big) \hat{\bm{v}}_{t-1}, \quad \hat{\bm{v}}_{t} \leftarrow \frac{\hat{\bm{v}}_{t}}{\|\hat{\bm{v}}_{t}\|}.
\end{equation}
For this problem, one may consider the function $\mbox{Loss} \lb \hat{\bm{v}}_{t},\bm{v}_0 \rb=\sin^2\lb \hat{\bm{v}}_{t},\bm{v}_0 \rb$, where $\bm{v}_0$ is the leading eigenvector of $\bm{\Sigma}$.
It is well known that, under mild conditions, Oja's iterates $\{ \hat{\bm{v}}_{t} \}_{t \geq 0}$ converge to $\bm{v}_0 $ \citep{oja1982simplified,oja1985stochastic}. The key advantage of Oja's algorithm, which operates in an online manner by processing one observation at a time, over batch methods is its memory and computational efficiency: it requires only \( O(p) \) memory and $O(np)$ time, where \( p \) is the dimension of the data.
The convergence properties of Oja's iterates  have the been extensively studied over the past decades. In particular, assuming an i.i.d. and Markovian sequence of data points from a centered distribution with covariance matrix $\bm{\Sigma}$, Oja's iterate at a {\it fixed number of iterations} has been shown to concentrate tightly around the principal eigenvector of $\bm{\Sigma}$. See, e.g., \cite{huang2021streaming,kumar2023streaming} and Section~\ref{PCA} for additional references and a short discussion.

Generalizing existing sharp fixed-time concentration bounds for Oja's algorithm to  time-uniform 
bounds, i.e. holding simultaneously over an infinite time course,  appears to be non-trivial. The main difficulty lies in the form \eqref{eq:oja.intro} of the updates, which consist of matrix product updates. Unlike matrix addition, matrix products do not not directly lead to an exponential supermartingale process, a representation that enables concentration. Indeed, the best known  matrix product concentration bounds, due to \cite{matrix.product}, were obtained using different means, namely uniform smoothness properties of the Schatten trace classes. As a result, to the best of our knowledge, techniques to derive sharp any-time concentration bounds for Oja's algorithm are not available in the literature.

In this article, we propose a new technique for constructing anytime-valid concentration bounds for sequential stochastic algorithms for which the loss function at the current iterate satisfies certain recursive inequalities that can be viewed as an adaptation of those defining
{\it almost supermartingale} processes put forward in the seminal paper \citep{robbins1971convergence}; see \eqref{almost-supermartingale} for a precise definition. Thus, our settings do not require a well-defined martingale or supermartingale structure. In particular, Oja's algorithm is one example falling within this class. 
The key difference from the original almost supermartingale process of Robbins and Siegmund is that the form of the recursion we assume will guarantee not only almost sure convergence of the process, but almost sure convergence to zero. Remarkably, this simple adaptation is enough to deduce any-time concentration bounds. 

Let us briefly introduce our framework, deferring a rigorous formulation to Section~\ref{sec: main results}.  
We consider a nonnegative process $\la L_t; t \geq 0 \ra$ adapted to some filtration $\la \mathcal{F}_t; t \geq 0 \ra$ satisfying a recursive inequality of the form  
\begin{align} \label{rec constraint}
L_t \leq (1 - \eta_t) L_{t-1} + U_t,
\end{align}
where $\la \eta_t; t \geq 1 \ra$ is a deterministic sequence of {\it stepsizes} taking values in $(0,1)$, and $\{ U_t, t \geq 0 \}$ is an adapted noise process whose magnitude can be controlled by polynomials in the step sizes and the previous iterate, i.e. such that 
\begin{align} \label{hierachy}
 \begin{cases}
      \left|  \mathbb{E} \left[ U_t \,\middle|\, \mathcal{F}_{t-1} \right]  \right| &\lesssim \displaystyle\sum_{i=1}^{m} \eta_t^{a_i} L_{t-1}^{b_i}, \\[0.5em]
      |U_t|  &\lesssim  \displaystyle\sum_{i=1}^{m} \eta_t^{c_i} L_{t-1}^{d_i},
 \end{cases}  
\end{align}
where $\{a_i, b_i, c_i, d_i\}_{1 \leq i \leq m}$ are positive constants and the symbol $\lesssim$ indicates equality up to some deterministic quantities, depending on the problem parameters. One should think of $L_t$ as the value of the loss function at the $t$th iterate of a sequential algorithm, i.e. in the notation introduced above, $L_t =  \mathrm{Loss}\!\left(\hat{\bm{\theta}}_t, \bm{\theta_*} \right)$. 
Inequality \eqref{rec constraint} quantifies how much the loss changes at time $t$, after one iteration with step size $\eta_t$. The noise process $U_t$ can be thought of as an extra term that prevents $L_t$ from being a supermartingale and needs to be controlled. Concrete examples falling under this framework -- that is, satisfying conditions \eqref{rec constraint} and \eqref{hierachy} -- include Oja's algorithm, stochastic gradient descent and  Robbins-Monroe scheme for stochastic approximations, which we will analyze in Section~\ref{sgd}. A general strategy used in the literature on the convergence analysis of iterative algorithms is to leverage inequalities of the form \eqref{rec constraint} and \eqref{hierachy} to control the process $U_t$ appropriately and finally demonstrate that the loss function decreases in probabilistic sense. In this paper we take a further step and formulate of a general technique to deduce not only stochastic convergence to zero of the loss process but also any-time concentration bounds for processes satisfying \eqref{rec constraint} and \eqref{hierachy}. Specifically, we make the following contributions.
 


\begin{itemize}
   \item In Theorem~\ref{master} we show that the information provided by \eqref{rec constraint} and \eqref{hierachy} is sufficient to yield a time-uniform concentration bound vanishing at the rate $O(\log \log t / t)$, using step sizes $\eta_t = \Omega(1/t)$, provided that $L_0$ is sufficiently small with high probability. 
\item The LIL-style $O(\log \log t / t)$ concentration rate is optimal among all processes that satisfy \eqref{rec constraint} and \eqref{hierachy} with step sizes $\eta_t = \Omega(1/t)$. This is the content of Theorem~\ref{lower bound}. 

\item We apply our results to derive time-uniform bounds for SGD in the strongly convex case or under the Polyak–Lojasiewicz condition (see Section \ref{sgd}), two versions of Oja’s algorithm (see Section \ref{PCA}) and Robbins-Monroe approximation scheme (see Section \ref{Robbins-Monroe}). To the best of our knowledge, ours are the first time-uniform concentration bounds for the output of Oja-style algorithms. Our any-time concentration bounds are as sharp as the best known fixed-time bound established by \cite{huang2021streaming} for streaming PCA, save for the unavoidable $\log \log t/t$ term.  
\item Interestingly, when applied to SGD, for which martingale-based time-uniform concentration bounds up to a fixed time horizon already exist \citep[see, e.g.][]{rakhlin2011making}, our technique delivers slightly 
sharper bounds, in addition to being applicable to an infinite time course; see Section~\ref{explicit}.
\item Our analysis is tailored to problems with bounded noise, a setting that is commonly assumed in the literature. Extensions to the cases involving unbounded noise, in which the coefficients on the right-hand side in \eqref{hierachy} is random and has sub-Gaussian tail, are presented in supplement. 

\end{itemize}
The rest of the paper is organized as follows. Background and related results are discussed in Section~\ref{prelim}. The main results are presented in Section~\ref{sec: main results}, with applications to various algorithms given in Sections~\ref{sgd}, \ref{PCA}, and \ref{Robbins-Monroe}. The lower bound is presented in Section~\ref{Robbins-Monroe}. Conclusions, remarks, and directions for future research are provided in Section~\ref{sec:conclusion}. All proofs and technical results are deferred to the remaining sections of the paper.

\color{black}

\subsection{Related results} \label{related results}

The analysis of stochastic gradient algorithms has been studied extensively over the past few decades. Recent advances in this direction have established non-asymptotic bounds and optimality results in a variety of settings. An incomplete list of recent works on convergence and limit theorems for the last iterate includes \cite{carter2025statistical,duchi2021asymptotic,li2022root,chen2024online,madden2024high}.

In contrast, time-uniform convergence results are much less understood. To the best of the authors' knowledge, this topic has only been explored in \cite{rakhlin2011making,feng2023anytime}. The results in \cite{rakhlin2011making} provide time-uniform bounds for the iterates of SGD in a strongly convex setting over a fixed time horizon, while \cite{feng2023anytime} establishes related results for stochastic momentum algorithms using Ville’s inequality. To the best of our knowledge, no existing work provides a time-uniform analysis of Oja's algorithm. Even for the last iterate, most existing results are limited to the case of bounded data.

\section{Main results} \label{sec: main results}


\subsection{Preliminaries} \label{prelim}
In order to illustrate the rationale behind our formulation and results, it  is first useful to recall the classical Robbins-Siegmund's lemma \citep{robbins1971convergence}, a result about the almost-sure convergence of non-negative almost-supermartingales. In detail, suppose $\la L_t; t\geq 1 \ra$ is a sequence of non-negative stochastic process adapted to a filtration $\la \mathcal{F}_t; t \geq 1 \ra$ such that 
\begin{align} \label{almost-supermartingale}
    \mb{E} \lb L_{t+1} | \mathcal{F}_{t} \rb \leq (1+a_t) L_t + b_t - c_t
\end{align}
where $a_t, b_t, c_t$ are non-negative random weights that are adapted to $\mathcal{F}_t$. A stochastic process that satisfies \eqref{almost-supermartingale} is called an ``almost-supermatingale''. It is clear that when $a_t, b_t, c_t$ are identically zero, \eqref{almost-supermartingale} reduces to the classical concept of a supermartingale.
Another useful variant of \eqref{almost-supermartingale}, introduced in \cite{lai2006extended}, takes the form
\begin{align} \label{Lai-as}
    L_{t+1} \leq (1+ a_t) L_t + b_{t+1} - c_{t+1} + w_{t} \xi_t
\end{align}
where $\la \xi_t, \mathcal{F}_t \ra$ is a martingale difference sequence and $a_t, b_t, c_t, w_t$ are non-negative random weights adapted to $\mathcal{F}_t$. 

Robbins-Siegmund's lemma states that any almost-supermatingale in the sense of \eqref{almost-supermartingale} converges almost surely to a limit $X$ on the event 
\[
\la  \sum_{t=1}^{\infty} a_t < \infty, \sum_{t=1}^{\infty} b_t < \infty \ra.
\]
The same type of convergence (but with more relaxed assumptions) also holds for \eqref{Lai-as}, though without the requirement that $\sum_{t=1}^\infty b_t < \infty$. 

The most common applications of Robbin-Siegmund's lemma arise when $L_t:= \mc{L}(\hat{\theta}_t, \theta)$ for some loss function $\mc{L}$, and $\hat{\theta}_t$ being a sequence of estimator that is obtained via the stochastic gradient descent (SGD) algorithms; see, for examples, \cite{robbins1970statistical,robbins1971convergence,lai2006extended} and the references therein. It is often the case that $X_t$ serves as the Lyapunov function associated with the ODE of the SGD, which is obtained by letting the step sizes diminish to zero in the noiseless setting. However, as we will see later, there are examples when the Lyapunov functions are hard to construct explicitly, and other choices of $L_t$ are needed. 

For both formulations \eqref{almost-supermartingale} and \eqref{Lai-as} of the quasi-supermartingale process, the almost sure limit $X$ needs not be zero. Indeed, without any additional assumptions, this limit can be arbitrary. For example, consider a simple situation in the deterministic setting with 
\[
L_t:= \prod_{k=1}^{t} \lb 1 + \frac{1}{k^2} \rb.
\]
It is easy to check that $L_t$ converges to a positive limit. Moreover, by changing the terms $1/k^2$ to $Y_k/k^2$, for some sequence of random variables $\la Y_k, k \geq 1 \ra$, one can get any positive limiting distribution. 
To ensure almost sure convergence to a zero limit, we need slightly different conditions than \eqref{almost-supermartingale} and \eqref{Lai-as}. 
In particular, the ideal recursion should encode the information that $L_t$ is decreasing on average. For that reason, we will consider a specific class of non-negative, almost-supermartingale as follows.

\begin{lemma}[Simplified Robbins-Siegmund's lemma] \label{simplified lemma}
    Let $\la L_t \ra_{t\geq 1}$ be an adapted, non-negative process with respect to the filtration $\la \mathcal{F}_t ; t \geq 1 \ra$. Suppose there exists a positive, deterministic sequence $\la \eta_t; t\geq 1 \ra$ and a non-negative, adapted sequence $\beta_t$ such that 
    \begin{align} \label{simplified Robbins-Siegmund}
            \mb{E} \lb  L_{t+1}| \mathcal{F}_t \rb \leq (1-\eta_t)L_t + \beta_t
    \end{align}
    where $\la \eta_t, \beta_t \ra_{t\geq 1}$ satisfies
    \begin{align*}
       & \sum_{t=1}^{\infty} \eta_t = \infty ; \\
       & \sum_{t=1}^{\infty} \beta_t < \infty \ \text{almost surely}; \\
       & \lim_{t \to \infty} \beta_t/\eta_t = 0.
    \end{align*}
    Then, $\lim_t L_t = 0 $, almost surely. 
\end{lemma}
Let us compare  \eqref{simplified Robbins-Siegmund} with \eqref{almost-supermartingale} and \eqref{Lai-as}. In the formulation \eqref{simplified Robbins-Siegmund}, we have assumed that the negative term $c_t$ is at least  $-\eta_t L_t$. In other words, on average, the decrease from $L_{t}$ to $L_{t+1}$ should be a non-negligible proportion of $L_t$. Moreover, the extra gain obtained from $\beta_t$ is required to be of smaller order than $\eta_t$. 

In the applications we consider here, the sequence $\eta_t$ is the step size and $\beta_t$ is the error induced after one update of the algorithms. It turns out that recursion of the form \eqref{quantitative R-S lemma} or \eqref{simplified Robbins-Siegmund} hold for many classes of problems, which will be explained in Sections \ref{PCA}, \ref{sgd} and \ref{Robbins-Monroe}. Note that the assumptions in Lemma \ref{simplified lemma} implies that the sequence of step sizes $\eta_t$'s satisfies
\begin{align} \label{admissible}
    \sum_{t} \eta_t = \infty \ \text{and} \ \sum_{t} \eta_t^2 < \infty. 
\end{align}
It is easy to see that the best possible choice of step sizes such that \eqref{admissible} holds are those that are of order $\Omega(1/t)$. 

Below we show that if the recursive conditions \eqref{simplified Robbins-Siegmund} hold for process $\{ L_t \}$ itself -- and not just for the process of the conditional expectations $\{ \mb{E} \lb  L_{t}| \mathcal{F}_{t-1} \rb \}$  -- then it is possible to obtain time-uniform concentration bounds. This is the main contribution of the paper.


\subsection{A quantitative Robbins-Siegmund's Lemma and time-uniform concentration bounds} 

In this section, we formulate a quantitative version of \ref{simplified lemma} in the form of a anytime-valid concentration bounds for the entire process $\{ L_t \}$ vanishing at the optimal LIL rate $O\left( \frac{ \log\log t}{t} \right)$ rate. Towards that goal, we impose the following conditions.
{\begin{assumption}[Recursion]\label{assumption:main}
 Let $\la  L_t; t \geq 1 \ra$ be  a non-negative process and $\la  U_t; t \geq 0 \ra$ a {\it noise} process, both adapted to the filtration $\la  \mc{F}_t; t \geq 0  \ra$. Let $\la \eta_t; t \geq 1\ra$ be a positive deterministic sequence of stepsizes such that 
  $$\frac{C}{t} <  \eta_t < 1/ C$$ 
  for all $t \geq 1$ and some constant $C_1 > 0$. Almost surely, for some positive constants $C_1<C$, $C_2, C_3$ and $\la A_i; B_i a_i; b_i; c_i; d_i \ra_{i=1}^{m}$, the  recursive conditions
 \begin{align}
 &  L_{t} \leq (1-C_1\eta_t) L_{t-1} + U_t , \label{quantitative R-S lemma}\\
      &  \Big| \mb{E} \lb U_t | \mathcal{F}_{t-1}  \rb \Big| \leq  C_2 \cdot \eta_t^2 + \sum_{i=1}^{m} A_i \cdot \eta_t^{1+a_i} L_{t-1}^{b_i},    \label{cond-error} \\
       & |U_t| \leq C_3 \eta_t \sqrt{L_{t-1}} + \sum_{i=1}^{m} B_i \cdot \eta_t^{1/2+c_i} \cdot L_{t-1}^{d_i},  \label{as-bound}
    \end{align}
    hold. 
\end{assumption}}


    
    The settings in~\eqref{quantitative R-S lemma},~\eqref{cond-error}, and~\eqref{as-bound} may initially seem unnatural. However, they arise frequently in the analysis of iterative algorithms:

\begin{itemize}
    \item \textbf{Stochastic Gradient Descent (SGD)}: see Section~\ref{sgd} for a detailed discussion. Roughly speaking, if we define $L^{SGD}_t \equiv \|x_t - x^{*}\|^2$, then~\eqref{quantitative R-S lemma},~\eqref{cond-error}, and~\eqref{as-bound} hold in the form
    \begin{align*} 
        L^{SGD}_t \leq \left(1 - 2\lambda \eta_t \right) L^{SGD}_{t-1} + 2\eta_t Y^{SGD}_t + B^2 \eta_t^2,
    \end{align*}
    where $\mathbb{E} \left[Y^{SGD}_t \mid \mathcal{F}_{t-1} \right] = 0$, $|Y^{SGD}_t| \leq B \sqrt{L^{SGD}_{t-1}}$, and $B$ is a constant independent of $t$.

    \item \textbf{Oja's Algorithm}: a variant of SGD tailored for PCA, which will be discussed in detail in Section~\ref{PCA}. If we set $L^{Oja}_t \equiv \sin^2 \left(\hat{\bm{v}}_t, \bm{v} \right)$, where $\hat{\bm{v}}_t$ is the output of the algorithm and $\bm{v}$ is the true principal eigenvector, then
    \[
        L^{Oja}_t \leq \left(1 - 2\rho \eta_t \right) L^{Oja}_{t-1} + 2\rho \eta_t \lb L^{Oja}_{t-1} \rb^2 + Q^{Oja}_t + 5B^4 \eta_t^2 + 2\eta_t^3 B^6,
    \]
    where $\mathbb{E} \left[ Q^{Oja}_t \mid \mathcal{F}_{t-1} \right] = 0$, $|Q^{Oja}_t| \leq B^2 \sqrt{L^{Oja}_{t-1}}$, and $B$, $\rho$ are constants independent of $t$.

Another variant of the Oja's algorithm is the Krasulina-Oja's algorithm; see Section \ref{PCA} for more details. In the same setting as above, and with $L^{Kra}_t$ being the loss function between the iterates of Krasulina-Oja's algorithm and the true parameter, we have
    \[
        L^{Kra}_t \leq \left(1 - 2\rho \eta_t \right) L^{Kra}_{t-1} + 2\rho \eta_t \lb L^{Kra}_{t-1} \rb^2 + Q^{Kra}_t + 5B^4 \eta_t^2,
    \]
    where $\mathbb{E} \left[ Q^{Kra}_t \mid \mathcal{F}_{t-1} \right] = 0$, $|Q^{Kra}_t| \leq B^2 \sqrt{L^{Kra}_{t-1}}$, and $B$, $\rho$ are constants independent of $t$. Note that, unlike Oja's algorithm, this version does not include a term of order $\eta_t^3$.

    \item  \textbf{Robbins-Monro scheme:} One aims to find the root of a univariate function $M(x)$ based on noisy observations $Y(x)$; see Section~\ref{Robbins-Monroe} for details. Under the sub-polynomial condition on $M$, Proposition~\ref{robbins-monroe constraint} shows that the loss $L^{RM}_t:=\| x_t -x^* \|^2$ associated with the iterates of the Robbins-Monro scheme satisfies the recursion
    \[
    L^{RM}_t \leq \lb 1 - R_1 \eta_t \rb L^{RM}_{t-1} + Q^{RM}_t + \eta_t^2 \cdot P \lb L^{RM}_{t-1} \rb
    \]
    where $\mb{E} \lb Q^{RM}_t | \mathcal{F}_{t-1} \rb=0$, $|Q^{RM}_t| =O \lb \eta_t \sqrt{L^{RM}_{t-1}} \rb$ and $P$ is some polynomial with positive coefficients.
\end{itemize}

Let us explain the roles of the parameters appearing in \eqref{cond-error} and \eqref{as-bound}. The two terms \( C_2 \eta_t^2 \) in \eqref{cond-error} and \( C_3 \eta_t \sqrt{L_{t-1}} \) are ubiquitous; they arise in all examples considered in this paper. The remaining terms in \eqref{cond-error} and \eqref{as-bound} are problem-dependent and may or may not appear, depending on the specific context.

We will later see that in the simplest setting where \( L_t \) corresponds to the squared error loss of iterates produced by SGD in the strongly convex regime, equations~\eqref{cond-error} and~\eqref{as-bound} reduce to~\eqref{cond-error 2} and~\eqref{as-bound 2}, which will be studied in greater detail in Section~\ref{explicit}. However, for more complex problems, especially in the non-convex case, additional terms with different exponents can arise in~\eqref{cond-error} and~\eqref{as-bound}, as illustrated by the example of Oja's algorithm discussed above.

The following general result shows that even when the noise $U_t$ is extremely complicated, given a sufficiently good initialization \( L_0 \), one can always construct a time-uniform bound  with asymptotically optimal width, scaling as \( \log \log t / t \).

\begin{thm} \label{master}
Suppose the positive process $\{L_t\}_{t \geq 0}$ satisfies the recursion \eqref{quantitative R-S lemma} with noise process $\{U_t\}_{t \geq 1}$ and filtration $\{\mathcal{F}_t\}_{t \geq 0}$ such that conditions \eqref{cond-error} and \eqref{as-bound} hold. Assume additionally that
\begin{align*}
    \min_{1\leq i \leq m} \la \lb a_i + b_i \rb \wedge \lb c_i + d_i \rb \ra > 1.
\end{align*}
Then, for any $\delta \in (0,e^{-2})$, there exist constants $M, r > 0$ (independent of $t$) and a step-size sequence $\{\eta_t\}_{t \geq 1}$ such that
\begin{align} \label{conf-seq}
\mathbb{P}\left( \forall t \geq 0:\ L_t \leq M \cdot \frac{\log(\delta^{-1}) + \log\log(t+10)}{t+10} \right) \geq 1 - 2\delta,
\end{align}
provided that
\begin{align} \label{small enough}
\mathbb{P}(L_0 \leq r) \geq 1 - \delta.
\end{align}
Moreover, the step sizes satisfy $\eta_t \asymp 1/t$ as $t \to \infty$. 
\end{thm}

A few comments and explanations are in order.

{ \it On the choice of the step sizes $\eta_t$.} { In the statement of Theorem \ref{master} above, the term $\eta_t \asymp 1/t$ is interpreted in the sense that there exists a positive constant $C^{'}$ independent of $t$, but depends on $\delta$ and other parameters in \eqref{quantitative R-S lemma}, \eqref{cond-error} and \eqref{as-bound}, such that
\[
\frac{C^{'}}{t} \leq \eta_t \leq \frac{1}{C^{'}t}.
\]
We do not have an explcit construction for $\eta_t$ in the abstract settings of Theorem \ref{master}. However, explicit construction of the step sizes is possible if the bounds on the right-hand sides of \eqref{cond-error} and \eqref{as-bound} are simpler. Such settings are investigated in Section \ref{explicit} below.
}

{\it On the constants $M$ and $r$.} The constants $M$ and $r$ in the statement of Theorem \ref{master} is hard to construct explicitly based on the settings in \eqref{cond-error} and \eqref{as-bound}. However, the proof of Theorem \ref{master} reveals that for small $\delta$, $M$ and $r$ scale like $O \lb \mbox{polylog} \lb \delta^{-1} \rb   \rb $ and $O\lb 1/ \mbox{polylog} \lb  \delta^{-1}  \rb \rb $, respectively, with respect to the parameters in \eqref{cond-error} and \eqref{as-bound}. 

{\it On the form of our assumptions.}
Let us elaborate on the assumptions of Theorem~\ref{master}. The condition~\eqref{cond-error} requires that the conditional mean of the noise can be controlled as a function of the step sizes and the previous iterate. The $O(\eta_t^2)$ term in \eqref{cond-error} is standard and commonly appears in analyses of SGD. We allow smaller-order terms in \eqref{cond-error} to account for more intricate examples, such as Oja's algorithm, which will be discussed in Section~\ref{PCA}. Similarly, the condition~\eqref{as-bound} is used to control the magnitude of the noise based on the previous iteration.
Heuristically, for a sequence of admissible step sizes in the sense of \eqref{admissible}, the optimal convergence rate (ignoring logarithmic factors) is of order $O(1/t)$. Substituting $L_{t-1} = O(1/t)$ and $\eta_t = O(1/t)$ into the right-hand side of~\eqref{cond-error}, we obtain
\begin{align*}
    \left| \mathbb{E} \left[ U_t \mid \mathcal{F}_{t-1} \right] \right| \leq O(t^{-2}) + \sum_{k=1}^{m} O(t^{-1 - a_k - b_k}) = O(t^{-2}),
\end{align*}
provided that $a_k + b_k > 1$ for all $k \in \{1, \dots, m\}$. Thus, condition~\eqref{cond-error} ensures that the conditional expectation of the noise $U_t$ decays at rate $O(t^{-2})$. Similarly, condition~\eqref{as-bound} implies that the conditional variance of $U_t$ is of order $O(t^{-3})$ under the same assumption. Informally, this suggests that
\[
L_T \lesssim \sum_{t > T} \left| \mathbb{E} \left[ U_t \mid \mathcal{F}_{t-1} \right] \right| + \sqrt{ \sum_{t > T} \mathrm{Var} \left( U_t \mid \mathcal{F}_{t-1} \right) } \lesssim T^{-1}
\]
with high probability, for a single iterate. Of course, this is only a heuristic argument, and rigorous justification requires more substantial analysis. The $\log \log T$ factor arises from a  { peeling} argument, similar in spirit to the classical law of the iterated logarithm.

{\it Extension to random constants.}
In~\eqref{cond-error} and~\eqref{as-bound}, we assume that all parameters are constant. However, it is possible to extend Theorem~\ref{master} to a slightly more general setting where $A_i$ and $B_i$ are positive random variables, provided their conditional moment generating functions are sub-exponential with deterministic sub-exponential parameters. The step sizes stated in the statement of Theorem~\ref{master} are piecewise constant and asymptotically of order $O(1/t)$. However, in this abstract setting, we are unable to provide an explicit expression for the sequence $\{\eta_t; t \geq 1\}$. A special and important case where the step sizes can be made explicit is presented in Subsection~\ref{explicit}.


{\it Optimality of the bound.} The LIL-style rate $\log \log t / t$ in \eqref{conf-seq} is sharp and matches that of the classical law of the iterated logarithm for partial sum of i.i.d. random variables. The same rate has been obtained in \cite{howard2021time} in the context of empirical Bernstein's inequality. We also show that for stochastic approximation with squared-error loss, the same rate is optimal among all the choices of step sizes of order $\Omega \lb  1/t  \rb$; see Section \ref{Robbins-Monroe} for more details. 
There have been several attempts to obtain quantitative versions of Robbins-Siegmund's lemma; see \citep{neri2024quantitative, karandikar2024convergence, liu2022almost, sebbouh2021almost} and the references therein. However, none of these works successfully capture the $\log \log t$ term in the bound. To fully capture this subtle logarithmic term, we employ a novel, sharp error bound in short-time intervals and a ``stitching'' argument. 

{\it Relationship with supermartingales.}
A natural but naive approach would be to construct a non-negative supermartingale directly from the recursive estimate \eqref{quantitative R-S lemma}. To the best of our knowledge, such a construction is not applicable to the settings in Theorem \ref{master} due to the highly nonlinear nature of \eqref{cond-error} and \eqref{as-bound}. In the special case where the noises $U_t$'s are deterministic and of order $O(\eta_t^2)$, it is indeed possible to construct a non-negative supermartingale directly; see \cite{gladyshev1965stochastic, lai2006extended} and the references therein. However, this construction breaks down when the noises are random and more complicated, as it involves an infinite product starting at time $n$ and running to infinity, which in turn leads to measurability issues.

{\it On the necessity of condition \eqref{small enough}.}
Interestingly, though the conditions stated in Theorem \ref{master} guarantee that $L_t$ has an almost sure limit when $a_i$'s are greater than $1$,  they do {\it not} guarantee the convergence of $L_t$ to $0$ almost surely. For example, consider the recursion
\begin{align} \label{counter example}
L_{t} \leq (1-\eta_t)L_{t-1} + 2\eta_t L_{t-1}^2 + \eta_t^2.
\end{align}
It is easy to see that the recursion above is of the form \eqref{quantitative R-S lemma} and satisfies \eqref{cond-error}, \eqref{as-bound} with $C_2 \equiv 1$, $C_3 \equiv 0$, $\lb A_1, a_1, b_1 \rb \equiv (2,0,2)$, $\lb B_1, c_1, d_1 \rb \equiv \lb 2, 1/2, 2 \rb $ and $\lb  B_2, c_2, d_2 \rb \equiv (1,3/2,0)$. However, if one takes a Bernoulli random variable $Y$ such that $P \lb Y=1 \rb=1/10$, then the process $L_t \equiv Y$ clearly satisfies \eqref{counter example} but
\[
\mb{P} \lb \lim _{t \to \infty} L_t = 0 \rb = 9/10.
\]
The problem with this counter example is that $L_0$ can be large with non-negligible probability. Therefore, to rule out such counter examples, one should only expect a time-uniform bound of the form \eqref{conf-seq} under \eqref{small enough}, where $r$ is sufficiently small and can be determined by the parameters in \eqref{cond-error} and \eqref{as-bound}. As a result, Theorem \ref{master} captures the stable dynamic of a general class of iterated algorithms: once the output of the algorithm gets close to the minimizer, it will stay close to the minimizer within a neighborhood of optimal width $\log \log t/ t$.

\subsection{An explicit step size construction} \label{explicit}

Although Theorem~\ref{master} provides a general affirmative result regarding the time-uniform bound with optimal width, the parameters involved in the construction cannot be made explicit as a function of time~$t$. To address this issue, we provide an explicit construction in a simpler setting that covers all the applications considered in the subsequent sections.

{ For simplicity, let us now focus on a simplified form of \eqref{cond-error} and \eqref{as-bound}. Suppose the process $\la L_t ; t \geq 0 \ra$ satisfies \eqref{quantitative R-S lemma} such that the noise process $U_t$ satisfying
 \begin{align}
      &  \Big| \mb{E} \lb U_t | \mathcal{F}_{t-1}  \rb \Big| \leq  C_2 \cdot \eta_t^2;    \label{cond-error 2} \\
       & |U_t| \leq C_3 \eta_t \sqrt{L_{t-1}} +  C_2 \cdot \eta_t^{2}.  \label{as-bound 2}
    \end{align}
Conditions \eqref{cond-error 2} and \eqref{as-bound 2} correspond to the most common setting of SGD under a strongly convex function. This context will be further analyzed in Section \ref{sgd} below. We first state a maximal inequality with explicit parameters.

As we will see shortly in Section \ref{sgd}, not only the maximal inequality in Proposition \ref{pre-planned 1} can be used to derive time-uniform bounds, it can also improve the existing last iteration bound for the classical SGD. Let $\la L_t ; t \geq 0 \ra$ satisfies \eqref{quantitative R-S lemma} such that the noise process $U_t$ satisfying \eqref{cond-error 2} and \eqref{as-bound 2}, then
}


    \begin{prop}[maximal inequality] \label{pre-planned 1}
 Assume that $\eta_t = 2/ \lb C_1 \lb t+L \rb \rb$ for all $t \in [t_0+1,t_1]$ and some positive integer $L\geq 3$. Suppose $A \subset \la L_{t_0} \leq a \ra$ such that $\mb{P} \lb A \rb >0$. Then, for any $\delta>0$, one has
    \[
    \mb{P} \lb \exists t \in [t_0+1, t_1]: L_t \geq \frac{M \lb t_1-t_0 \rb \log \lb \delta^{-1}\rb}{\lb t+3 \rb^2} \Big| A \rb \leq \frac{\delta}{\mb{P} \lb A \rb}
    \]
    where 
    \[
    M = \frac{31.5 \times (L-1)}{L}\max \la \frac{ aL(L-1)}{ \log \lb \delta^{-1} \rb(t_1-t_0)}; \frac{C_2}{C_1^2 \log \lb \delta^{-1} \rb}; \frac{C_2}{C_1^2\sqrt{ \log \lb \delta^{-1} \rb }}; \frac{C_3^2}{C_1^2}  \ra.
    \]
\end{prop}

Proposition \ref{pre-planned 1} gives a maximal inequality for the sequence $\la L_t \ra $ over the time interval $[t_0,t_1]$ when the step sizes are set to be  $\eta_t = 2/ \lb C_1 \lb t+L \rb \rb$. If one sets $t_0 \equiv 0$ and suppose that $L_{t_0} \leq a$ almost surely, then for all $\delta$ sufficiently small, we have
\begin{align} \label{last iterate}
L_{t_1} \leq 31.5 \times \frac{L-1}{L} \cdot \frac{ C_3^2 \log \lb \delta^{-1} \rb }{C_1^2 t_1}
\end{align}
with probability at least $1-\delta$.

In the context of SGD, \eqref{last iterate} implies that the last iterate is of order $1/t$, which is an improvement of some results in literature. This will be discussed in more details in Section \ref{sgd}. Proposition \ref{pre-planned 1} allows us obtain the time-uniform bound in the same settings.

\begin{thm} \label{conf}
    Consider the settings in Proposition \ref{pre-planned 1}. Then,
    \[
        \mb{P} \lb \exists t \geq 0: L_t \geq 31.5 \times K \max \la \frac{aL}{\log \lb \delta^{-1} \rb}; \frac{C_2}{C_1^2}; \frac{C_3^2}{C_1^2} \ra \frac{\log \lb \delta^{-1} \rb + 2\log \log (t+9)}{t+L} \rb \geq 1-2\delta
    \]
    if
    \[
    \mb{P}\lb    L_0 \leq a   \rb \geq 1- \delta,
    \]
    and
    \begin{align} \label{K}
    K:= \max \la L-2; 32 \ra  \cdot \lb \mathbf{1}_{\la L \geq 32 \ra} + 32 \cdot \mathbf{1}_{\la L\leq 31 \ra} \rb.
    \end{align}
\end{thm}

Theorem~\ref{conf} provides a time-uniform bound for any process $\{L_t\}_{t \geq 0}$ satisfying conditions~\eqref{quantitative R-S lemma}, \eqref{cond-error 2}, and \eqref{as-bound 2}, given an initial upper bound $a$ on $L_0$ that holds with probability at least $1 - \delta$. In contrast to Theorem~\ref{master}, we do not require $a$ to be sufficiently small, and both the constants and step sizes in the construction are made explicit. This should not be surprising: under the conditions \eqref{quantitative R-S lemma},  \eqref{cond-error 2}, and \eqref{as-bound 2}, one can show that $L_t$ converges to $0$ almost surely for any admissible sequence of step sizes $\la \eta_t; t \geq 1 \ra$ in the sense of \eqref{admissible}. The requirement that $a$ be sufficiently small, as in Theorem~\ref{master}, arises only in the \textit{inhomogeneous} setting, where the right-hand side contains terms of order $\eta_t^{k}$ for some $k < 2$.

\section{Applications} 

As anticipated, we will exemplify the general time-uniform concentration bounds of Theorem~\ref{master} in three important applications.

\subsection{SGD algorithms}
\label{sgd}

We consider standard settings for the analysis of SGD algorithms, which we briefly review; see, e.g., \cite{nemirovski2009robust,rakhlin2011making} (see also the references therein). Throughout, $\mathcal{X}$ denotes a compact, convex set in $\mb{R}^d$ with non-empty interior and $F: \mc{X} \to \mb{R}$ a real valued function on it. We will consider several assumptions on $F$. 
\begin{itemize}
    \item {\it Strong convexity.} $F$ is said to be $\lambda$-strongly convex ($\lambda$-SC) if for all $\bm{x},\bm{y}\in \mc{X}$
    \begin{equation} \label{strong convex}
  \tag{$\lambda$-SC}
          F(\bm{y}) \geq F(\bm{x}) + \langle \bm{y}-\bm{x}, \nabla F( \bm{x}) \rangle + \frac{\lambda}{2} \| \bm{y}-\bm{x}\|^2.
\end{equation}

   \item {\it Smoothness.} $F$ is said to be $\mu$-smooth if for all $\bm{x},\bm{y}\in \mc{X}$
     \begin{equation} \label{smoothness}
     \tag{$\mu$-smooth}
    F(\bm{y}) \leq F(\bm{x}) + \langle \bm{y}-\bm{x}, \nabla F(\bm{x}) \rangle + \frac{\mu}{2} \| \bm{y}-\bm{x}\|^2.
     \end{equation}

 \item {\it Polyak-Lojasiewicz condition.} $F$ is said to satisfy the Polyak-Lojasiewicz (PL) condition with parameter $\tau>0$ if for all $\bm{x},\bm{y}\in \mc{X}$
   \begin{equation} \label{PL condition}
       \tag{P-L} 
       \| \nabla F(\bm{x}) \|^2 \geq \tau \lb  F(\bm{x}) - F(\bm{x}^{*}) \rb 
   \end{equation}
    where $\bm{x}^{*} \in \mbox{int} \lb \mc{X} \rb$ is the global minimum. 
\end{itemize}

We are interested in solving the optimization problem
\begin{align} \label{optimization}
    \min_{\bm{x} \in \mc{X}}  F\lb \bm{x} \rb 
\end{align}
where  $ F\lb \bm{x} \rb: = \mb{E} \lb f(\bm{x}, \bm{\xi}) \rb $,  $\bm{\xi}$ is a random vector from a distribution $P_{\bm{\xi}}$ supported in $\mb{R}^{d_1}$ and $ f(\bm{x}, \bm{\xi})$ is a convex function on $\mathcal{X} \times \mathbb{R}^{d_1}$, differentiable in $\bm{x}$ for almost surely all realizations of $\bm{\xi}$. Under this assumption, it is easy to check that $F$ is a convex differentiable function on $\mc{X}$. 

In the SGD framework, one has access to a series of unbiased estimators of the true gradient of $F$ by sampling $\bm{\xi}$ repeatedly. Denote by $\Pi_{\mc{X}}: \mb{R}^d \to \mc{X}$ the projection map
\[
\Pi_{\mc{X}} \lb \bm{y} \rb:= \argmin_{\bm{x} \in \mc{X}} \| \bm{x} - \bm{y} \|.
\]
It is well-known that $\Pi_{\mc{X}}$ is well-defined and is a $1$-Lipchitz map when $\mc{X}$ is a compact convex set. The projected SGD algorithm for solving the problem \eqref{optimization} is defined iteratively as
\begin{align} \label{projected-sgd}
\bm{x}_{t} := \Pi_{\mc{X}} \Big( \bm{x}_{t-1} - \eta_{t}  g \lb \bm{x}_{t-1}, \bm{\xi}_{t} \rb  \Big)
\end{align}
where $\la \bm{\xi}_t; t \geq 1 \ra$ is a sequence of i.i.d. samples from $P_{\bm{\xi}}$, $\bm{x}_0 \in \mc{X}$ is a starting point and $g \in \mathbb{R}^d$ and noise vector satisfying
\[
\mb{E} \lb g \lb \bm{x}_{t-1}, \bm{\xi}_{t} \rb \Big| \bm{\xi}_1, \bm{\xi}_2,\dots, \bm{\xi}_{t-1} \rb = \nabla F(\bm{x}_{t-1}),
\]
and 
\[
\max \la \|g - \nabla F \|; \|g \| \ra  \leq B
\]
almost surely, for some positive constant $B>0$.
\subsection{Time-uniform bounds under (\ref{strong convex})} Suppose the function $F$ satisfies \eqref{strong convex}. Letting $L^{SGD}_t =\| x_{t} - x^{*}\|^2$, we immediately have the following recursion.
\begin{prop} \label{SC recursive}
    For any $t \geq 1$, we have 
    \begin{align} \label{sgd-recursive1}
           L^{SGD}_t \leq \lb 1 - 2\lambda \eta_t \rb L^{SGD}_{t-1} +2\eta_t Y^{SGD}_t + B^2 \eta_t^2
    \end{align}
    where 
     $\mb{E} \lb  Y^{SGD}_t | \mc{F}_{t-1} \rb =0 $ and $|Y_t| \leq B \sqrt{L^{SGD}_{t-1}}$.
\end{prop}
From Proposition~\ref{SC recursive}, we observe that the process $\la L_t; t \geq 0 \ra$ satisfies a recursion of the form~\eqref{quantitative R-S lemma}, with a noise process that satisfies conditions~\eqref{cond-error 2} and~\eqref{as-bound 2}, along with the corresponding constants
\begin{align*}
    C_1 \equiv 2\lambda; \ C_2 \equiv  B^2; \ C_3 \equiv 2B.
\end{align*}

{ An application of Theorem \ref{conf} immediately yields time-uniform concentration bound for SGD.}
\begin{corollary} \label{sgd-conf}
Set the step sizes $\eta_t = 1/ \lb  \lambda( t+32) \rb$, then 
    \[
    \mathbb{P} \left( \forall t \geq 1: L^{SGD}_t \leq 1008 \cdot \frac{B^2}{\lambda^2} \cdot \frac{\log(\delta^{-1}) + 2\log \log (t+9)}{t+32} \right) \geq 1-\delta
    \]
    for any $\delta >0$.
\end{corollary}
Corollary~\ref{sgd-conf} is a direct application of Theorem~\ref{conf} and the fact that 
\[
L^{SGD}_t \leq \frac{B^2}{\lambda^2}
\]
almost surely, for all $t\geq 1$. This simple fact follows from the boundedness assumption and the strong convexity assumption; see Lemma 5 in \cite{rakhlin2011making} for a proof. We can see that Corollary~\ref{sgd-conf}  provides a time-uniform bound with optimal width with a larger constant than the convergence rate for the last iterate and an extra $\log \log$ factor, which is optimal and can not be removed.

It is helpful to compare the previous bound with an analogous one given in Proposition 1 of \cite{rakhlin2011making}, which, under the same settings,  states that for a fixed time $T \geq 3$, and setting $\eta_t = \frac{1}{\lambda t}$,
\begin{align} \label{sgd rakhlin}
\mb{P} \lb \forall 1 \leq t \leq T: L^{SGD}_t \leq 624 \cdot \frac{B^2}{\lambda^2} \frac{ \lb  \log \lb \delta^{-1}\rb + \log \log T  \rb}{ t} \rb \geq 1-\delta.
\end{align}
The previous bound is, of course, not anytime-valid, as it requires the specification of a terminal time $T$. In contrast, \ref{sgd-conf} provides time-uniform guarantees while exhibiting the same dependence on the model-related parameters $\frac{B^2}{\lambda^2}$.

{
Not only our method gives any-time concentration bound, it can also yield insight on the last iteration of SGD. To be more precise, we have
\begin{corollary} \label{bound for last iterate}
 Set $\eta_t = 1/(\lambda(t+3))$. Then, for any $\delta>0$ and $t\geq 1$,
\[
\mb{P} \lb L^{SGD}_t \leq  \frac{21B^2}{\lambda^2}  \cdot \frac{\log \lb \delta^{-1} \rb}{t+3} \rb \geq 1-\delta.
\]
\end{corollary}
In comparison with~\eqref{sgd rakhlin}, our result in Corollary~\ref{bound for last iterate} improves the convergence rate of the last iterate by a factor of $\log \log T$. In the next subsection, we will consider a non-convex setting involving the popular P-L condition. 
}


\color{black}
{
 Let us end this subsection with an example regarding least square.
\subsubsection*{Sequential ridge regression via SGD}
To illustrate the effectiveness and generality of our time-uniform bounds, we show how Corollary~\ref{sgd-conf} immediately yields  anytime concentration bounds of optimal width based on the SGD iterates for the problem of estimating in a sequential manner the parameters of a standard regression linear model. 
In detail, we assume that we observe sequentially a stream of i.i.d. of response/covariates pairs $(\{ y_t, \bm{x}_t)\}_{t=1,2,\ldots} $ obeying the linear regression model 
\begin{equation}
\label{eq:linear.model}
y_t = \langle \bm{\theta}^{*}, \bm{x}_t \rangle + \xi_t.
\end{equation}
Above, the unknown parameter $\bm{\theta}^*$ belongs to a known convex set  $\Omega \subset \mb{R}^d$ such that $\mathrm{diam}(\Omega) = D$ and the random covariates and noise are bounded: for all $t$,
\[
\max\!\big\{ \| \bm{x}_t \|,\, | \xi_t | \big\} \le B
\quad \text{almost surely.}
\]
We are interested in estimating $\bm{\theta}^{*}$ using the sequential ridge regression estimator obtained as the output of the SGD iterates applied to the minimization problem 
\begin{align*}
     \min_{\bm{\theta} \in \Omega } \la \frac{1}{2} \mathbb{E}\!\left[ \big( y - \langle \bm{\theta}, \bm{x} \rangle \big)^2 \right] + \frac{\lambda}{2} \| \bm{\theta} \|^2 \ra,
\end{align*}
where $\lambda \geq 0$ is a fixed regularization parameter and the pair $(y,\bm{m})$ is a draw from the data generating distribution obeying the linear model \eqref{eq:linear.model}. The case $\lambda = 0$ corresponds to the least squares minimization.

In this setting, the SGD algorithm initialized at an arbitrary $\hat{\bm{\theta}}_0$ yields the update rule
\begin{equation}\label{eq:SGD.ridge}
    \hat{\bm{\theta}}_t
= \mbox{Proj}_{\Omega} \lb  \hat{\bm{\theta}}_{t-1}
- \eta_t \bm{x}_t \big( \bm{x}_t^\top \hat{\bm{\theta}}_{t-1} - y_t  \big) + \lambda \hat{\bm{\theta}}_{t-1} \rb,
\end{equation}
for a given sequence of step sizes $\{\eta_t\}_{t \ge 1}$. 
Then, assuming that 
\[
\lambda_{\min} := \lambda_{\min}\!\big( \mathbb{E}[\bm{x}_1 \bm{x}_1^\top] \big) > 0,
\]
with the choice of step sizes $\eta_t = 2 / \big( \lambda_{\min}(t + 32) \big)$, Corollary~\ref{sgd-conf} implies that, with probability at least $1-\delta$  uniformly over all $t \geq 1$,
\begin{align} \label{Ridge any time}
    \| \hat{\bm{\theta}}_t - \bm{\theta}^* \|^2
    \le \frac{\lambda^2 \| \bm{\theta}^* \|^2}{ \lambda_{\rm min} ^2} +
    1008 \cdot \frac{B_1^2}{\lambda_{\min}^2}
    \cdot
    \frac{\log(\delta^{-1}) + 2\log \log (t + 9)}{t + 32}, 
\end{align}
where
\[
B_1:= B^2D+B^2 +\lambda D + \lambda \| \bm{\theta}^* \|.
\]
Note that when the penalty parameter $\lambda$ is zero, the bound does not depend on $\| \bm{\theta}^* \|$. This is not unexpected, since the SGD procedures in this case converges to the ordinary least squares estimator, which is unbiased in well-specified settings.

To the best of our knowledge, the above bound is new in the SGD literature, both because of its time-uniform form and because of its optimal rate of decay in $t$.
It is illustrative to compare the above bound with existing anytime concentration bound for the online ridge estimator obtained by \cite{abbasi2011improved,abbasi2011online}; see also \cite{lattimore2020bandit,catoni2017dimension,martinez2025vector} and the references therein.
Let
\begin{align*}
    \hat{\bm{\theta}}_t^{\mathrm{Ridge}} &:= \big( \bm{x}_{1:t}^\top \bm{x}_{1:t} + \lambda \bm{I} \big)^{-1} \bm{x}_{1:t}^\top \bm{y}_{1:t}, \\
    \bar{\bm{V}}_t &:= \lambda \bm{I} + \sum_{s=1}^t \bm{x}_s \bm{x}_s^\top, \\
    \| \hat{\bm{\theta}}_t^{\mathrm{Ridge}} - \bm{\theta}^* \|_{\bar{\bm{V}}_t}
    &:= \sqrt{ \big\langle \hat{\bm{\theta}}_t^{\mathrm{Ridge}} - \bm{\theta}^*,\, 
    \bar{\bm{V}}_t \big( \hat{\bm{\theta}}_t^{\mathrm{Ridge}} - \bm{\theta}^* \big) \big\rangle }.
\end{align*}
Here, $\bm{x}_{1:t}$ denotes the matrix with rows $\bm{x}_1^\top, \dots, \bm{x}_t^\top$.
Then, Theorem~2 in \cite{abbasi2011improved} establishes that
\begin{align} \label{Ridge anytime 2}
\mathbb{P}\!\left(
    \forall t \ge 1:\;
    \| \hat{\bm{\theta}}_t^{\mathrm{Ridge}} - \bm{\theta}^* \|_{\bar{\bm{V}}_t}
    \le
    B \sqrt{d \log \!\left( \frac{1 + tB^2 / \lambda}{\delta} \right)}
    + \sqrt{\lambda}\, \| \bm{\theta}^* \|
\right)
\ge 1 - \delta.
\end{align}
Though both concentration bounds \eqref{Ridge any time} and \eqref{Ridge anytime 2} are anytime-valid, they differ in several important aspects.
First, the bound by \cite{abbasi2011improved} is with respect to the Malahanobis distance using  the sample covariance matrix, as it is based on a self-normalization approach \cite{de2007pseudo,de2009self} that is specifically well-suited to the linear regression settings and remain valid  under minimal assumptions on the data generating distribution; in particular, the data need not be i.i.d. and can be chosen arbitrarily in a data-adaptive manner. 
In contrast, our bound was obtained as corollary to  a general time-uniform concentration bound for almost-super martingale-type processes (namely, Theorem \ref{conf}) that was not specifically designed for linear regression. 

The conditions we impose for our bound, namely that  the errors and the covariates are
bounded and drawn in an i.i.d. manner,  are of course restrictive compared to ones assumed by \cite{abbasi2011improved}. However, it is worth mentioning that they are in fact standard conditions in the literature on SGD,  a methodology that should not be expected to work well in data-adaptive settings where the next observation is chosen in a predictable manner\footnote{Indeed, consistency of the SGD iterates is predicated on the assumption that the gradient of the target function can be estimated unbiasedly using a new observation, which rules out virtually any non-i.i.d. data streams.}

A second important difference is that our bound \eqref{Ridge any time} is concerned with a sequential, SGD-based estimator $\hat{\bm{\theta}}_t$ given in \eqref{eq:SGD.ridge}, a one-pass procedure with a computational complexity $O(d)$ per iteration. In contrast, the more expensive online ridge estimator in \eqref{Ridge anytime 2} is based on the inverse of the regularized Gram matrix $ \bar{\bm{V}}_t $ at every $t$, which, using the Sherman-Morrison formula,  can computed as a rank-one-update with $O(d^2)$ operations at each iteration. 

 Finally, in terms of scaling in $t$, the (squared) bound in \eqref{Ridge anytime 2} is of order $(\log t) / t$, while our bound in \eqref{Ridge any time} is of order $(\log \log t) / t$, thereby recovering the optimal logarithmic dependence, as suggested by Theorem~\ref{lower bound} below.    Recently, in the context of deriving anytime valid concentration inequalities for self-normalized multivariate processes,  \cite{whitehouse2023time} have sharpened the results in \cite{abbasi2011improved,abbasi2011online} and recovered the  optimal scaling $(\log \log \lb t)/t \rb$ for i.i.d. data; see Theorem 5.2 of  therein. Thus, our result implies that, in the i.i.d settings in which SGD is usually deployed, the one-pass SGD-based estimator \eqref{eq:SGD.ridge} achieves the same rate as the estimator computed using the full covariance matrix, while being more efficient in terms of computational cost.

}





\subsection{Uniform error bound under (\ref{smoothness}) and (\ref{PL condition})}

We now turn to the analysis of SGD in non-convex problems, where we assume that (\ref{smoothness}) and (\ref{PL condition}) hold. In this setting, we are interested in estimating the minimal value $F(\bm{x}^{*})$ of the objective function at any minimizer $\bm{x}^{*}$.

Given the sequence $\{ \bm{x}_t,t \geq 0\}$ of projected SGD iterations in \eqref{projected-sgd} we let 
\[
L^{PL}_t:= F(\bm{x}_t) - F(\bm{x}^{*}).
\]
Similar to Proposition \ref{SC recursive}, the corresponding recursion in this case is given next.

\begin{prop} \label{P-L recursive}
    For any $t \geq 1$, we have 
    \begin{align} \label{sgd-recursive2}
            L^{PL}_{t+1} \leq (1-\tau \eta_t) L^{PL}_t + \eta_t Y^{PL}_t + 2 \mu B^2 \cdot \eta_t^2
    \end{align}
    where 
    \[
    Y^{PL}_t := \langle \nabla F(\bm{x}_t), \nabla F(\bm{x}_t)-\hat{g}(x_t) \rangle.
    \]
    Moreover, $\mb{E} \lb  Y^{PL}_t | \mc{F}_{t-1} \rb =0 $ and $|Y^{PL}_t| \leq B\sqrt{\mu} \cdot \sqrt{L^{PL}_t}$.
\end{prop}
Then, Proposition \ref{pre-planned 1} and  Theorem \ref{conf} yield the following fixed and anytime-valid concentration bounds.


{
\begin{corollary} \label{any-time for PL}
Recall $\tau$ in \eqref{PL condition}. Consider the optimization problem \eqref{optimization} with the iterates \eqref{projected-sgd}. 

(i) Set  the step sizes $\eta_t=2/\lb \tau(t+3) \rb$. For any $\delta \in (0,e^{-4})$ and $t\geq 3/ \log \lb \delta^{-1} \rb$, 
\[
\mb{P} \lb  F \lb \bm{x}_t \rb - F \lb \bm{x}^* \rb \geq \frac{21\mu B^2}{\tau^2} \cdot \frac{\log \lb \delta^{-1} \rb}{t+3}  \rb \leq \delta.
\]

(ii) Set  the step sizes $\eta_t=2/\lb \tau(t+32)\rb$. Then, uniformly over all time $t \geq 0$, 
\[
 F \lb \bm{x}_t \rb - F \lb \bm{x}^* \rb \leq 1008 \cdot \max \la \frac{128B^2}{\tau \log \lb \delta^{-1} \rb} ; \frac{2B^2\mu}{\tau^2} \ra \cdot \frac{ \log \lb  \delta^{-1} \rb + 2\log \log (t+9)}{t+32}
\]
with probability at least $1-\delta$. 
\end{corollary}

Corollary \ref{any-time for PL} is an immediate application of Proposition \ref{pre-planned 1} and \ref{conf} and the fact that for all $t \geq 0$,
\[
 F \lb \bm{x}_t \rb - F \lb \bm{x}^* \rb \leq \frac{1}{\tau} \| \nabla F \lb \bm{x}_t \rb \|^2 \leq \frac{4B^2}{\tau}
\]
which follows from \eqref{PL condition}.

To the best of our knowledge, the time-uniform bound for the error under \eqref{PL condition} and \eqref{smoothness} in Corollary \ref{any-time for PL} is the first in the literature. Regarding the convergence of the last iterate, a similar result but without explicit constants was recently obtained by \cite{madden2024high}.

}



\subsection{Oja's algorithm for PCA} \label{PCA}

Principal Component Analysis (PCA) is one of the most fundamental problems in multivariate analysis, where the goal is to estimate a low-dimensional subspace from observed data. 
In this section, we focus on the sequential setting, in which estimation is performed as data points arrive one by one. 
Accordingly, the results are presented in the fixed-dimensional setting. 
However, as we shall see, the effect of dimensionality appears only through a logarithmic factor, making the proposed method well-suited for large datasets.

Let $\bm{X}_1, \bm{X}_2, \dots, \bm{X}_n \in \mathbb{R}^p$ be i.i.d. centered random vectors with covariance matrix $\bm{\Sigma}$ such that:  

\begin{itemize}
    \item $\| \bm{X}_i \| \leq B$ almost surely.  
    \item $\lambda_1 > \lambda_2$, where $\{ \lambda_i \}_{i=1}^{p}$ are the eigenvalues of $\bm{\Sigma}$ sorted in descending order.  
\end{itemize}  

The goal of PCA is to estimate the principal component, i.e., the eigenvector associated with the largest eigenvalue $\lambda_1$. In optimization terms, this corresponds to solving  
\begin{align} \label{principal vector}
    \bm{v}_1 := \argmin_{\| \bm{v} \|=1} -\bm{v}^{\top} \bm{\Sigma} \bm{v}.
\end{align}
In statistical literature, a natural approach is to use the leading eigenvector of the empirical covariance matrix,  
\[
\hat{\bm{\Sigma}}_n := \frac{1}{n} \sum_{k=1}^{n} \bm{X}_k \bm{X}_k^{\top}.
\]  
However, this method is computationally expensive for large datasets and it requires storing all data points, which is often impractical due to memory constraints.  

Oja's algorithm, initially proposed in \cite{oja1982simplified,oja1985stochastic}, is an efficient method for solving \eqref{principal vector} that requires only \( O(p) \) memory and has a runtime of \( O(np) \). Oja's idea was to recast the problem \eqref{principal vector} as an optimization problem and apply stochastic gradient descent (SGD) to solve it iteratively. Interestingly, the objective function is concave, and hence highly non-convex. Nevertheless, the optimization landscape is simple and free of spurious local minima, making this one of the simplest non-convex problems that can be solved efficiently using SGD, both in theory and in practice.

An incomplete list of early results includes \cite{jain2016streaming,de2015global,hardt2014noisy,balsubramani2013fast,shamir2016convergence,mitliagkas2013memory,balcan2016improved}. More recent works, which establish optimal rates and extensions to various settings, can be found in \cite{huang2021streaming,lunde2021bootstrapping,kumar2024oja,allen2017first,lowprecPCA2025,Kumar25entry}; see also the references therein. However, to the best of our knowledge, there are no existing results on deriving any-time bounds or confidence sequences for the sequence of estimators $\la \hat{\bm{v}}_t; t \geq 1 \ra$ defined above.
A natural approach to obtaining such bounds is to apply matrix concentration inequalities. However, such methods often yield suboptimal convergence rates in $t$. This is not surprising: achieving the optimal convergence rate for Oja’s algorithm requires a two-stage step-size schedule (see \cite{huang2021streaming} for details): a chaotic phase with constant step sizes, during which the algorithm explores the sphere to locate the region containing the eigenvector, followed by a stable phase in which the algorithm achieves the $O(t^{-1})$ convergence rate. Consequently, concentration inequalities applied to a single step-size schedule of the form $\eta_t = A/(t+B)$ are unlikely to yield the optimal rate.

To describe Oja's algorithm, we first reformulate \eqref{principal vector} as a stochastic optimization problem, where the population covariance matrix is replaced by its empirical counterpart:
\begin{align} \label{principal vector oja}
    \bm{\hat{v}}_n := \argmin_{\| \bm{v} \|=1} -\bm{v}^{\top} \hat{\bm{\Sigma}}_n \bm{v}.
\end{align}
Note that the optimization problem \eqref{principal vector oja} is \textit{non-convex}. However, the landscape of the objective function is relatively simple: its population version \eqref{principal vector} has no spurious local minima. In fact, this is one of the simplest non-convex problems that can be solved efficiently. In what follows, we describe the original Oja's algorithm and a common variant, known as the Krasulina-Oja algorithm.

 \textbf{(i) Oja's algorithm.} Initialize a vector \( \hat{\bm{v}}_0 \sim \text{Uni} \big( \mathbb{S}^{p-1} \big) \). Let \( \{ \eta_t \}_{t\geq 1} \) be a sequence of step sizes. The update rule of Oja's algorithm is given by:
\begin{equation} \label{Oja}
  \tag{Oja}
      \hat{\bm{v}}_{t} := \big( \bm{I}_p + \eta_{t}\bm{X}_{t} \bm{X}_{t}^{\top} \big) \hat{\bm{v}}_{t-1}, \quad \hat{\bm{v}}_{t} \leftarrow \frac{\hat{\bm{v}}_{t}}{\|\hat{\bm{v}}_{t}\|}.
\end{equation}
Under the standard assumption \eqref{admissible} on the step sizes, it can be shown that Oja's algorithm converges almost surely to either \( \bm{v}_1 \) or \( -\bm{v}_1 \). Moreover, under \eqref{admissible}, the normalization step can be omitted in practice, and it suffices to normalize only the output in the last iteration.

 \textbf{ (ii) Krasulina-Oja's algorithm.} Initialize a vector \( \hat{\bm{v}}_0 \sim \text{Uni} \big( \mathbb{S}^{p-1} \big) \). Let \( \{ \eta_t \}_{t\geq 1} \) be a sequence of step sizes. The update rule of Krasulina-Oja's algorithm is given by:
\begin{equation} \label{krasulina}
    \tag{Kra-Oja}
        \hat{\bm{v}}_{t} = \hat{\bm{v}}_{t-1} + \eta_t \underbrace{y_t \left[ \bm{X}_t - \frac{y_t}{\| \hat{\bm{v}}_{t-1}\|^2}\cdot \hat{\bm{v}}_{t-1}  \right]}_{\bm{z}_t}
\end{equation}
where 
\begin{align} \label{yt}
    y_t:= \bm{X}_t^\top \hat{\bm{v}}_{t-1}.
\end{align}
The assumption that the $\bm{X}_i$’s are centered is not essential and can be removed with a simple modification. For example, in \eqref{Oja}, instead of updating with a single data point at each iteration, one may use two data points and modify the update rule as
\[
    \hat{\bm{v}}_{t} 
    := \left[ \bm{I}_p + \frac{\eta_{t}}{\sqrt{2}} \big( \bm{A}_t - \bm{B}_t \big) \big( \bm{A}_t - \bm{B}_t \big)^{\top} \right] \hat{\bm{v}}_{t-1},
\]
where $\bm{A}$ and $\bm{B}$ are two independent copies of $\bm{X}_1$ and are independent of the history up to time $t-1$.
This modification ensures that the update is centered, regardless of the original mean. 
For the sake of simplicity, we will assume throughout the rest of the paper that the data are centered.

The difference between \eqref{krasulina} and \eqref{Oja} is the second order correction term in $\bm{z}_t$. One advantage of such a term is that the update is orthogonal to the previous iteration. It is also known that under the admissible condition \eqref{admissible}, the algorithm \eqref{krasulina} output a sequence $\la \hat{\bm{v}}_t; t \geq 1 \ra$ such that 
\[
\lim_{t \to \infty} \sin^2 \lb \frac{\hat{\bm{v}}_{t}}{\|\hat{\bm{v}}_{t}\|}, \bm{v}_1 \rb =0.
\]
Here the sine-squared error loss $\sin^2(\bm{x}, \bm{y})$ is defined as
\[
\sin^2 (\bm{x}, \bm{y}) = 1-\frac{\langle \bm{x}, \bm{y} \rangle}{\| \bm{x} \| \cdot \| \bm{y} \|}
\]
for $\bm{x}, \bm{y} \in \mb{R}^p$.

The convergence analysis of \eqref{krasulina} and \eqref{Oja} has been extensively studied in the literature. Convergence results for the principal eigenvector can be found in \cite{de2015global,shamir2016convergence,shamir2016fast}, along with additional references therein. In contrast, results for $k$-PCA are relatively scarce; see \cite{allen2017first,huang2021streaming} and the references therein. Recent developments in Oja's algorithm include a bootstrap approximation for the sine-squared error \cite{lunde2021bootstrapping} and a variant designed specifically for sparse PCA settings \cite{kumar2024oja}.

Despite the extensive literature on the convergence analysis of the final iteration of Oja's algorithm, little is known about proving time-uniform bounds for its iterates. There are two fundamental difficulties in establishing such a result. First, the dynamics of the algorithm are unstable during the initial iterations, providing little useful information about the principal eigenvector. Second, the outputs in \eqref{Oja} and \eqref{krasulina} involve products of independent random matrices, which lack an immediately exploitable martingale structure. Indeed, while a {\it fixed} number of products of independent random matrices does exhibit concentration properties, as recently demonstrated by \cite{matrix.product}, it is not obvious how to extend such sharp bounds to a random number of products. 

It is worth noting that the $(1-\delta)$ time-uniform bound for Oja's algorithm cannot be obtained by simply applying a union bound over fixed-time convergence results. The difficulty arises because the step sizes themselves depend on $\delta$, which alters the solution path if one attempts to apply last-iterate convergence results sequentially over a collection of probability levels.
 In contrast, the step-size construction in Section~\ref{sgd}, corresponding to the SGD case, is independent of $\delta$.

The first difficulty can be resolved by using the smoothing technique in \cite{huang2021streaming} to show that after an initial uncertainty phase of order $T_0: =  \Omega \lb B \delta^{-2} (\lambda_1- \lambda_2)^{-2} \rb$ steps, the output will be within a ball of radius $(T-T_0)^{-1/2}$ of the principal eigenvector. However, the second difficulty appears to be intractable using existing techniques, as constructing a non-negative supermartingale solely based on the matrix product remains difficult.

We will now demonstrate how Theorem \ref{conf} can be used to resolve the second difficulty. Let us first start with an useful observation: without loss of generality, one can assume $\bm{\Sigma}$ to be a diagonal matrix, i.e.
\begin{align} \label{diagonal}
    \bm{\Sigma} = \bm{D}=  \mbox{diag} \la \lambda_1, \lambda_2,\dots,\lambda_p \ra.
\end{align}
Indeed, for an arbitrary $\bm{\Sigma}$, one can write $\bm{\Sigma}= \bm{O D } \bm{O}^{\top} $, for some orthogonal matrix $\bm{O}$. Thus, the rotated data $\bm{OX}_1, \dots, \bm{OX}_n$ has covariance matrix $\bm{D}$. Furthermore, note that \eqref{Oja} implies
\[
\bm{O}\hat{\bm{v}}_{t} := \big( \bm{I}_p + \eta_{t}\bm{X}_{t} \bm{X}_{t}^{\top} \big) \bm{O}\hat{\bm{v}}_{t-1}.
\]
This means that the output with respected to the rotated data is also rotated by the same matrix $\bm{O}$. Thus, if the algorithm converges for $\bm{D}$, it also converges for $\bm{\Sigma}$. The same argument applies for \eqref{krasulina} as well. Therefore, we assume \eqref{diagonal} and $\bm{v}_1= \bm{e}_1$ from now on.
The next lemma characterizes the improvement per iteration when one applies either \eqref{Oja} or \eqref{krasulina} on a new data point.

At the moment, our analysis requires the data to be bounded. Although our framework for proving Theorems \ref{master} and \ref{conf} can be extended to cases where the coefficients are sub-Gaussian or even sub-exponential, the coefficients in the recursion of Oja's algorithm are sub-Weibul, even for Gaussian data. This arises due to the presence of terms of order $B^4$, which correspond to $\| \bm{X}_t \|^4$.


\begin{prop} \label{oja-err bound}
  { Let $\rho:= \lambda_1-\lambda_2$ be the eigengap. }   Suppose $\bm{\hat{v}}_n$ be the output from either \eqref{Oja} or \eqref{krasulina}. Put
\[
L^{Oja}_t:= \sin^2 \lb \frac{\hat{\bm{v}}_{t}}{\|\hat{\bm{v}}_{t}\|}, \bm{v}_1 \rb .
\]
Then, 
\begin{align} \label{oja recursive}
    L^{Oja}_t \leq \begin{cases}
       &  \lb  1-2 \rho \eta_t \rb L^{Oja}_{t-1} + 2\rho \eta_t \lb L_{t-1}^{Oja} \rb^2 + Q^{Oja}_t + 4B^4 \eta^2_t \quad \text{for} \ \eqref{krasulina}, \\
        &    \lb  1-2 \rho \eta_t \rb L^{Oja}_{t-1} + 2\rho \eta_t \lb L_{t-1}^{Oja} \rb^2 + Q^{Oja}_t + (5B^4+2\eta_tB^6) \eta_t^2 \quad \text{for} \ \ \eqref{Oja}.
    \end{cases}
\end{align}
for all $t \geq 1$, where  $Q_t$ satisfies
\begin{align*}
    \mb{E} \lb Q^{Oja}_t \Big| \mathcal{F}_{t-1} \rb = 0 \quad \text{and} \quad 
    |Q^{Oja}_t| \leq 8B^2 \eta_t \sqrt{L_{t-1}}.
\end{align*}
\end{prop}

The recursions for both \eqref{krasulina} and \eqref{Oja} share the same form, differing only in the term involving $B$ in $\eta_t^2$. However, a key distinction arises between \eqref{oja recursive} and \eqref{sgd-recursive1}: the Oja algorithm includes an additional term, $\eta_t S_{t-1}^2$, which is absent in standard SGD settings. Heuristically, $S_{t-1}$ can be large in the initial steps, leading to minimal improvement in the next iterations. This is the reason why most analyses of Oja’s algorithm are divided into two phases: an initial unstable phase, where the algorithm’s dynamics are erratic, followed by a stable phase, where the dynamics decay at the optimal $O(t^{-1})$ rate. Since the analyses of \eqref{Oja} and \eqref{krasulina} are the same, we will only focus on \eqref{krasulina} in what follows. Our main result in this section is the following theorem.


\begin{thm} \label{oja-conf}
  { Let $\rho:= \lambda_1-\lambda_2$ be the eigengap. }Recall the settings in Proposition \ref{oja-err bound}. Suppose
    $\mb{P} \lb  L_0 \leq 1/4 \rb \geq 1-\delta^3$
    and define 
    $$L= \max \la \left[  128B^4 \log \lb \delta^{-1} \rb^2/ \rho^2  \right]; 32 \ra.$$ 
    Choose step sizes as $\eta_t =  \frac{2}{\rho(t + L)}$. Then,
{
with probability at least  $1-2(e+1)\delta$,
\[
L^{Oja}_t \leq \max \la  \frac{252L}{\log \lb \delta^{-1} \rb}; \frac{1008B^4}{\rho^2}\ra  \cdot \frac{\log \lb \delta^{-1} \rb + 2\log \log (t+9)}{t+L} 
\]
uniformly over all all $t \geq 0$}.
\end{thm}

The constant $L$ in the step sizes $\eta_t$ control the magnitude of each update uniformly and guarantees that, with high probability, all iterates do not leave a neighborhood of radius $1/2$ (with respect to the sine squared distance) around the true eigenvector.

It should be noted that Theorem~\ref{oja-conf} does not follow directly from Theorem~\ref{conf}. Instead, we derive Theorem~\ref{oja-conf} by analyzing an auxiliary process that coincides with $L_t$ on certain ``good events'' and is set to zero otherwise. This construction necessitates the assumption $\mathbb{P}(L_0 \leq 1/4) \geq 1 - \delta^3$, which introduces a cubic dependence on $\delta$, in contrast to the linear dependence required in Theorem~\ref{master}. The cubic dependence on $\delta$ can be relaxed to a $(1+\varepsilon)$-dependence for any $\varepsilon > 0$, at the cost of introducing an additional constant factor $C_\varepsilon > 0$ in the bound, where $C_\varepsilon$ depends only on $\varepsilon$. For the sake of clarity and simplicity, we state the time-uniform bound with a cubic dependence on $\delta$.

{
Let us compare the any-time-valid bound of Theorem~\ref{oja-conf} with the fixed-time bound of Theorem~2.3 in \cite{huang2021streaming}. Below, the notation $\Tilde{\Theta}(a_n)$ means asymptotically equivalent up to logarithm factors of $a_n$.  Their results state that for
\begin{align*}
T_0= \Tilde{\Theta} \lb \frac{B^4}{\delta^2 \rho^2} \rb; \ \beta= \Tilde{\Theta} \lb \frac{B^4}{\delta^2 \rho^2} \rb; \eta_t = \begin{cases}
    \Tilde{\Theta} \lb \frac{1}{\rho T_0} \rb \ &\text{for} \ t \leq T_0; \\
    \Tilde{\Theta} \lb \frac{1}{\rho(t+\beta -T_0)}  \rb  \ &\text{for} \ t \geq T_0 + 1,
\end{cases}
\end{align*}
we have 
\[
\sin^2 \lb \frac{\hat{\bm{v}}_{t}}{\|\hat{\bm{v}}_{t}\|}, \bm{v}_1 \rb \lesssim  \sqrt{\frac{\beta+1}{\beta+T}}
\]
for all fixed $T>T_0$, with probability at least $1-\delta$.
}


Both results share the same dependence on $B^4 / \rho^2$, up to an additional factor of order $\mathrm{polylog}\!\big(\delta^{-1}\big)$. The extra $\log \log t$ term is expected and represents the mild cost of achieving a time-uniform bound. { The step-size schedule adopted here is also similar to that used in the stable phase of \cite{huang2021streaming} in the stable phase $T > T_0$, which exhibits a factor of order $\rho^{-1}$.}


To get a time-uniform bound from Theorem \ref{oja-conf} with a ``cold" start, it is crucial to determine the first time $T_0$ at which $L_{T_0}$ falls below $1/4$, as required by Theorem \ref{oja-conf}. Determining \( T_0 \) depends on the specific structure of the problem, and sharp estimates generally cannot be obtained from an almost-supermartingale of the form \eqref{quantitative R-S lemma}. This is not surprising, as the form of \eqref{quantitative R-S lemma} reduces the analysis of the dynamics of a multi-dimensional system to that of a one-dimensional process. As a result, information about the initial phase of the dynamics is lost. 
To derive $T_0$, we follow the scheme proposed in \cite{huang2021streaming}. Recall that $\hat{\bm{v}}_t$ is the output of the Oja's algorithm according to either \eqref{Oja} or \eqref{krasulina}. The starting time $T_0$ has to be chosen in such a way that 
\[
\mb{P} \lb  \sin^2 \lb \hat{\bm{v}}_{T_0}, \bm{v}_1 \rb \leq \frac{1}{4}  \rb \geq 1- \delta^3.
\]
 According to \cite{huang2021streaming}, we get 
\[
T_0 = O \lb \frac{B^4}{\rho^2} \cdot \log \lb \frac{B^2}{\rho \delta} \rb \rb +  \underbrace{\Tilde{\Theta} \lb \frac{B^4}{\delta^6 \rho^2} \rb}_{H_0}
\]
where the corresponding sequence of step sizes $\eta_t$ can be determined as
  \begin{align*}
            \eta_t = \begin{cases}
               & \Tilde{\Theta} \lb \frac{1}{\rho {H_0}} \rb, \ \text{if $t \leq { H_0}$}, \\
               & \Theta \lb \frac{1}{\rho \lb \beta + t - { H_0} \rb}  \rb , \ \text{if $t \geq { H_0}$}.
            \end{cases}
        \end{align*}   
If a uniformly distributed initialization $\hat{\bm{v}}_0$ is used, then the bound in Theorem~\ref{oja-conf} should be shifted by a factor of $T_0$, that is, by replacing $t$ with $t + T_0$. Note that $T_0$ is the sum of two terms. The second term, denoted $H_0$, corresponds to the initial chaotic phase during which the algorithm explores a sufficiently small region likely to contain the true minimum. In this phase, \cite{huang2021streaming} employs a constant step-size scheme. The first term in the expression for $T_0$ corresponds to the stable phase, where the error decays at the optimal rate and the step sizes are chosen to be of order $1/t$.


\subsection{Robbins-Monroe scheme and a lower bound} \label{Robbins-Monroe}

In this section, we show that the rates obtained in Theorem \ref{master} is sharp by constructing an explicit example that one can not improve the rate for any choice of step sizes that asymptotically of order $\Omega(1/t)$. The goal of this section is solely to demonstrate the sharpness of the rates, and we therefore stick to a one-dimensional example. { We will establish a lower bound in the context of the Robbins–Monro scheme, a stochastic gradient–type algorithm for estimating the root of an equation based on observing unbiased realizations of the underlying function. It was introduced in \cite{robbins1951stochastic} and has been studied extensively in the last few decades. Basic properties and convergence analyses of the  Robbins–Monro scheme can be found in \cite{borkar2008stochastic}; see also \cite{kushner2003stochastic,kushner2010stochastic} and the references therein.}

Consider the classical Robbin-Monroe scheme in dimension one, that is to find the root $\theta$ of the function $M(x)$ through an unbiased estimator $Y(x)$ of $M(x)$. { The Robbins-Monroe scheme works by initializing a value $X_0$ and update recursively by the rule 
\[
X_{k+1}:= X_k -\eta_{k+1} Y(X_{k}).
\] }
We make the following assumptions on $Y$, $M$ and $\la X_n; n \geq 1\ra$:
\begin{itemize}
    \item { $Y(X_n)$ is distributed as $Y(x)$ conditionally on 
    $$X_n=x, X_{n-1}, Y(X_{n-1}), X_{n-2},\dots,X_1, Y(X_1).$$}
    In other words,
    \begin{align} \label{cond dist}
        Y(X_n) \Big| X_n=x, X_{n-1}, Y \lb X_{n-1} \rb, \cdots, Y(X_0), X_0 \stackrel{d}{=} Y(x).
    \end{align} 
    \item $\mb{E} Y(x)= M(x)$ and $\mbox{Var} \lb {Y(x)} \rb = \sigma^2(x) $ and $\sup_{x} \sigma^2(x) <\infty$.
    \item $M(\theta) = 0$ and  $ M^{'}(x) \geq R >0$ for almost everywhere $x \in \mb{R}$.
    \item $M(x)$ is sub-polynomial, that is 
    \begin{align} \label{sub-poly}
        |M(x)| \leq P \lb |x-\theta| \rb
    \end{align} 
    for some polynomial $P$ with positive coefficients.
\end{itemize}

The sub-polynomial condition on $M$ is more general than the commonly imposed sub-linear condition in the literature. The sub-linear condition corresponds to the special case where $P(x) = Ax + B$ for some positive constants $A$ and $B$. In general, the Robbins-Monro scheme may not yield consistent solutions under the sub-polynomial condition. However, Theorem~\ref{master} guarantees consistency provided that the initialization is sufficiently close to the root $\theta$.

Put $L^{RM}_t:= | X_t - \theta |^2$. In what follows, we consider the additive noise structure $Y(x)=M(x)+\xi_x$, for a collection of i.i.d. centered, unit variance random variables $\xi_x$ such that 
\begin{align} \label{R1}
|\xi| \leq R_1
\end{align}
for some constant $R_1>0$. We will assume that the distribution of $\xi_x$'s is continuous to avoid the technical issue that the iterates form a loop with positive probability. 

It is easy to check that such a noise structure satisfies \eqref{cond dist}. With this setting, we have the recursion 

\begin{prop} \label{robbins-monroe constraint}
    Under the settings described above, for all $t\geq 1$,
    \[
    L^{RM}_t \leq \lb 1 - 2R \eta_t \rb L^{RM}_{t-1} + Q^{RM}_t + 2\eta_t^2 \cdot \lb P^2 \lb \sqrt{L^{RM}_{t-1}} \rb + R_1^2 \rb
    \]
    where $\mb{E} \lb Q^{RM}_t | \mathcal{F}_{t-1} \rb=0$, $|Q^{RM}_t| \leq 2\eta_t R_1 \sqrt{S_{t-1}}$ and $P$ is the polynomial in \eqref{sub-poly}.
\end{prop}
Proposition \ref{robbins-monroe constraint} asserts that the iterates of the Robbins-Monroe scheme satisfy a recursion of a form compatible with our main Assumption \ref{assumption:main}. { Thus, the results in Theorem~\ref{master} also hold in this case. Since our focus is on proving a lower bound, we will not attempt to derive an anytime bound in this setting. } 

Next, we show that in the simple case where $P$ is linear, the typical extreme magnitude of the iterates under squared-error loss is exactly of order $\log \log t / t$. 

\begin{thm} \label{lower bound}
    Suppose $P(x)=Ax+B$ for some positive constants $A$ and $B$, where $P(x)$ is as in \eqref{sub-poly}. Assume that  the sequence of step sizes $\la \eta_t; t \geq 1 \ra$ satisfy 
    \[
    \frac{L_1}{t} \leq \eta_t \leq \frac{L_2}{t}
    \]
    for some positive constants $L_1, L_2$. Then, 
    \[
    \mb{P} \lb \limsup_{t \to \infty} \lb \frac{t\cdot L^{RM}_t}{\log \log t} \rb  \geq L \rb     = 1
    \]
     where
    \begin{align} \label{L}
        L:= \frac{\sqrt{L_1}}{4\left[ 1 + L_2\cdot\log \lb 8 \rb   \cdot M^{'}(\theta) \right]}.
    \end{align} 
\end{thm}

Theorem \ref{lower bound} implies that the typical extreme magnitude of $| X_t - \theta |^2$ cannot be asymptotically smaller than order $\log \log t / t$, for any choice of step sizes of order $\Omega(t^{-1})$. In other words, the construction given in the proof of Theorem \ref{master} yields a time-uniform bound of optimal width in a fairly general setting.

\section{Outline of the argument}

Let us demonstrate our technique by explaining the proof of Proposition \ref{pre-planned 1}. Assume $t_0=0$ for simplicity. By choosing step sizes of the form $\eta_t = 2/C_1(t+L)$, the recursion can be rewritten as
\begin{align*}
    L_t \leq b_t \big( X_{0} + U^{*}_t \big),
\end{align*}
where
\begin{align*}
    b_t &:= \frac{L(L-1)}{(t+L-1)(t+L)}, \\
    U^{*}_t &:= \frac{U_t}{b_t} + \frac{U_{t-1}}{b_{t-1}} + \dots + \frac{U_1}{b_1} 
    = \frac{U_t}{b_t} + U^{*}_{t-1}.
\end{align*}
Note that $U^*_t$ is a weighted sum of $U_t$. To obtain a time-uniform bound on the interval $[0,t_1]$, it suffices to derive a maximal inequality of the form
\[
\mathbb{P} \Big( \max_{1\leq k \leq t_1} U^{*}_k \geq x_\delta \Big) \leq \delta.
\]
for some $x_\delta >0$ depending on $\delta$. To do this, we exploit the recursive structure of the problem. Heuristically, we have
\[
L_0 \ \text{small} 
\;\; \underbrace{\implies}_{\text{by \eqref{cond-error} + \eqref{as-bound}}} \;\; U_1 \ \text{small} 
\;\; \underbrace{\implies}_{\text{by \eqref{quantitative R-S lemma}}} \;\; L_1 \ \text{small} 
\;\; \underbrace{\implies}_{\text{by \eqref{cond-error} + \eqref{as-bound}}} \;\; \dots
\]
In order to fully exploit this heuristic, we make use of the following simple observation:

\begin{lemma}
    Let $(M_t)$ be an adapted process with the corresponding filtration $\la \mathcal{F}_t; t \geq 1 \ra$. Let $\varepsilon > 0$ and let $\tau$ be a stopping time such that
    \begin{align*}
        \mathbb{P} \Big( \tau \geq t+1 \,\Big|\, \max_{1 \leq k \leq t} M_t \leq \varepsilon \Big) = 1
    \end{align*}
    for all $1 \leq t \leq T$. Then
    \[
    \mathbb{P} \Big( \max_{1 \leq t \leq T} M_t > \varepsilon \Big)
    = \mathbb{P} \Big( \max_{1 \leq t \leq T} M_{\tau \wedge t} > \varepsilon \Big).
    \]
\end{lemma}

The lemma above has been used in the analysis of Oja's algorithm by \cite{chou2020ode}. To use it, define the stopping time
\[
\tau = \inf \{ t \geq t_0 : L_t > r_t \}
\]
for some sequence $\la r_t \ra$ that serves as tuning parameters to obtain sharp thresholds.

Observe that the stopped process $U^{*}_{\tau \wedge t}$ has manageable conditional mean and magnitude, given by
\begin{align*}
     U^*_{\tau \wedge t} - U^*_{\tau \wedge (t-1)}
     &= \mathbf{1}_{\{ \tau \geq t \}} \cdot \big( U^*_{t} - U^*_{t-1} \big) \\
     &= \frac{ \mathbf{1}_{\{ \tau \geq t \}} }{b_t} \cdot U_t.
\end{align*}
Thus, using \eqref{cond-error 2} and \eqref{as-bound 2}, and noting that $L_t \leq r_t$ on the event $\{ \tau \leq t \}$, we get that 
\begin{align*}
   | U^*_{\tau \wedge t} - U^*_{\tau \wedge (t-1)}| &= \frac{ \mathbf{1}_{\la  \tau \geq t \ra}}{b_t} \cdot \left| U_t \right| 
\leq  \frac{2C_3(t+L-1)\sqrt{r_{t-1}}}{C_1 (L-1)L } + \frac{4C_2(t+L-1)}{C_1^2(t+L)(L-1)L}
\end{align*}
and
\begin{align*}
   \left| \mb{E} \left[ U^*_{\tau \wedge t} - U^*_{\tau \wedge (t-1)} | \mathcal{F}_{t-1}  \right] \right| &\leq \frac{4 C_2}{C_1^2 \lb t  +L \rb} \cdot \frac{t+L-1}{(L-1)L}.
\end{align*}
A standard application of Azuma–Hoeffding’s inequality along with a telescoping argument then yields
\[
\mathbb{P} \Big( \max_{1\leq k \leq t_1} U^{*}_{k \wedge \tau} \geq y(\delta,r_t) \Big) \leq \delta
\]
for some function $y(\delta,r_t)$ of $\delta$ and $r_t$. We then tune the sequence $\la r_t \ra$ in a way such that the lemma above applies to deduce that
\[
\mathbb{P} \Big( \max_{1\leq k \leq t_1} U^{*}_{k} \geq y(\delta,r_t) \Big)
= \mathbb{P} \Big( \max_{1\leq k \leq t_1} U^{*}_{k \wedge \tau} \geq y(\delta,r_t) \Big) \leq \delta.
\]
The conclusion of Proposition \ref{pre-planned 1} follows immediately.

The proof of Theorem \ref{master} also follows from the same scheme, but requires a different choice of the tuning sequence $r_t$ to accommodate piecewise-constant step sizes. The general form of \eqref{cond-error} and \eqref{as-bound} makes it quite challenging to construct an explicit step-size schedule. Finally, note that the $\log \log$ term comes from a peeling argument: we apply Proposition \ref{pre-planned 1} repeatedly on small intervals.

One major advantage of the approach above is that it can yield a maximal inequality by fully exploiting the hierarchical structure of the recursion, while requiring only simple concentration inequalities in the analysis. We believe that a similar type of maximal inequality can also be obtained by using the mixture method, as in \cite{howard2020time,howard2021time} (see also the references therein). However, it is not clear to us which choices of the prior would yield a comparable bound involving only $C_3^2/C_1^2$ and $C_2/C_1^2$.

\section{Conclusion and discussion} \label{sec:conclusion}
In this paper, we develop a new technique based on a variant of the classical Robbins--Siegmund lemma \cite{robbins1971convergence} for obtaining time-uniform bounds for a broad class of almost supermartingale-like processes, where the remainder terms at each iteration are bounded by a polynomially growing function of the step sizes and previous iterates. We prove that the rates achieved by our technique are optimal and provide explicit constructions of step size schedules for processes with simple landscapes (i.e., without spurious local minima). As applications, we use the method to derive time-uniform bounds for projected SGD, Oja’s algorithm, and the Robbins-Monro scheme. We conclude by outlining several potential extensions of our framework that present interesting problems for future research:

\begin{itemize}
    \item \textbf{Streaming $k$-PCA.} In our analysis of Oja’s algorithm, we considered only the estimation of the principal eigenvector. A natural extension is to handle the streaming $k$-PCA problem by deriving a recursive estimate analogous to Proposition~\ref{oja-err bound}. We believe the argument will similar, but more technically involved. The key step would be to derive an analog of Proposition \ref{oja-err bound}.

    \item \textbf{Relaxed moment conditions.} Our results rely in a fundamental way on the bounds \eqref{cond-error} and \eqref{as-bound} of our recursive Assumption \ref{assumption:main}, which deliver the optimal $\log \log t / t$ scaling. A natural weakening of Assumption \ref{assumption:main} is to allow for the coefficients $A_i$'s and $B_i$'s in  \eqref{cond-error} and \eqref{as-bound} to be random.   Determining the moment conditions on those coefficients that still guarantee a $\log \log t / t$ rate is  both interesting and challenging. Although our results remain valid (possibly with worse constants) when the coefficients are conditionally sub-Gaussian (see Appendix A in the supplement for more details), it remains unclear whether the same holds under sub-exponential, or even the more challenging sub-Weibull, settings. 
    
    Notably, in the simplified case of \eqref{cond-error 2} and \eqref{as-bound 2}, the same rate (with possibly worse  constants) still holds when $C_3$ is sub-Gaussian and $C_2$ is sub-exponential---matching the SGD rate under sub-Gaussian error recently setting studied in \cite{madden2024high}. { We believe that settings involving noise distributions with tails heavier than sub-exponential can be addressed using a combination of truncation and path coupling arguments. This represents an interesting direction for future work.}

    \item \textbf{Low-rank constrained least squares.} Our technique can be applied to other problems of similar structure. For example, the low-rank constrained least squares framework of \cite{de2015global} encompasses several important tasks, including matrix completion, phase retrieval, and subspace tracking.  

    \item \textbf{Stochastic heavy-ball algorithms.} Another promising direction is to extend our results to stochastic heavy-ball methods; see, for example, \cite{sebbouh2021almost,liu2023last} and references therein. The main technical challenge is that these algorithms involve an additional momentum term at each step, leading to a two-dimensional process rather than a scalar one.  
\end{itemize}

\section{Acknowledgments}
We thank Shubhanshu Shekhar for productive discussions during the initial phase of the project and Aaditya Ramdas for helpful comments and for pointing out the reference \cite{whitehouse2023time}. TP and PS gratefully acknowledge NSF grants 2217069 and 2019844.

\section{Proof of Theorem \ref{master}}
The following result is a variant of Lemma 6.9 in \cite{chou2020ode}, previously used in the analysis of certain Oja-type algorithms, and plays a crucial role in our proof. It shows that the crossing probability of a process can be controlled by that of a stopped process via a suitably chosen stopping time.
 \begin{lemma} \label{stopped}
    Let $M_t$ be an adapted process with the corresponding filtration $\la \mathcal{F}_t; t \geq 1 \ra$. Let $\ve>0$ and $\tau$ be a stopping time with respect to $\la \mathcal{F}_t; t \geq 1 \ra$ such that
    \begin{align} \label{stopped cond}
    \mb{P} \lb  \tau \geq t+1 \Big| \max_{1 \leq k \leq t} M_t \leq  \ve  \rb =1
    \end{align}
    for all $1 \leq t \leq T$. Then,
    $$\mb{P} \lb \max_{1 \leq t \leq T} M_t > \ve \rb = \mb{P} \lb \max_{1 \leq t \leq T} M_{\tau \wedge t}  > \ve \rb.$$
\end{lemma}

\noindent \textbf{Proof of Lemma \ref{stopped}.} One direction is obvious because 
\[
 \max_{1 \leq t \leq T} M_t \geq \max_{1 \leq t \leq \tau \wedge T} M_t =  \max_{1 \leq t \leq T} M_{\tau \wedge t}
\]
which implies 
\[
\mb{P} \lb \max_{1 \leq t \leq T} M_t > \ve \rb \geq \mb{P} \lb \max_{1 \leq t \leq T} M_{\tau \wedge t}  > \ve \rb.
\]
To prove the other direction, write 
\begin{align*}
    \mb{P} \lb \max_{1 \leq t \leq T} M_{\tau \wedge t}  > \ve \rb  &= \sum_{k=1}^{T} \mb{P} \lb \tau = k, \max_{1 \leq i \leq k} M_i > \ve  \rb + \mb{P} \lb \tau \geq T+1, \max_{1 \leq i \leq T} M_i > \ve   \rb \\
    &=  \sum_{k=1}^{T} \left[ \mb{P} \lb \tau = k  \rb - \underbrace{\mb{P} \lb \tau = k, \max_{1 \leq i \leq k} M_i \leq \ve  \rb}_{=0} \right] + \mb{P} \lb \tau \geq T+1, \max_{1 \leq i \leq T} M_i > \ve   \rb \\
    &= \sum_{k=1}^{T} \mb{P} \lb \tau = k  \rb + \mb{P} \lb \tau \geq T+1, \max_{1 \leq i \leq T} M_i > \ve   \rb \\
    &= \mb{P} \lb \max_{1 \leq i \leq T} M_i > \ve   \rb  \underbrace{- \mb{P} \lb \tau \leq T,  \max_{1 \leq i \leq T} M_i > \ve   \rb  + \mb{P} \lb \tau \leq T \rb}_{\geq 0}.
\end{align*}
The proof is completed. $\hfill$ $\square$



\begin{lemma}[Azuma-Hoeffding's maximal inequality] \label{freedman}
    Let $M_t$ be a process adapted to the filtration $\la \mathcal{F}_t; t \geq 0 \ra$ and such that $M_0=0$ almost surely. For every $t$, let $\mu_t$ and $\sigma_t^2$ satisfy 
    \begin{align}
     |M_t - M_{t-1}| \leq c_t \label{M-asbound}\\
        \left|  \mb{E} \lb M_t - M_{t-1} \Big| \mathcal{F}_{t-1} \rb   \right| \leq \mu_t. \label{M-condmean} 
    \end{align}
    Then, for all $T \geq 1 $ and $\delta \in (0,1)$,   
    $$\mbox{P} \lb \exists t \in [1,T]: M_t \geq \sqrt{2} V_T(\delta) + \sum_{i=1}^{t} \mu_i     \rb \leq \delta, $$
    where
    \begin{align}
        V_T(\delta) :=  \sqrt{\log \lb \frac{1}{\delta}\rb  \sum_{t=1}^{T} c_t^2}  \label{var profile}
    \end{align}
\end{lemma}

\noindent \textbf{Proof of Lemma \ref{freedman}.} 
Define $d_0:=0$ and $d_k:= M_{k}-M_{k-1} $. It is easy to check that 
$$M_k= \sum_{i=1}^{k} \left[ d_i - \mb{E} \lb d_i | \mathcal{F}_{i-1} \rb \right] + \sum_{i=1}^{k}  \mb{E} \lb d_i | \mathcal{F}_{i-1} \rb. $$ 
Let $h_i:= d_i - \mb{E} \lb d_i | \mathcal{F}_{i-1} \rb$. It is easy to check that $\la h_i; 1\leq i \leq T \ra$ forms a martingale difference sequence and 
\begin{align*}
    \Big| \mb{E} \lb d_i | \mathcal{F}_{i-1} \rb \Big| &\leq \mu_i; \\
     \mb{E} \lb e^{\lambda h_i} | \mathcal{F}_{i-1} \rb &\leq e^{\lambda^2 c_i^2/2}
\end{align*}
almost surely, for all $\lambda >0$. Note that the second estimate in the display above follows from Chernoff's bound and the fact that $|h_i|\leq 2c_i$, which is due to \eqref{M-asbound}. Thus, we have 
\begin{align*}
    & \mb{P} \lb \exists t \in [1,T]: M_t \geq \sqrt{2}V_T(\delta) + \sum_{i=1}^{t} \mu_i  \rb\\ 
    \leq & \mb{P} \lb \exists t \in [1,T]: \sum_{k=1}^{t} h_k \geq \sqrt{2}V_T(\delta) \rb \\
    = & \mb{P} \lb \exists t \in [1,T]: \sum_{k=1}^{t} h_k \geq \sqrt{2 \log \lb \frac{1}{\delta}\rb \cdot \sum_{k=1}^{T} c_k^2 }  \rb \\
    \leq & \exp \lb- \frac{2\log \lb \frac{1}{\delta}\rb \cdot \sum_{k=1}^{T} c_k^2}{2\sum_{k=1}^{T} c_k^2} \rb = \delta
\end{align*}
where the first inequality follows from \eqref{M-condmean}, the second inequality follows from Azuma-Hoeffding's inequality. The proof is completed. 

$\hfill$ $\square$

\subsection{Step 1: Maximal inequality}
The following result is our main technical result. It gives a sharp, short-time maximal inequality for a stochastic process that satisfies \eqref{quantitative R-S lemma} with constant step sizes in a finite interval. 

\begin{lemma} \label{maximal short time}
    Suppose $X_0, X_1,\dots,X_n$ is a sequence of random variables that satisfy \textcolor{violet}{ALE: Do you mean Assumption 1?} \eqref{quantitative R-S lemma} with constant step sizes $\eta_n=\eta$. Let $r \in (0,1)$ and $\delta \in (0,1)$, and put 
     \begin{align} \label{D_delta}
     D_{\delta}(r) :=  \frac{8D}{(1-C_1 \eta)^n} \cdot \lb \sqrt{\log \lb \frac{1}{\delta} \rb}  \cdot \lb \sqrt{\eta r} + \sum_{i=1}^{m} \eta^{c_i} r^{d_i} \rb + \eta  + \sum_{i=1}^{m} \eta^{a_i} r^{b_i} \rb,
     \end{align}
     where
     \begin{align} \label{D}
     D:= \max_{1 \leq i \leq m} \la \frac{A_i}{C_1} \ra \vee \max_{1 \leq i \leq m} \la \frac{\sqrt{m+1}\cdot B_i}{\sqrt{C_1}} \ra \vee \frac{C_2}{C_1} \vee \frac{C_3\sqrt{m+1}}{\sqrt{C_1}}.
     \end{align}
     Then, for any $\ve_0 >0$, we have that
    $$\mb{P} \lb \exists t \in [1,n] : X_t  \geq \lb 1- C_1\eta \rb^{t} \lb \ve_0 + D_{\delta} \rb, A  \rb \leq\delta $$
    whenever $\ve_0 + D_{\delta}(r) \leq r $ and $A \subset \la X_{0} \leq \ve_0 \ra$ is such that $\mb{P} \lb A \rb > 0$.
\end{lemma}

\noindent \textbf{Proof of Lemma \ref{maximal short time}.} For each $t \geq 1$, define 
\begin{align*}
    U^{*}_t &:= \lb 1 - C_1\eta \rb^{-t} \Big[ U_t + \lb 1 - C_1\eta \rb \cdot U_{t-1} + \cdots + \lb 1 - C_1\eta \rb^{t} U_{0} \Big]  
\end{align*}
and set $U^{*}_0=0$.
Note that with this notation, we have
\begin{align} 
    U^{*}_t  &= \lb 1 - C_1\eta \rb^{-t}U_{t} + U^{*}_{t-1}; \label{recursive-U*} \\
    X_{t}  & \leq (1- C_1 \eta)^{t} \lb X_0 + U^{*}_t \rb. \label{recursive-U*-2}
\end{align}
Our strategy is to establish via Lemma \ref{stopped} the maximal inequality
\begin{align} \label{maximal ineq}
 \mb{P} \lb \max_{1 \leq t \leq n} U^{*}_{ t}   \geq D_{\delta} | A \rb \leq \frac{\delta}{\mb{P}(A)}.
\end{align}
Towards that end, we will first show that 
\begin{align} \label{cond-maximal}
    \mb{P} \lb \max_{1 \leq t \leq n} U^{*}_{\tau \wedge t}   \geq D_{\delta} | A \rb \leq \frac{\delta}{\mb{P}(A)},
\end{align}
where $\tau$ is the stopping time 
\[
\tau = \inf\la t \geq 0: X_t \geq r  \ra.
\]
Note that, for any arbitrary stopping time $\tau$, we have the relation 
\[
 U^{*}_{\tau \wedge t} - U^{*}_{\tau \wedge (t-1)} = \mathbf{1}_{\la \tau \geq t \ra} \lb U^{*}_{t} - U^{*}_{t-1} \rb.
\]
Thus,
\begin{align*} 
     \left| U^{*}_{ \tau \wedge t} - U^{*}_{ \tau \wedge (t-1)}  \right|  &= \left| \mathbf{1}_{\la \tau \geq t \ra} \lb U^{*}_{t} - U^{*}_{t-1} \rb \right| \leq c_t  ;  \\
    \left| \mb{E} \lb  U^{*}_{\tau \wedge t} - U^{*}_{\tau \wedge (t-1)} \Big| \mathcal{F}_{t-1} \rb \right| & = \Big| \mathbf{1}_{\la \tau \geq t \ra} \underbrace{\mb{E} \lb U^{*}_{t} - U^{*}_{t-1} \Big| \mathcal{F}_{t-1} \rb}_{\text{by using } \ \eqref{recursive-U*}}  \Big|\\     
    &= (1- C_1\eta)^{-t}  \left| \mb{E} \lb U_t | \mathcal{F}_t \rb  \right|
    \leq  \mu_t,
\end{align*}
where 
\begin{align*}
    \mu_t &:= (1- C_1\eta)^{-t} \cdot \lb C_2\eta^2 + \sum_{k=1}^{m} A_k\cdot  \eta^{1+a_k} r^{b_k} \rb; \\
    c_t &:= (1- C_1\eta)^{-t} \cdot \lb C_3\eta\sqrt{r} + \sum_{k=1}^{m} B_k \eta^{1/2+c_k} r^{d_k} \rb.
\end{align*}
By Lemma \ref{freedman}, we obtain that
\begin{align*}
    \mb{P} \lb \exists t \in [1,n]: U^{*}_{\tau \wedge t} \geq 8 \sqrt{ \log \lb \delta^{-1}  \rb \cdot\sum_{k=1}^{n} c^2_k}  + \sum_{k=1}^{t} \mu_k \rb \leq \delta.
\end{align*}

Next, one can easily check that
\begin{align*}
    \sum_{k=1}^{t} \mu_k &= \lb C_2 \cdot \eta^2 + \sum_{k=1}^{m} A_k\cdot  \eta^{1+a_k} r^{b_k} \rb \sum_{k=1}^{t} (1-C_1\eta)^{-k} \\
    &=  \lb C_2 \cdot \eta^2 + \sum_{k=1}^{m} A_k\cdot  \eta^{1+a_k} r^{b_k} \rb (1-C_1 \eta)^{-1} \frac{1}{1/(1-C_1\eta)-1} \left[  \frac{1}{(1-C_1\eta)^t} -1  \right] \\
    &\leq \lb C_2\eta^2 +  \sum_{k=1}^{m} A_k\cdot  \eta^{1+a_k} r^{b_k}\rb \frac{(1-C_1 \eta)^{-t}}{C_1 \eta}\\
    &\leq   D \lb \eta + \sum_{k=1}^{m}\eta^{a_k} r^{b_k}\rb (1-C_1 \eta)^{-n}
\end{align*}
where $D$ is defined in \eqref{D}. Similarly, 
\begin{align*}
     \sum_{k=1}^{n} c_k^2 & \leq (m+1)\cdot \lb \eta^2r + \sum_{k=1}^{m} \eta^{2c_i+1} r^{2d_i} \rb   \cdot  \frac{ C_3^2 \vee \lb \max_{1\leq k \leq m} B_k^2  \rb}{2C_1\eta - C_1^2 \eta^2} \cdot (1-C_1\eta)^{-2n}\\
     &\leq \underbrace{(m+1)\cdot  \frac{ C_3^2 \vee \lb \max_{1\leq k \leq m} B_k^2  \rb}{C_1}}_{\leq D^2} \cdot \lb \eta r + \sum_{k=1}^{m} \eta^{2c_i} r^{2d_i} \rb \cdot (1-C_1\eta)^{-2n} \\
     &\leq D^2 \lb \eta r + \sum_{k=1}^{m} \eta^{2c_i} r^{2d_i} \rb \cdot (1-C_1\eta)^{-2n}.
\end{align*}
Consequently,
\begin{align*}
    \mb{P}_A \lb \exists t \in [1,n]: U^{*}_{\tau \wedge t} \geq D_{\delta} \rb 
    \leq & \frac{\delta}{\mb{P} \lb A \rb}
\end{align*}
where $\mb{P}_A$ denotes the conditional probability on $A$ and 
\begin{align*} 
    D_\delta:=  \frac{8D}{(1-C_1 \eta)^n} \cdot \lb \sqrt{\log \lb \frac{1}{\delta} \rb}  \cdot \lb \sqrt{\eta r} + \sum_{i=1}^{m} \eta^{c_i} r^{d_i} \rb + \eta  + \sum_{i=1}^{m} \eta^{a_i} r^{b_i} \rb.
\end{align*}
Note that the inequality above implies \eqref{cond-maximal} since
\begin{align*}
    (1-C_1\eta)^{-t} \leq (1-C_1\eta)^{-n}.
\end{align*}
To deduce \eqref{maximal ineq} from \eqref{cond-maximal}, we need to check that
\[
\mb{P}_A \lb \tau \geq k+1 \Big|  \max_{1 \leq t \leq k} U^{*}_t \leq  D_{\delta}  \rb =1.
\]
Note that Lemma \ref{stopped} is being applied to the conditional measure $\mb{P}_A$ instead of $\mb{P}$. The last display is true since conditionally on $\la \max_{1 \leq t \leq k} U^{*}_t \leq  D_{\delta} \ra \cap A$, \eqref{recursive-U*-2} gives
\begin{align*}
    X_k \leq \lb 1 - C_1\eta  \rb^{k} \lb  X_{0} + U^{*}_{k}  \rb \leq \lb 1 -C_1\eta  \rb^{k} \lb  \ve_0 + D_{\delta}  \rb \leq r
\end{align*}
 since $\ve_0 + D_{\delta} \leq r$ and $A \subset \la X_0 \leq \ve_0 \ra$. Consequently,
$$ \mb{P} \lb \max_{1 \leq t \leq n} U^{*}_{ t}   \geq D_{\delta} | A \rb \leq \frac{\delta}{\mb{P}(A)}.$$
The proof is completed. $\hfill$ $\square$

\subsection{Step 2: Induction argument} \label{stitching}
Let $h_i = h_0 \cdot 2^{-i}$ for some constant $h_0$ to be determined later. The value $h_0$ correspond to $r$ in the statement of Theorem \ref{master}. We want to derive a sequence of epochs $\la  t_i; i \geq 0 \ra$ such that with probability at least $1-2\delta$,  $L_{t} \lesssim h_i$ for all $i \geq 0$ and $t \in [t_{i-1}+1,t_i]$. 

To visualize the argument below, it is helpful to refer to Figure~\ref{adaptive stepsize}, where we plot the width of the confidence sequence over time, assuming $h_0 = 1$. The error dynamic of our analysis can be summarized as follows
\begin{align*}
    \underbrace{t_0=0}_{h_0}    \to  \underbrace{t_0+1}_{\leq (1-C_1\eta_0)2h_0} \to \dots \to \underbrace{t_1-1}_{\leq (1-C_1\eta_0)^{t_1-t_0-1} 2h_0} \to \underbrace{t_1}_{\leq 2h_0/4=h_1} \to \dots
\end{align*}

\begin{figure}[H] 
    \centering
    \begin{tikzpicture}[scale=1.35]
    \draw[->] (-0.5,0) -- (8,0) node[right] {time};
    \draw[->] (0,-0.5) -- (0,7) node[above] {width};
    
    \node[left] at (0,4) {2};
    \node[left] at (0,2) {$h_0=1$};
    \node[left] at (0,1) {$h_1=0.5$};
         \node[left] at (0,1/2) {$h_2=0.25$};

    \node[below] at (1,0) {1};
    \node[below] at (3,0) {$t_1$};
    \node[below] at (4,0) {$t_{1}+1$};
    \node[below] at (7,0) {$t_2$};
    \node[below,font=\footnotesize] at (7.5,0) { $t_2 +1$};
    \draw[dashed] (1,0) -- (1,4);
    \draw[dashed] (3,0) -- (3,1);
        \draw[dashed] (4,0) -- (4,2);
            \draw[dashed] (7,0) -- (7,0.5);

    \draw[dashed] (3,1) -- (0,1);
    \draw[dashed] (7,0.5) -- (0,0.5);
            \draw[dashed] (4,2) -- (0,2);
            \draw[dashed,blue,thick] (7.25,1) to[out=-60,in=-180] (8,0.125);
        \draw[dashed] (7.25,1) to (7.25,0);
            \draw[blue,thick] (0,2) -- (1,4);
            \draw[blue,thick] (1,4) to[out=-80,in=-180] (3,1);
            \draw[blue,thick] (3,1) -- (4,2);
            \draw[blue,thick] (4,2) to[out=-60,in=-180] (7,0.5);
            \draw[blue,thick] (7,0.5)--(7.25,1);
        \fill (0,0) circle (1.5pt);
                \fill (1,0) circle (1.5pt);
                                \fill (4,0) circle (1.5pt);
        \fill (3,0) circle (1.5pt);
                                        \fill (7,0) circle (1.5pt);
    \fill (3,1) circle (1.5pt);
        \fill (1,4) circle (1.5pt);
        \fill (4,2) circle (1.5pt);
    \fill (7,0.5) circle (1.5pt);
            \fill (0,1) circle (1.5pt);
              \fill (0,2) circle (1.5pt);
    \fill (7.25,1) circle (1.5pt);
    \fill (7.25,0) circle (1.5pt);

    
    \node[below left] at (-0.2,0) {$t_0 = 0$};
    
\end{tikzpicture}

    \caption{Errors induced by adaptive step sizes}  \label{adaptive stepsize}
\end{figure}
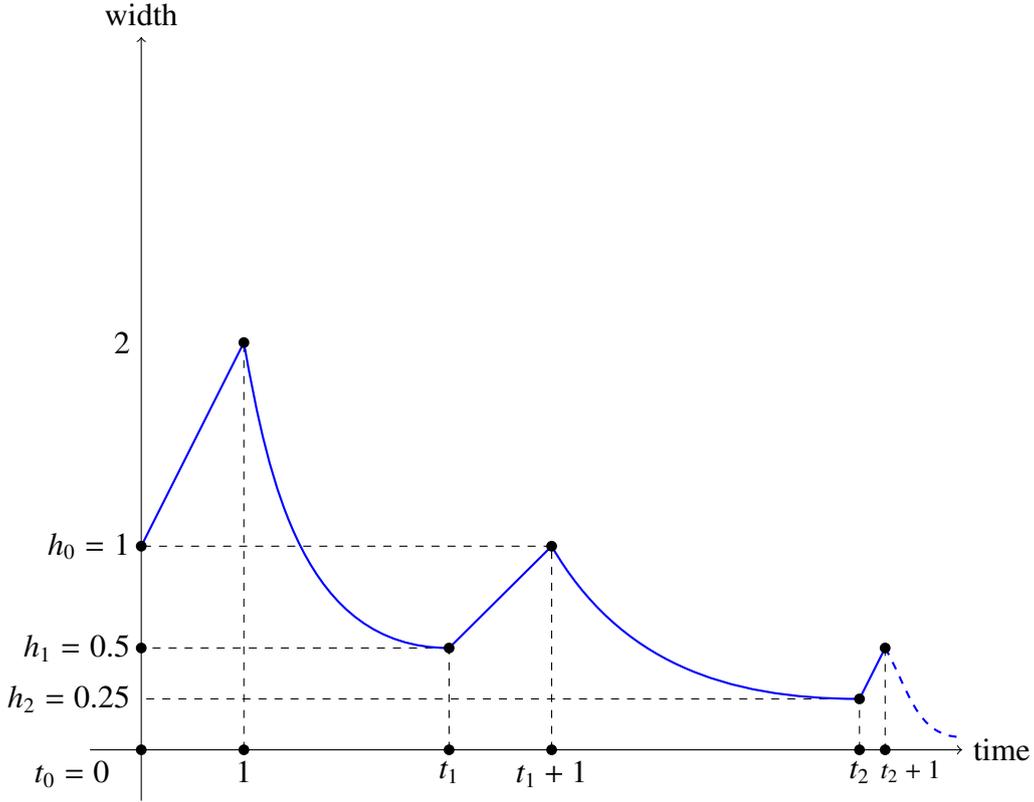

The width of the confidence sequence will be a function $f$ such that $a_i \leq f(t_i)$. Put
\begin{align*}
 t_0 :&=0, \\
 \delta_i :&=  \frac{\delta}{(i+10)^2}.
\end{align*}
 For a small constant $\kappa \leq 1$ to be chosen later, let us define
\begin{align}
    \eta_i :&= \frac{\kappa h_{i-1}}{\log \lb \delta_i^{-1} \rb} \label{eta}
\end{align}
and choose the epochs $t_i$ such that
\begin{align}
     \frac{1}{16} &\leq \lb 1 - C_1\eta_i \rb^{t_i-t_{i-1}} \leq \frac{1}{8} \label{t-constraint}.
\end{align}
The step sizes are chosen to be the constant $\eta_i$ on the interval $[t_{i-1}+1, t_i]$. In what follows, we will determine $\lb \kappa, h_0 \rb$ by applying Lemma \ref{maximal short time} repeatedly. We will also require
\begin{align} \label{upper eta}
\sup_{i \geq 1} \eta_i \leq \frac{1}{16C_1}.
\end{align}

\noindent \textit{ \underline{Step 2a: Maximal inequality.}}  Let us show that for all $i \geq 1$,
\begin{align} \label{interval bound}
       \mb{P} \lb \exists t \in [t_{i-1}+1,t_i] : L_t  \geq 4\lb 1- C_1\eta_i \rb^{t-t_{i-1}} h_i \Big| L_{t_{i-1} } \leq h_{i-1}  \rb \leq \frac{\delta_i}{\mb{P} \lb L_{t_{i-1}} \leq h_{i-1} \rb}
   \end{align}
   for some choices of $\kappa$ in \eqref{eta} and $a_0$  that depend only on the parameters in \eqref{cond-error} and \eqref{as-bound}.

To choose $\kappa, a_0$ such that \eqref{interval bound} holds, we apply Lemma \ref{maximal short time} to $\la L_t \ra$ to get
 $$ \mb{P} \lb \exists t \in [t_{i-1}+1,t_i] : L_t  \geq \lb 1- C_1\eta_i \rb^{t} \lb h_{i-1} + D(r_i) \rb \Big| L_{t_{i-1} } \leq h_{i-1}  \rb \leq \frac{\delta_i}{\mb{P} \lb L_{t_{i-1}} \leq h_{i-1} \rb} $$
    whenever 
    \begin{align} \label{r_i}
        h_{i-1} + D(r_i) \leq r_i.
    \end{align} 
    Here $D(r_i)$ is defined to be
    \[
    D(r_i) := \frac{8D}{(1-C_1 \eta_i)^{t_i - t_{i-1}}} \cdot \left[ \sqrt{\log \lb \frac{1}{\delta_i} \rb}  \cdot \lb \sqrt{\eta_i r_i} + \sum_{k=1}^{m} \eta_i^{c_k} r_i^{d_k} \rb + \eta_i  + \sum_{k=1}^{m} \eta_i^{a_k} r_i^{b_k} \right]
    \]
    where $D$ is the constant in \eqref{D}. 

    One does not have to worry about the existence of $r_i$ in \eqref{r_i}. Let us fix the choice $r_i = 2 h_{i-1}$ and choose $\kappa$ in \eqref{eta} and $a_0$ such that $h_{i-1} + D(r_i) \leq r_i $ for all $i\geq 1$. Observe that
    \begin{align*}
        \sqrt{\log \lb \delta_i^{-1} \rb} \sqrt{\eta_i r_i} &= \sqrt{2\kappa} h_{i-1}; \\
         \sqrt{\log \lb \delta_i^{-1} \rb} \eta_i^{c_k} r_i^{d_k}  &=  \sqrt{\log \lb \delta_i^{-1} \rb} \underbrace{\lb. \frac{\kappa}{2\log \lb \delta_i^{-1} \rb} \rb^{c_k}}_{\leq 1} 2^{c_k+d_k} \cdot h_{i-1}^{c_k+d_k-1} \cdot h_{i-1} \\
         &\leq \left[ 2^{c_k+d_k} \cdot \frac{\sqrt{\log \lb \delta_i^{-1} \rb}}{2^{(i-1)(c_k+d_k-1)}} \cdot h_0^{c_k+d_k-1} \right] h_{i-1}; \\
    \eta_i^{a_k} r_i^{b_k} & = \underbrace{\lb  \frac{\kappa}{2\log \lb \delta_i^{-1} \rb} \rb^{a_k}}_{\leq 1} \cdot 2^{a_k+b_k} \cdot h_{i-1}^{a_k+b_k-1} \cdot h_{i-1} \\
    &\leq \left[ 2^{a_k+b_k} \cdot h_0^{a_k+b_k-1} \right] \cdot h_{i-1}.
    \end{align*}
    Note that $\lb 1 -C_1 \eta_i \rb^{t_i -t_{i-1}} \geq 1/16$ due to \eqref{t-constraint}. Thus, to get \eqref{r_i}, we need to choose $K$ and $a_0$ such that 
    \begin{align*}
        128D \left[ \sqrt{2\kappa} + \sum_{k=1}^{m} 2^{c_k+d_k} \cdot \frac{\sqrt{\log \lb \delta_i^{-1} \rb}}{2^{(i-1)(c_k+d_k-1)}} \cdot h_0^{c_k+d_k-1} + \frac{\kappa}{\log \lb \delta_i^{-1} \rb} + \sum_{k=1}^{m} 2^{a_k+b_k} \cdot h_0^{a_k+b_k-1} \right] \leq 1.
    \end{align*}
    It suffices to pick $\kappa$ and $a_0$ such that 
    \begin{align*}
        \sqrt{2\kappa} &  \leq \frac{1}{(2m+2)128D}; \\
        h_0^{c_k+d_k-1} &\leq \frac{1}{(2m+2)128D} \cdot 2^{-c_k-d_k} \cdot \min_{i\geq 1} \la \frac{2^{(i-1)(c_k+d_k-1)}}{\sqrt{ \log \lb \delta_i^{-1} \rb}} \ra, \forall k \in [1,m]; \\
        \kappa & \leq \frac{1}{(2m+2)128D}; \\
        h_0^{a_k+b_k -1} &\leq \frac{2^{-a_k-b_k}}{(2m+2)128D}, \forall k \in [1,m].
    \end{align*}
   From the first and third inequalities, it is clear that one can choose $\kappa$ as a function of $m, D$:
   \[
   \kappa := \min \la 1;  \frac{1}{2}\left[  \frac{1}{(2m+2)128D} \right]^2; \frac{1}{(2m+2)128D}  \ra.
   \]
   Moreover, one can also choose $h_0$ such that the second and last inequalities also hold because 
   \[
    \min_{i\geq 1} \la \frac{2^{(i-1)(c_k+d_k-1)}}{\sqrt{ \log \lb \delta_i^{-1} \rb}} \ra  > 0 \ \text{and} \ (c_k+d_k-1) \wedge (a_k+b_k-1)>0
   \]
   for all $1 \leq k \leq m$. Such a specific choice is
   \[
   h_0 := \min_{1 \leq k \leq m} \la \exp \lb \frac{A_k \lb\delta \rb}{c_k+d_k-1} \rb \wedge \exp \lb \frac{B_k}{a_k+b_k-1} \rb \wedge \frac{\log \lb  \delta^{-1} \rb}{16 C_1 \kappa}\ra
   \]
   where the last term is due to the constraint on step sizes \eqref{upper eta} and
   \begin{align*}
       A_k(\delta):&= \log \lb \frac{1}{(2m+2)128D} \cdot 2^{-c_k-d_k} \cdot \min_{i\geq 1} \la \frac{2^{(i-1)(c_k+d_k-1)}}{\sqrt{ \log \lb \delta_i^{-1} \rb}} \ra \rb; \\
       B_k :&= \log \lb \frac{2^{-a_k-b_k}}{(2m+2)128D} \rb.
   \end{align*}

   Consequently, we have \eqref{interval bound} since $h_{i-1} + D(r_i) \leq r_{i} = 4h_i$. We now set $r \equiv h_0$, where $r$ is the constant in the conclusion of Theorem \ref{master}. From now on, we  assume that
   \[
   \mb{P} \lb  L_0 \leq h_0 \ra \geq 1- \delta.
   \]

   \noindent \textit{\underline{Step 2b: Valid coverage of the confidence sequence.}} Define 
   \begin{align*}
       r(i,t):= \begin{cases}  4\lb 1 - C_1\eta_i \rb^{t-t_{i-1}} h_i, \ &\text{for} \ t \in [t_{i-1}+1,t_i], \\
        h_i,  \ \ \   \  \ \ \ \ \ \ \ \ \ \ \ \ \ \ \  \ \ \ \ \ \ \ \ \ \ \  \ &\text{for} \ t=t_i.
       \end{cases}
   \end{align*}
   Note that $r\lb i, t_i \rb=h_i$ for all $i \geq 0$ by the defition above. Put
   \begin{align} \label{r*}
   r^{*}_t := \begin{cases} r(i,t) \ &\text{for} \ t \in  [t_{i-1}+1,t_i], \\
   h_0 \ &\text{for} \ t=0.
   \end{cases}
   \end{align}
   By using the union bound and \eqref{interval bound}, we have 
   \begin{align*}
       \mb{P} \lb  \forall t \geq 0: L_t \leq r^{*}_t \rb &\geq 1 - \mb{P} \lb L \leq a_0 \rb -\sum_{i=1}^{\infty}   \mb{P} \lb \la \exists t \in [t_{i-1}+1,t_i] : L_t  \geq r^{*}_t \ra \cap \la  L_{t_{i-1} } \leq r^{*}_{t_{i-1}} \ra  \rb \\
       &\geq 1 -\delta - \sum_{i=1}^{\infty} \delta_i \geq 1- 2\delta.
   \end{align*}
     
\noindent \textit{ \underline{Step 2c: Width of the confidence sequence.} } We will show that there exists a constant $M$ not depending on $t$ such that
\begin{align} \label{M}
r^{*}_t \leq M \cdot \frac{ \log(\delta^{-1}) + \log\log(t+10)}{t+10} 
\end{align}
for all $t\geq 0$, where $r^{*}$ is defined as in \eqref{r*}. 

Let us first show that 
\begin{align*}
    \sup_{k \geq 0} \la  h_k  \cdot \frac{t_k+10}{\log(\delta^{-1}) + \log\log(t_k+10)} \ra  < \infty.
\end{align*}
From \eqref{t-constraint}, we have 
\[
\frac{\log(16)}{\log \lb  \frac{1}{1-C_1 \eta_i} \rb} \geq t_i - t_{i-1} \geq \frac{\log(8)}{\log \lb  \frac{1}{1-C_1 \eta_i} \rb}.
\]
Note that there exists a universal constant $c_0 >0$ small enough such that 
\[
  \frac{c_0}{x} \leq \frac{1}{\log (1+x)} \leq \frac{1}{c_0 x}
\]
for all $x \in (0,1/15)$. 

Write 
\begin{align}
    t_k = \sum_{i=1}^{k} t_i - t_{i-1} &\leq \sum_{i=1}^{k}  \frac{\log(16)}{\log \lb  \frac{1}{1-C_1 \eta_i} \rb}  \\
    &=  \sum_{i=1}^{k} \frac{\log(16)}{\log \lb  1+ \frac{C_1 \eta_i}{1-C_1 \eta_i} \rb} \nonumber \\
    &\leq  \log(16) \cdot \sum_{i=1}^{k} \frac{1}{c_0} \cdot \frac{1-C_1\eta_i}{C_1\eta_i} \nonumber  \\
    &\leq \frac{\log(16)}{c_0 C_1 \kappa} \cdot   \sum_{i=1}^{k} h_{i-1}^{-1} \lb \log \lb \delta^{-1}  \rb + 2 \log(i+1) \rb \nonumber \\
    &= \frac{\log(16)h_0}{2c_0 C_1 \kappa} \cdot \sum_{i=1}^{k} 2^{i} \lb \log \lb \delta^{-1}  \rb + 2 \log(i+1) \rb \nonumber \\
    &\leq \frac{\log(16)h_0}{c_0 C_1 \kappa} \cdot 2^{k} \lb \log \lb \delta^{-1}  \rb + 2 \log(k+1) \rb. \label{tk upper}
\end{align}
Similarly,
\begin{align}
t_k &\geq \log(16)c_0 \cdot \sum_{i=1}^{k} \frac{1-C_1\eta_i}{C_1\eta_i} \nonumber \\
&\geq \frac{15 \log(16)c_0}{16} \cdot \sum_{i=1}^{k} h_{i-1}^{-1} \lb \log \lb \delta^{-1}  \rb + 2 \log(i+1) \rb \nonumber \\
&\geq  \frac{15 \log(16)c_0 h_0}{32} \cdot  2^{k} \lb \log \lb \delta^{-1}  \rb + 2 \log(k+1) \rb. \label{tk lower}
\end{align}
Now, observe that
\begin{align*}
      h_k  \cdot \frac{t_k+10}{\log(\delta^{-1}) + \log\log(t_k+10)} &\leq \frac{2^{-k} \cdot h_0 \left[ \frac{\log(16)h_0}{c_0 C_1 \kappa} \cdot 2^{k} \lb \log \lb \delta^{-1}  \rb + 2 \log(k+1) \rb \right]}{\log(\delta^{-1}) + \log\log(t_k+10)} \\
      &\leq \frac{\log(16)h_0^2}{c_0 C_1 \kappa}\cdot  \frac{  \log \lb \delta^{-1}  \rb + 2 \log(k+1)}{ \log(\delta^{-1}) + \log\log(t_k+10)}.
\end{align*}
If $k+1 \geq \sqrt{\delta} \exp  \lb \frac{2}{ \log(16)c_0h_0  }  \rb$, we have 
\begin{align*}
    &\log(\delta^{-1}) + \log\log(t_k+10)\\
    > &     \log(\delta^{-1}) + \log\log(t_k) \\
    \geq & \log(\delta^{-1}) + \log \left[  (k+1) \log 2 + \underbrace{\log \lb \frac{15 \log(16)c_0 h_0}{64} \lb \log \lb \delta^{-1}  \rb + 2 \log(k+1) \rb \rb}_{\geq 0} \right] \\
    \geq &  \log(\delta^{-1}) + \log (k+1) + \log \log 2.  
\end{align*}
Therefore,
\begin{align*}
          h_k  \cdot \frac{t_k+10}{\log(\delta^{-1}) + \log\log(t_k+10)} &\leq \frac{\log(16)h_0^2}{c_0 C_1 \kappa} \cdot  \frac{  \log \lb \delta^{-1}  \rb + 2 \log(k+1)}{ \log(\delta^{-1}) + \log (k+1) + \log\log 2} \\
          \leq \frac{\log(16)h_0^2}{2c_0 C_1 \kappa} 
\end{align*}
for all $k \geq \sqrt{\delta} \exp  \lb \frac{2}{ \log(16)c_0h_0  }  \rb -1$. 

Consequently,
\[
    \sup_{k \geq 0} \la  h_k  \cdot \frac{t_k+10}{\log(\delta^{-1}) + \log\log(t_k+10)} \ra = H < \infty
\]
for some $H>0$ that depends only on $h_0$ and the parameters in \eqref{cond-error} and \eqref{as-bound}. 

Now, suppose $t \in [t_{k-1}+1, t_k]$ for some $k\geq 1$. Recall $r^*_t$ in \eqref{r*}, we have
\begin{align*}
    r^{*}_t &= 4 \lb 1  - C_1 \eta_k \rb^{t-t_{k-1}} h_i \\
    &\leq 4h_i \\
    &\leq 4H \cdot \frac{\log \lb  \delta^{-1} \rb + \log \log \lb  t_k +10 \rb}{t_k + 10} \\
    &\leq 4H \cdot  \frac{\log \lb  \delta^{-1} \rb + \log \log \lb  t +10 \rb}{t + 10}
\end{align*}
by using Lemma \ref{decreasing}. Thus, the constant $M$ in \eqref{M} can be chosen to be $4H$.

\noindent{ \underline{\textit{Final step: Showing $\eta_t \asymp 1/t$.}} } Finally, we prove the statement regarding the asymptotic scaling of $\eta_t$. For $t \in [t_{k-1}+1,t_k]$, write
\begin{align*}
    \eta_t &= \frac{2\kappa h_0}{2^{k} \left[ \log \lb  \delta^{-1} \rb + 2  \log(k+1) \right]} \\
    & \leq \frac{2\kappa h_0}{\log (16) h_0} \cdot c_0C\kappa \cdot \frac{1}{t_k}  \\
    &\leq \frac{2c_0C \kappa^2}{\log(16)} \cdot \frac{1}{t},
    \end{align*}
where we use \eqref{tk upper} in the second line. For the lower bound, suppose $k \geq 3$. Then, by using \eqref{tk lower},
\begin{align*}
    \eta_t &= \frac{2\kappa h_0}{2^{k} \left[ \log \lb  \delta^{-1} \rb + 2  \log(k+1) \right]} \\
    & = \frac{\kappa h_0 }{2^{k-1} \left[ \log \lb  \delta^{-1} \rb + 2  \log k \right]} \cdot \underbrace{\frac{ \log \lb  \delta^{-1} \rb + 2  \log k }{ \log \lb  \delta^{-1} \rb + 2  \log(k+1)}}_{\geq 1/4} \\
    &\geq \frac{32 \kappa h_0}{15 \log(16) c_0 h_0} \cdot \frac{1}{t_{k-1}} \cdot \frac{1}{4} \\
    &\geq \frac{8 \kappa}{15 \log(16) c_0} \cdot \frac{1}{t}.
\end{align*}
The proof is completed. $\hfill$ $\square$



\section{Proof of Proposition \ref{pre-planned 1}} Rewrite \eqref{quantitative R-S lemma} as
\begin{align*}
    L_t \leq b_t \lb X_{t_0} + U^{*}_t  \rb
\end{align*}
where 
\begin{align*}
    b_t :&= \frac{(t_0+L-1)(t_0+L)}{(t_1+L-1)(t_1+L)} \\
    U^{*}_t :&= \frac{U_t}{b_t} + \frac{U_{t-1}}{b_{t-1}} + \dots + \frac{U_1}{b_1} =  \frac{U_t}{b_t} + U^{*}_{t-1} . 
\end{align*}
Let us give a tail bound on $\max_{t_0+1 \leq t  \leq  t_1}  U^*_k $.   Choose a time-dependent stopping rule $\tau$:
\[
\tau = \inf \la  t \geq t_0: X_t > r_t \ra
\]
for some sequence $\la r_t; t_0\leq t \leq t_1 \ra$ to be chosen later. 

Consider the stopped process 
\[
\max_{t_0+1 \leq t \leq t_1}  U^*_{\tau \wedge t}.
\]
Easily,
\begin{align*}
   | U^*_{\tau \wedge t} - U^*_{\tau \wedge (t-1)}| &= \frac{ \mathbf{1}_{\la  \tau \geq t \ra}}{b_t} \cdot \left| U_t \right| \\
   &\leq \frac{(t+L-1)(t+L)}{(t_0+L-1)(t_0+L)} \cdot \lb  \frac{2C_3\sqrt{r_{t-1}}}{C_1(t+L)} + \frac{4C_2}{C_1^2(t+L)^2}     \rb \\
   &= \underbrace{\frac{2C_3(t+L-1)\sqrt{r_{t-1}}}{C_1 (t_0+L-1)(t_0+L) } + \frac{4C_2(t+L-1)}{C_1^2(t+L)(t_0+L-1)(t_0+L)}}_{c_t}
\end{align*}
and
\begin{align*}
   \left| \mb{E} \left[ U^*_{\tau \wedge t} - U^*_{\tau \wedge (t-1)} | \mathcal{F}_{t-1}  \right] \right| &\leq \frac{4 C_2}{C_1^2 \lb t  +L \rb^2} \cdot \frac{(t+L-1)(t+L)}{(t_0+L-1)(t_0+L)}  \\
   &= \underbrace{ \frac{4C_2(t+L-1)}{C_1^2(t+L)(t_0+L-1)(t_0+L)}}_{\mu_t}.
\end{align*}
By using Azuma-Hoeffding maximal inequality, we have
\[
\mb{P} \lb  \max_{ t_0+ 1 \leq t \leq t_1}  U^*_{\tau \wedge t} \geq D \lb \delta \rb \Big| A \rb \leq \frac{\delta}{\mb{P} \lb A \rb}
\]
where
\begin{align*}
D(\delta) :&= \frac{4C_2(t_1-t_0)}{C_1^2(t_0+L-1)(t_0+L)} \\
&+ \sqrt{ 4\log \lb \delta^{-1} \rb \cdot \lb \frac{16C_2^2(t_1-t_0)}{C_1^4(t_0+L-1)^2(t_0+L)^2} + \sum_{t=t_0+1}^{t_1} \frac{4C_3^2(t+L-1)^2 r_{t-1}}{C_1^2(t_0+L-1)^2(t_0+L)^2}  \rb }.
\end{align*}
Now, let us tune the sequence $\la r_t \ra$ in a way such that 
\[
\mb{P} \lb \max_{t_0+1 \leq t \leq t_1}  U^*_{t} \geq D \lb \delta \rb \Big| A \rb = \mb{P} \lb  \max_{t_0+1 \leq t \leq t_1}  U^*_{\tau \wedge t} \geq D \lb \delta \rb \Big| A \rb.
\]
By Lemma \ref{stopped}, the display above holds if for any $1 \leq t \leq T$,
\[
\mb{P} \lb \tau \geq t+1 \Big| \max_{1 \leq k \leq t} U^{*}_k \leq D(\delta), A \rb =1.
\]
To get this, it suffices to choose $r_k$'s in a way such that 
\[
b_t \lb a_0 + D \lb \delta \rb \rb \leq r_t
\]
for all $ t_0 \leq t \leq t_1$.

For $t \in [t_0+1,t_1]$, let us take 
$$r_t=\frac{M \lb t_1-t_0 \rb \log \lb \delta^{-1} \rb}{(t+L)^2}$$ 
where the number $M$ is chosen to be  large enough and will be specified shortly. For any $t \in [t_0+1,t_1]$, we have
\begin{align*}
   & \frac{ (t_0+L-1)(t_0+L)}{(t+L-1)(t+L)} \lb a_0 + D(\delta) \rb  \\
    \leq & \frac{(t_0+L-1)(t_0+L)}{(t+L-1)(t+L)} \lb a_0 + \frac{4C_2(t_1-t_0)}{C_1^2(t_0+L-1)(t_0+L)} \rb \\
    + & 2\sqrt{\log \lb \delta^{-1} \rb } \cdot \frac{(t_0+L-1)(t_0+L)}{(t+L-1)(t+L)} \lb \frac{4C_2\sqrt{t_1-t_0}}{C_1^2(t_0+L-1)(t_0+L)} + \frac{2C_3\sqrt{M \log \lb \delta^{-1} \rb} (t_1-t_0)}{C_1(t_0+L-1)(t_0+L)} \rb
\end{align*}
To make sure that the display above is smaller than $r_t=\frac{M \lb t_1-t_0 \rb \log \lb \delta^{-1} \rb}{(t+L)^2}$, we need
\begin{align*}
 & a_0(t_0+L-1)(t_0+L) + \frac{4C_2(t_1-t_0)}{C_1^2} +  \frac{8C_2\sqrt{(t_1-t_0)\log \lb \delta^{-1} \rb}}{C_1^2} \\
 + & \frac{4C_3\sqrt{M} \log \lb \delta^{-1} \rb (t_1-t_0) }{C_1}  
\leq  \frac{t+L-1}{t+L} \cdot M \lb t_1 -t_0 \rb \log \lb \delta^{-1} \rb
\end{align*}
for all $t \in [t_0+1,t_1]$.

By observing that $(t+L-1)/(t+L) \geq (L-1)/L$ for all $t \geq 0$, it suffices to pick $M$ such that 
\begin{align*}
   x \cdot \frac{L}{L-1} M \lb t_1 -t_0 \rb \log \lb \delta^{-1} \rb &\geq a_0(t_0+L-1)(t_0+L); \\
    y \cdot  \frac{L}{L-1} M \lb t_1 -t_0 \rb \log \lb \delta^{-1} \rb &\geq  \max \la   
    \frac{4C_2(t_1-t_0)}{C_1^2}; \frac{8C_2\sqrt{(t_1-t_0)\log \lb \delta^{-1} \rb}}{C_1^2} \ra; \\
      z \cdot  \frac{L}{L-1} M \lb t_1 -t_0 \rb \log \lb \delta^{-1} \rb &\geq   \frac{4C_3\sqrt{M} (t_1-t_0)\log \lb \delta^{-1} \rb}{C_1}.
\end{align*}
for some positive numbers $x,y,z$ such that $x+y+z=1$.  Thus, we can pick $M$ as
\begin{align*}
M &= \frac{L-1}{L} \cdot \min_{\substack{x, y, z > 0 \\ x + y + z = 1}}  \la \max\la \frac{1}{x}; \frac{8}{y}; \frac{16}{z^2} \ra \ra \\ 
&\times \max \la \frac{ a_0 (t_0+L-1)(t_0+L)}{ \log \lb \delta^{-1} \rb(t_1-t_0)}; \frac{C_2}{C_1^2 \log \lb \delta^{-1} \rb}; \frac{C_2}{C_1^2\sqrt{\log \lb \delta^{-1} \rb}}; \frac{C_3^2}{C_1^2}  \ra.
\end{align*}
By Lemma \ref{19}, we have
\[
\min_{\substack{x, y, z > 0 \\ x + y + z = 1}}  \la \max\la \frac{1}{x}; \frac{8}{y}; \frac{16}{z^2} \ra \ra \leq 31.5.
\]
Thus, we can pick
\[
M = \frac{31.5(L-1)}{L} \cdot \max \la \frac{ a_0 (t_0+L-1)(t+L)}{ \log \lb \delta^{-1} \rb(t_1-t_0)}; \frac{C_2}{C_1^2 \log \lb \delta^{-1} \rb}; \frac{C_2}{C_1^2\sqrt{\log \lb \delta^{-1} \rb}}; \frac{C_3^2}{C_1^2}  \ra.
\]
The proof is completed. $\hfill$ $\square$

\section{Proof of Theorem \ref{conf} }

Let us pick $t_k=\alpha^{k+1}-\alpha$, where $\alpha=\max \la L-1, 33  \ra$. We will apply Lemma \ref{pre-planned 1} repeatedly on the intervals $[t_{k-1}+1, t_k]$, for $k \geq 1$ and $\delta_i=\delta/(i+1)^2$.  With a specified $a_0$, define 
\begin{align*}
    M_k :&= 31.5 \times\max \la \frac{ a_{k-1} (t_{k-1}+L-1)(t_{k-1}+L)}{ \log \lb \delta_k^{-1} \rb(t_{k}-t_{k-1})}; \frac{C_2}{C_1^2};  \frac{C_3^2}{C_1^2}  \ra ;  \\
    r(k,t) :&= \frac{M_k \lb t_k-t_{k-1} \rb \log \lb \delta_k^{-1} \rb}{(t+L)^2} \ \text{for} \ t \in [t_{k-1}+1,t_k];   \\
    a_{k+1} :&= r(k+1,t_{k+1}).
\end{align*}
It is useful to think of this as $a \to M \to r \to a $. By Lemma \ref{pre-planned 1}, for any $k\geq 1$, we have
\[
\mb{P} \lb  \exists t \in [t_{k-1}+1, t_k]: L_t \geq r \lb k,t \rb, L_{t_{k-1}} \leq a_{k-1}  \rb \leq \delta_k.
\]
Therefore, by using the union bound,
\[
\mb{P} \lb  \forall k \geq 1: \max_{t \in [t_{k-1}+1, t_k]} \frac{L_t}{r(k,t)} \leq 1 \rb \geq 1 -\sum_{k\geq 0} \delta_i \geq 1-2\delta.
\]
Next, let us show that by induction on $k$ such that 
\[
a_k \leq A \cdot\frac{ \log \lb  \delta_k^{-1}  \rb}{t_k+L}
\]
where
\[
A:=  31.5 \times \max \la \frac{a_0L}{\log \lb \delta^{-1} \rb}; \frac{C_2}{C_1^2}; \frac{C_3^2}{C_1^2} \ra.
\]
For $k=0$, the statement above is true since 
\[
a_0 \leq  \frac{31.5 \times a_0L}{\log \lb \delta^{-1} \rb}  \cdot\frac{\log \lb  \delta_0^{-1}  \rb}{t_0+L} = 31.5 \times a_0.
\]
Assume the statement is true for $k$, let us consider $k+1$.  Note that 
\begin{align*}
    a_{k+1} &\leq A \cdot\frac{\log \lb  \delta_{k+1}^{-1}  \rb}{t_{k+1}+L} \\
 \iff  \frac{M_{k+1} \lb t_{k+1}-t_{k} \rb \log \lb \delta_{k+1}^{-1} \rb}{(t_{k+1}+L)^2} &\leq A \cdot\frac{\log \lb  \delta_{k+1}^{-1}  \rb}{t_{k+1}+L} \\
 \iff  M_{k+1} \cdot  \frac{t_{k+1}-t_k}{t_{k+1}+L} &\leq A \\
 \iff  \frac{ 31.5 \times a_{k} (t_{k}+L-1)(t_{k}+L)}{\log \lb \delta_{k+1}^{-1} \rb(t_{k+1}-t_{k})} \cdot  \frac{t_{k+1}-t_k}{t_{k+1}+L}  &\leq A.
\end{align*}
By the induction hypothesis,
\begin{align*}
\frac{ 31.5 \times a_{k} (t_{k}+L-1)(t_{k}+L)}{ \log \lb \delta_{k+1}^{-1} \rb(t_{k+1}-t_{k})} \cdot  \frac{t_{k+1}-t_k}{t_{k+1}+L} &=  \frac{ 31.5 \times a_{k} (t_{k}+L-1)(t_{k}+L)}{ \log \lb \delta_{k+1}^{-1} \rb (t_{k+1}+L)} \\
       &\leq 31.5 \times \cdot \frac{A\log \lb \delta_{k}^{-1} \rb}{\log \lb \delta_{k+1}^{-1} \rb} \cdot \frac{(t_{k}+L-1)(t_{k}+L)}{(t_{k+1}+L)(t_k+L)} \\
      &\leq 31.5 \times A \cdot  \frac{1}{\alpha} \leq A.
\end{align*}
Note that this also gives
\begin{align*}
\sup_{k \geq 1} M_k &\leq 31.5 \times \max \la A \cdot \sup_{k \geq 1} \la \frac{\log \lb  \delta_{k-1}^{-1}  \rb}{ \log \lb  \delta_k^{-1}  \rb} \cdot \frac{\lb t_{k-1}+ L-1 \rb \lb t_{k-1}+ L \rb}{(\alpha-1)\lb t_{k-1}+ L \rb\lb t_{k-1} +\alpha  \rb} \ra; \frac{C_2}{C_1^2}; \frac{C_3^2}{C_1^2} \ra \\
& \leq A
\end{align*}
where we use Lemma \ref{tk technical} to get 
\[
\frac{t_{k-1}+ L-1 }{(\alpha-1)\lb t_{k-1} +\alpha  \rb} \leq \frac{1}{31.5}
\]
in the last line. To get the final convergence rate, note that for $t \in [t_{k-1}+1,t_k]$,
\begin{align*}
    r(k,t) &\leq A \cdot \frac{(\alpha-1)\lb t_{k-1}+\alpha\rb \left[ \log \lb \delta^{-1} \rb + 2\log (k+1) \right]}{(t+L)^2} \\
    &\leq A\cdot (\alpha-1) \cdot \lb \mathbf{1}_{\la L \geq 32 \ra} + 33 \cdot \mathbf{1}_{\la L\leq 31 \ra} \rb \cdot \frac{\log \lb \delta^{-1} \rb + 2\log (k+1)}{t+L}.
\end{align*}
Moreover, for all $k\geq 1$,
\begin{align*}
    \log (k+1) &\leq \log \log \lb 33^{k}-33+ 10 \rb \\
    &= \log \log \lb t_{k-1} + 10 \rb \leq \log \log \lb t+9 \rb.
\end{align*}
Thus, for $t \notin \la t_k; k\geq 0 \ra$, we have
\[
r(k,t) \leq KA \cdot \frac{\log \lb \delta^{-1} \rb + 2\log \log (t+9)}{t+L}.
\]
where
\[
K:=  (\alpha-1) \cdot \lb \mathbf{1}_{\la L \geq 32 \ra} + 33 \cdot \mathbf{1}_{\la L\leq 31 \ra} \rb.
\]
Finally, at one of the epochs $t_k$ ($k\geq 1$), we have
\begin{align*}
 a_k &\leq  A \cdot  \frac{\log \lb \delta^{-1} \rb + 2\log (k+1)}{t_k+L} \\
 &\leq A \cdot  \frac{\log \lb \delta^{-1} \rb + 2\log \log (t_k+9)}{t_k+L}.
\end{align*}
Consequently,
\[
\mb{P} \lb \exists t \geq 0: L_t \geq 31.5 \times K \max \la \frac{a_0L}{\log \lb \delta^{-1} \rb}; \frac{C_2}{C_1^2}; \frac{C_3^2}{C_1^2} \ra \frac{\log \lb \delta^{-1} \rb + 2\log \log (t+9)}{t+L} \rb \geq 1-2\delta.
\]
The proof is completed. $\hfill$ $\square$

\section{Proof of Theorem \ref{oja-conf}}

In what follows, $S_t$ stands for $L^{Oja}_t$ for simplicity.  We split the proof of Theorem \ref{oja-conf} into two steps.

\noindent \underline{\textit{ Step 1: Confidence sequence for a modified process}}. Define the process
\begin{align*}
        S^{*}_t :&=  S_0;  \\
    S^{*}_t :&= S_t \cdot    \prod_{i=0}^{t-1} \mathbf{1}_{\la S_{i} \leq 1/2   \ra}.
\end{align*}
Observe that for any $t\geq 1$,
\begin{align}
        S^{*}_t &\leq \prod_{i=0}^{t-1} \mathbf{1}_{\la S_{i} \leq 1/2   \ra} \cdot \left[ \lb  1-2 \rho \eta_t \rb S_{t-1} + 2\rho \eta_t S_{t-1}^2 + Q_t + 4B^4 \eta^2_t  \right] \nonumber \\
        &\leq \lb  1-2 \rho \eta_t \rb S^{*}_{t-1} \cdot  \mathbf{1}_{\la S_{t-1} \leq 1/2   \ra} + 2 \rho \eta_t \lb S^{*}_{t-1} \rb^2 \cdot  \mathbf{1}_{\la S_{t-1} \leq 1/2   \ra} \nonumber \\
        &+ Q_t \cdot  \prod_{i=0}^{t-1} \mathbf{1}_{\la S_{i} \leq 1/2   \ra} + 4B^4 \eta_t^2 \nonumber \\
        &\leq \lb  1-2 \rho \eta_t \rb S^{*}_{t-1} + \rho \eta_t S^{*}_{t-1} + Q^{*}_t +  4B^4 \eta_t^2 \nonumber \\
        &=  \lb  1-\rho \eta_t \rb S^{*}_{t-1} + Q^{*}_t +  4B^4 \eta_t^2 \nonumber
\end{align}
where
\[
 Q^{*}_t :=  Q_t \cdot  \prod_{i=0}^{t-1} \mathbf{1}_{\la S_{i} \leq 1/2   \ra}.
\]
Note that we still have 
\begin{align*}
    \mb{E} \lb Q^{*}_t | \mathcal{F}_{t-1}  \rb = 0 
\end{align*}
and 
\[
\left|  Q^{*}_t  \right| \leq \eta_t B^2 \sqrt{        S^{*}_{t-1}}. 
\]
Therefore, by applying Theorem \ref{conf} with the constants
    $C_1 \equiv \rho$, 
    $C_2 \equiv 4B^4$, 
    $C_3 \equiv  B^2$ and 
    $a \equiv 1/4$,
and the choice of step sizes $\eta_t=2/\rho(t+L)$, we obtain that
\begin{align} \label{conf S*}
\mb{P} \lb \forall t \geq 0:  S^{*}_t \leq 31.5 \times K  \cdot \max \la \frac{L}{4 \log \lb \delta^{-1} \rb}; \frac{B^4}{\rho^2}\ra \frac{\log \lb \delta^{-1} \rb + 2\log \log (t+9)}{t+L} \rb \geq 1-2\delta
\end{align}
where $K=\max \la L-2; 32 \ra  \cdot \lb \mathbf{1}_{\la L \geq 32 \ra} + 32 \cdot \mathbf{1}_{\la L\leq 31 \ra} \rb$. In other words, \eqref{conf S*} reduces to
\[
\mb{P} \lb \forall t \geq 0:  S^{*}_t \leq 1008  \cdot \max \la \frac{L}{4 \log \lb \delta^{-1} \rb}; \frac{B^4}{\rho^2}\ra \frac{\log \lb \delta^{-1} \rb + 2\log \log (t+9)}{t+L} \rb \geq 1-2\delta,
\]
for $L \geq 32$.
Note that the statement above holds because 
\[
\mb{P} \lb S^{*}_0 \geq 1/4 \rb = \mb{P} \lb S_0 \geq 1/4 \rb \geq 1- \delta^3 \geq 1-\delta.
\]

\noindent \underline{\textit{ Step 2: The original process is close to the modifed process}}.  Let us show that if one sets $$L= \left[128B^4 \log \lb \delta^{-1} \rb^2/ \rho^2 \right]$$
and $\eta_t=2/\rho(t+L)$, then 
\[
\mb{P} \lb \forall t \geq 0: S_t= S^{*}_t \rb \geq 1-2e\delta.
\]
By using Proposition \ref{oja-err bound}, for any $t\geq 1$,  we have that 
\begin{align*}
    \mb{E} \left[  \exp \lb \lambda S_t \rb | \mathcal{F}_{t-1} \right] &\leq \mb{E} \left[   \exp \lb \lambda S_{t-1} + \lambda Q_t + 4B^4 \lambda  \eta_t^2\rb  | \mathcal{F}_{t-1} \right] \\
    &\leq \exp \lb S_{t-1} + 4B^4 \eta_t^2 \rb  \mb{E} \left[ \exp \lb \lambda Q_t \rb | \mathcal{F}_{t-1} \right] \\
    &\leq  \exp \lb \lambda S_{t-1} + 4 \lambda B^4 \eta_t^2 \rb  \exp \lb 32 \lambda^2 B^4 \eta_t^2 \rb,
\end{align*}
where the last estimate follows from Chernoff's inequality. 

Thus, for any $\delta, \lambda >0$,
\[
    \mb{E} \left[  \exp \lb \lambda S_t \rb | \mathcal{F}_{t-1} \right] \leq  \exp \lb \lambda S_{t-1} + 4 \lambda B^4 \eta_t^2 + 32\lambda^2 B^4 \eta_t^2 \rb
\]
Consequently, for any $\lambda >0$, the process
\[
\la V_t \ra_{t\geq 0} := \la \exp \left[ \lambda S_{t} + \lb  4 \lambda B^4+ 32\lambda^2 B^4\rb \cdot \sum_{k \geq t+1}  \eta_k^2  \right] \ra_{t\geq 0}
\]
forms a non-negative super martingale. Here we define $V_0:= \lambda S_0$. 

Moreover, for any $t\geq 1$,
\begin{align*}
    \sum_{k \geq t+1}  \eta_k^2 \leq \sum_{k\geq1}  \eta_k^2 &= \sum_{k=1}^{\infty} \frac{4}{\rho^2(k+L)^2} \\
    &\leq \frac{4}{\rho^2} \int_{L+1}^{\infty} \frac{1}{x^2} dx = \frac{4}{\rho^2(L+1)}.
\end{align*}

Thus, Ville's inequality yields 
\begin{align*}
    \mb{P} \lb  \sup_{t \geq 0} S_t \geq 1/2 \rb &\leq  \inf_{\lambda>0} \la \mb{P} \lb \sup_{t\geq 0} V_t \geq \frac{\lambda}{2} - \frac{16\lambda B^4(1+8\lambda)}{\rho^2(L+1)} \rb \ra \\
    &\leq \inf_{\lambda > 0} \la \exp \lb -\frac{\lambda}{2} + \frac{16\lambda B^4(1+8\lambda)}{\rho^2(L+1)} \rb \mb{E} V_0 \ra \\
    &= \inf_{\lambda > 0} \la  \exp \lb -\frac{\lambda}{2} + \frac{32\lambda B^4(1+8\lambda)}{\rho^2(L+1)} \rb  \cdot \mb{E} \lb \exp \lb \lambda S_0 \rb \rb \ra.
\end{align*}
To bound the last term, note that
\begin{align*}
    \mb{E} \lb \exp \lb \lambda S_0 \rb \rb \leq e^{\lambda/4} \lb 1 - \delta^3 \rb + e^{\lambda} \cdot \delta^3.
\end{align*}
Consequently,
\begin{align*}
     \mb{P} \lb  \sup_{t \geq 0} S_t \geq 1/2 \rb &\leq \inf_{\lambda>0} \la \exp \lb  \frac{32\lambda B^4(1+8\lambda)}{\rho^2(L+1)} \rb \cdot \left[ e^{-\lambda/4} \lb 1 - \delta^3 \rb + e^{\lambda/2} \cdot \delta^3 \right]  \ra.
\end{align*}
By setting 
\[
\lambda= \frac{4}{3} \log \lb \frac{1-\delta^3}{\delta^3} \rb \leq 4\log \lb \delta^{-1} \rb,
\]
we get 
\begin{align*}
     \mb{P} \lb  \sup_{t \geq 0} S_t \geq 1/2 \rb &\leq \exp \lb 128 \frac{B^4 \log \lb \delta^{-1} \rb^2}{\rho^2(L+1)} \rb \cdot \left[  2\delta \cdot \lb 1 -\delta^3  \rb^{2/3}  \right] \\
     & \leq 2e\delta 
\end{align*}
whenever
\[
L+1 \geq \frac{128 B^4 \log \lb \delta^{-1} \rb^2}{\rho^2}.
\]
Thus,
\begin{align*}
    \mb{P} \lb \forall t \geq 0: S_t= S^{*}_t \rb &\geq 1- \mb{P} \lb  \sup_{t \geq 0} S_t \geq 1/2 \rb \\
    &\geq 1 -2e\delta.
\end{align*}
The proof is completed by combining the expression above and \eqref{conf S*}. $\hfill$ $\square$

\section{Proof of Theorem \ref{lower bound}}

In what follows, $S_t$ stands for $L^{RM}_t$ for simplicity. Let us split the proofs in a few steps. 

\noindent \underline{\textit{ Step 1: $X_n$ converges to $0$ almost surely}}. Recall that for $S_t=(X_t-\theta)^2$, Proposition \ref{robbins-monroe constraint} gives the recursion
\[
S_t \leq (1-2R \eta_t)S_{t-1} + Q_t + 2\eta_t^2 \left[ \lb A \sqrt{S_{t-1}}+B \rb^2 + R_1^2 \right].
\]
 By taking the conditional expectation over $\mc{F}_{t-1}$ and note that $\mb{E} \lb Q_t | \mc{F}_{t-1} \rb=0$, we have
\begin{align*}
    \mb{E} \lb S_{t} | \mc{F}_{t-1} \rb & \leq (1-2R \eta_t)S_{t-1} + 2\eta_t^2 \lb A^2 S_{t-1} + 2AB \sqrt{S_{t-1}} + B^2+ R_1^2 \rb \\
    &\leq (1-2R \eta_t)S_{t-1} + 2\eta_t^2 \lb A^2 S_{t-1} + A^2 + B^2 S_{t-1} + R_1^2  \rb \\
    &= S_{t-1} \lb 1 + 2\eta_t^2 (A^2+B^2) \rb - 2R\eta_tS_{t-1} + 2\eta_t^2 \lb A^2+R_1^2 \rb
\end{align*}
By applying the Robbins-Siegmund lemma in the form of \eqref{almost-supermartingale}, with  
\[
a_t = 2\eta_t^2 (A^2 + B^2), \quad
b_t = 2\eta_t^2 (A^2 + R_1^2), \quad
c_t = 2R\eta_t S_{t-1},
\]  
and noting that  
\begin{align*}
    \sum_{t \geq 1} a_t &\leq 2(A^2 + B^2)L_2^2 \sum_{t=1}^{\infty} \frac{1}{t^2} < \infty, \\
    \sum_{t \geq 1} b_t &\leq 2(A^2 + R_1^2)L_2^2 \sum_{t=1}^{\infty} \frac{1}{t^2} < \infty,
\end{align*}
we conclude that \( S_t \) converges almost surely to a non-negative random variable \( S \), and  
\[
\sum_{t \geq 1} c_t = \sum_{t = 1}^{\infty} \eta_t S_{t-1} < \infty
\]  
almost surely.

If $P(S>\ve)>\ve$ for some $\ve>0$, then 
\[
 \sum_{t=1}^{\infty} \eta_t S_{t-1}  \geq L_1 \sum_{t=1}^{\infty} \frac{S_{t-1}}{t} \geq \frac{L_1 \ve}{2} \cdot \limsup_{k \to \infty} \sum_{i=k}^{\infty} \frac{1}{k} =\infty
\]
on the event $\la S>\ve \ra$. Thus,
\[
\mb{P} \lb \sum_{t=1}^{\infty} \eta_t S_{t-1} =  \infty \rb \geq \mb{P} \lb S>0 \rb > \ve,
\]
which is a contradiction. Therefore, $X_t$ converges almost surely to $\theta$.

\noindent \underline{\textit{ Step 2: A moderate-deviation bound}}. For simplicity, write $Z_k= Y(X_k) - M(X_k)$. Note that $\la Z_k; k \geq 1\ra$ are i.i.d. with mean $0$ and unit variance. Moreover, for all $n \geq 1$, we have
\begin{align} \label{decomp}
    X_{2^{n+1}} = X_{2^n+1} - \sum_{k=2^n+1}^{2^{n+1}} \eta_k M(X_k) - \sum_{k=2^n+1}^{2^{n+1}} \eta_k Z_k.
\end{align}
We will show that 
\begin{align} \label{mod dev}
    \mb{P} \lb \sum_{k=2^n+1}^{2^{n+1}} \eta_k Z_k \geq M \sqrt{\frac{\log\log 2^n}{2^n}} \rb \geq \frac{1}{n}
\end{align}
for all $n$ sufficiently large, where 
\[
M:= \sqrt{L_1/2}.
\]
We will make use of the following moderate deviation result from \cite{chen2013stein} (Proposition 4.5 in the same paper for more details)
\begin{lemma}[\cite{chen2013stein}] \label{mod deviation}
    Let \(\{\xi_i: 1 \le i \le n\}\) be independent random variables satisfying \(\mathbb{E}\,\xi_i = 0\) and \(\mathbb{E}\,\mathrm{e}^{t_n|\xi_i|} < \infty\) for some \(t_n > 0\) and all \(1 \le i \le n\). Assume furthermore that
\begin{equation}\label{eq:xi-variance}
\sum_{i=1}^n \mathbb{E}\,\xi_i^2 = 1.
\end{equation}
Define
\[
W = \sum_{i=1}^n \xi_i,
\qquad
\gamma = \sum_{i=1}^n \mathbb{E}\!\bigl[\,|\xi_i|^3 \mathrm{e}^{x|\xi_i|}\bigr].
\]
Then, for \(0 \le x \le t_n\),
\[
\frac{\mathbb{P}(W \ge x)}{1 - \Phi(x)}
= 1 + O(1)\, (1 + x^3)\,\gamma\,\mathrm{e}^{4x^3 \gamma},
\]
where \(\Phi\) is the CDF of a standard Gaussian and the \(O(1)\) term is universal.
\end{lemma}
Put 
\begin{align*}
    \sigma_n^2 :&= \mbox{Var} \lb  \sum_{k=2^n+1}^{2^{n+1}} \eta_k Z_k \rb =   \sum_{k=2^n+1}^{2^{n+1}} \eta_k^2 .
\end{align*}
It is easy to check that
\[
  L_1 \sum_{k=2^n+1}^{2^{n+1}} \frac{1}{k^2} \leq \sigma_n^2 \leq L_2 \sum_{k=2^n+1}^{2^{n+1}} \frac{1}{k^2}.
\]
Thus,
\begin{align*}
\sigma_n^2 \geq L_1 \int_{2^n+1}^{2^{n+1}} \frac{1}{x^2} &= L_1 \lb \frac{1}{2^n+1} -\frac{1}{2^{n+1}} \rb \\
&= L_1 \cdot \frac{2^{n+1}-2^n-1}{2^{n+1} \lb 2^n+1     \rb} \\
&= L_1 \cdot \frac{2^n-1}{2^{n+1}(2^n+1)},
\end{align*}
and 
\begin{align*}
    \sigma_n^2 \leq L_2 \int_{2^n}^{2^{n+1}} \frac{1}{x^2} = L_2 \cdot \lb \frac{1}{2^n} -\frac{1}{2^{n+1}} \rb = \frac{L_2}{2^n}.
\end{align*}
Consequently, for all $n \geq 10$,
\[
\frac{L_1}{2^{n+1}} \leq \sigma_n^2 \leq \frac{L_2}{2^n}.
\]
We now employ Lemma \ref{mod deviation} to get 
\begin{align*}
    & \frac{\mb{P} \lb  \sum_{k=2^n+1}^{2^{n+1}}  \eta_k Z_k \geq M \sqrt{\frac{\log\log 2^n}{2^n}} \rb}{1- \Phi \lb  \frac{M\sqrt{\log \log 2^n}}{\sigma_n 2^{n/2}}\rb} \\
    = & 1 + O(1)\lb  1 +  x_n^3   \rb \gamma_n \exp \lb 4x_n^3\gamma_n \rb
\end{align*}
where
\begin{align*}
    x_n &:= \frac{M}{\sigma_n \cdot  2^{n/2}} \cdot \sqrt{\log\log 2^n}  ; \\
    \gamma_n &:= \sum_{k=2^n+1}^{2^{n+1}}  \frac{\eta_k^3}{\sigma_n^3} \cdot \mb{E} \lb |Z_k|^3 \cdot \exp\lb \frac{x_n \eta_kZ_k}{\sigma_n} \rb  \rb.
\end{align*}
It is easy to see that
\[
x_n \leq \frac{M\sqrt{2}}{\sqrt{L_1}} \sqrt{\log \lb n + \log 2 \rb}.
\]
Let us estimate the size of $\gamma_n$. Observe that for all $k \in \left[ 2^n+1, 2^{n+1} \right]$, we have
\[
\eta_k \leq \frac{L_2}{k} \leq \frac{L_2}{2^n+1}.
\]
Thus, for all $n$ sufficiently large,
\begin{align*}
    \gamma_n &\leq \sigma_n^{-3} \cdot \max_{2^n+1\leq k \leq 2^{n+1}} \la  \mb{E} \lb |Z_k|^3 \cdot \exp\lb \frac{x_n \eta_kZ_k}{\sigma_n} \rb  \rb  \ra \cdot  \sum_{k=2^n+1}^{2^{n+1}} \frac{L_2^3}{k^3} \\ 
    &\leq \frac{2^{3/2(n+1)}}{L_1^{3/2}} \cdot  \max_{2^n+1\leq k \leq 2^{n+1}} \la R_1^3 \cdot \underbrace{\exp \lb  \frac{M\sqrt{2}}{\sqrt{L_1}} \sqrt{\log \lb n + \log 2 \rb} \cdot \frac{L_2 R_1}{2^n+1} \cdot \sqrt{\frac{2^{n+1}}{L_1}} \rb }_{\leq 1 \ \text{when $n$ is large}}   \ra \\
    &\times L_2^3 \cdot  \sum_{k=2^n+1}^{2^{n+1}} \frac{1}{k^3} \\
    &\leq \frac{2^{3/2(n+1)} L_2^3}{L_1^{3/2}} \cdot \lb  R_1^3 e \rb \cdot \int_{2^n}^{2^{n+1}} \frac{1}{x^3} dx
    = \frac{3 L_2^3 R_1^3 e}{2L_1^{3/2}} \cdot \frac{2^{3/2(n+1)}}{2^{2n}} = \frac{3\sqrt{2}L_2^3 R_1^3 e }{L_1^{3/2}}
    \cdot 2^{-n/2}.
\end{align*}
Therefore, for all $n$ large enough, we have
\begin{align*}
    \frac{ \mb{P} \lb  \sum_{k=2^n+1}^{2^{n+1}}  \eta_k Z_k \geq M \sqrt{\frac{\log\log 2^n}{2^n}} \rb}{1- \Phi \lb  \frac{M\sqrt{\log \log 2^n}}{\sigma_n 2^{n/2}}\rb} \geq 1/2.
\end{align*}
By using the tail bound $1 - \Phi(x) \geq \lb 2 \pi \rb^{-1/2} \cdot (2x)^{-1} \cdot e^{-x^2/2}$ (which holds for all $x$ suffciently large), we obtain 
\begin{align*}
    \mb{P} \lb  \sum_{k=2^n+1}^{2^{n+1}}  \eta_k Z_k \geq M \sqrt{\frac{\log\log 2^n}{2^n}} \rb &\geq \frac{1}{4\sqrt{2\pi}} \cdot \frac{\sigma_n \cdot 2^{n/2}}{M \sqrt{\log \log 2^n}} \cdot \exp \lb -\frac{M^2 \log \log 2^n}{2 \sigma_n^2 2^n} \rb \\
    &\geq \frac{\sqrt{L_1}}{8\sqrt{\pi}} \cdot \frac{1}{M\sqrt{\log \lb  n + \log 2 \rb}} \cdot (n+\log 2)^{- \frac{M^2}{L_1}} \\
    &=  \frac{\sqrt{L_1}}{8\sqrt{\pi}} \cdot \frac{1}{M\sqrt{\log \lb  n + \log 2 \rb}} \cdot \lb  n+\log 2 \rb^{-1/2} \\
    &\geq \frac{1}{n}
\end{align*}
for all $n$ sufficiently large where we have used the fact that $M=\sqrt{L_1/2}$ in the third display.

\noindent \underline{\textit{ Step 3: Wrap-up}}. From \eqref{mod dev}, we conclude that
\[
\sum_{n=1}^{\infty}     \mb{P} \lb  \sum_{k=2^n+1}^{2^{n+1}}  \eta_k Z_k \geq M \sqrt{\frac{\log\log 2^n}{2^n}} \rb \geq \liminf_{k \to \infty} \lb \sum_{i=k}^{\infty} \frac{1}{i} \rb = \infty.
\]
Note that the events $\la    \sum_{k=2^n+1}^{2^{n+1}}  \eta_k Z_k \geq M \sqrt{\frac{\log\log 2^n}{2^n}}  \ra$ are independent, so Borel–Cantelli lemma yields 
\begin{align} \label{Borel-Cantelli}
  \sum_{k=2^n+1}^{2^{n+1}}  \eta_k Z_k \geq \sqrt{\frac{L_1}{2}} \cdot \sqrt{\frac{\log\log 2^n}{2^n}}
\end{align}
eventually always with probability one. 

Recall $L$ from \eqref{L}. Define the event 
\[
A= \la  \limsup_{t \to \infty} \lb \sqrt{\frac{t}{\log \log t}}\cdot  |X_t-\theta| \rb \leq L  \ra.
\]
If $\mb{P}(A)=0$, then the proof is completed. Assume $\mb{P}(A)>0$. From Step 1, we know that $X_t$ converges almost surely to $0$, which leads to 
\begin{align*}
    |M(X_k)| = |M(X_k) - M(\theta)| \leq 2M'(\theta) \cdot |X_k - \theta| 
\end{align*}
for all $k$ sufficiently large, with probability one. 

Now, on the event $A$, we have
\begin{align*}
   & \limsup_{n \to \infty} \la \sqrt{\frac{2^n}{\log \log 2^n}} \cdot \left| X_{2^n+1} - X_{2^{n+1}}  \right| \ra \\
    \leq &  \limsup_{n \to \infty} \la   \sqrt{\frac{2^n}{\log \log 2^n}} \cdot \left| X_{2^n+1} - \theta  \right| +  \sqrt{\frac{2^n}{\log \log 2^n}} \cdot \left| X_{2^{n+1}} - \theta  \right|  \ra \\
    \leq &   L   \cdot \limsup_{n \to \infty} \la     \sqrt{\frac{2^n}{\log \log 2^n}} \cdot  \lb \sqrt{ \frac{\log \log \lb 2^ n +1  \rb}{2^n + 1}} + \sqrt{\frac{\log \log \lb 2^{n+1}  \rb}{2^{n+1}}}  \rb \ra \\
    \leq &  L   \cdot \lb 1 + \frac{1}{\sqrt{2}} \rb < 2L.
\end{align*}
Moreover, from \eqref{decomp}, we have 
\begin{align*}
    & \sqrt{\frac{2^n}{\log \log 2^n}} \cdot \left|  \sum_{k=2^n+1}^{2^{n+1}}  \eta_k Z_k \right| \\
    = & \sqrt{\frac{2^n}{\log \log 2^n}} \cdot \left| X_{2^{n+1}} - X_{2^n+1} + \sum_{k=2^n+1}^{2^{n+1}}  \eta_k M(X_k) \right| \\
    \leq &  \sqrt{\frac{2^n}{\log \log 2^n}} \cdot \left| X_{2^n+1} - X_{2^{n+1}}  \right| + \sqrt{\frac{2^n}{\log \log 2^n}} \cdot \left| \sum_{k=2^n+1}^{2^{n+1}}  \eta_k M(X_k)  \right|,
\end{align*}
which leads to
\begin{align*}
   & \limsup_{n \to \infty} \la  \sqrt{\frac{2^n}{\log \log 2^n}} \cdot \left|  \sum_{k=2^n+1}^{2^{n+1}}  \eta_k Z_k \right| \ra \\
     \leq & \limsup_{n \to \infty} \la  \sqrt{\frac{2^n}{\log \log 2^n}} \cdot \left| X_{2^n+1} - X_{2^{n+1}}  \right| \ra \\
     + & \limsup_{n \to \infty} \la  \sqrt{\frac{2^n}{\log \log 2^n}} \cdot 2L_2 M^{'}(\theta) \cdot 2L \cdot \sum_{k=2^n+1}^{2^{n+1}} \sqrt{\frac{\log \log k}{k}}  \cdot \frac{1}{k}\ra.
\end{align*}
On the event $A$, by noting that the function $x \to \log \log(x)/x$ is decreasing, we get
\begin{align*}
     & \limsup_{n \to \infty} \la  \sqrt{\frac{2^n}{\log \log 2^n}} \cdot \left|  \sum_{k=2^n+1}^{2^{n+1}}  \eta_k Z_k \right| \ra \\
     \leq & 2L + 4LL_2M^{'}(\theta) \sum_{k=2^n+1}^{2^{n+1}}  \frac{1}{k} \\
     \leq & 2L \lb 1 + 2L_2M^{'}(\theta) \log 2 \rb = \frac{\sqrt{L_1}}{2} .
\end{align*}
Therefore, on the event $A$,
\[
\sqrt{\frac{2^n}{\log \log 2^n}} \cdot   \sum_{k=2^n+1}^{2^{n+1}}  \eta_k Z_k \leq \frac{\sqrt{L_1}}{2} 
\]
eventually always, which is a contradiction to \eqref{Borel-Cantelli}. The proof is completed. $\hfill$ $\square$

\begin{appendix}
\section{Extension to the sub-Gaussian coefficients settings} \label{subG coeff}
\subsection{Extension of Theorem \ref{master}}
In this section, we extend Theorem \ref{master} to the settings where the coefficients $A_i, B_i$ in \eqref{cond-error} and \eqref{as-bound} are allowed to be random, positive random variables. The only coefficient that has to be deterministic is $C_1$ in \eqref{quantitative R-S lemma}. It turns out that in such settings, the conclusion of Theorem \ref{master} is still valid under a sub-Gaussian moment condition on these coefficients.

Let us  first restate $\eqref{cond-error}$ and \eqref{as-bound} as  
\begin{align}
          &  \Big| \mb{E} \lb U_t | \mathcal{F}_{t-1}  \rb \Big| \leq  \mb{E} \lb \mc{C}^{(t)}_2 | \mc{F}_{t-1} \rb \cdot \eta_t^2 + \sum_{i=1}^{m} \mb{E} \lb \mc{A}^{(t)}_{i} | \mc{F}_{t-1} \rb \cdot \eta_t^{1+a_i} L_{t-1}^{b_i};    \label{cond-error2} \\
       & |U_t| \leq \mc{C}^{(t)}_3 \eta_t \sqrt{L_{t-1}} + \sum_{i=1}^{m} \mc{B}^{(t)}_{i} \cdot \eta_t^{1/2+c_i} \cdot L_{t-1}^{d_i}.  \label{as-bound2}
\end{align}
where $\la  \mc{C}^{(t)}_2;  \mc{C}^{(t)}_3; \mc{A}^{(t)}_{i}; \mc{B}^{(t)}_{i}; 1 \leq i \leq m; t \geq 1 \ra$ are positive random variables. 

\begin{assump} \label{SE Bi}
 There exists $\sigma>0$ such that

\noindent \textit{(i)}. 
\begin{align} \label{Sub expo Bi}
    \max \la \mb{E} \left[ \exp \lb \frac{|\mc{C}^{(t)}_3|^2}{\sigma^2} \rb \Big| \mc{F}_{t-1}  \right]; \left[ \exp \lb \frac{|\mc{B}^{(t)}_1|^2}{\sigma^2} \rb \Big| \mc{F}_{t-1}  \right];\dots; \left[ \exp \lb \frac{|\mc{B}^{(t)}_m|^2}{\sigma^2} \rb \Big| \mc{F}_{t-1}  \right] \ra   \leq 2
\end{align}
almost surely.

\noindent \textit{(ii)}. We have 
\begin{align} \label{bounded Ai}
    \max_{t \geq 1} \la   \mb{E} \lb \mc{C}^{(t)}_2 | \mc{F}_{t-1} \rb; \mb{E} \lb \mc{A}^{(t)}_{1} | \mc{F}_{t-1}\rb;\dots;\mb{E} \lb \mc{A}^{(t)}_{m} | \mc{F}_{t-1}\rb  \ra \leq \sigma
\end{align}
almost surely.
\end{assump}
Condition~\eqref{Sub expo Bi} states that the time-dependent coefficients in~\eqref{as-bound2} are conditionally sub-Gaussian with respect to the sigma-fields \( \mathcal{F}_{t-1} \). It is further required that the conditional sub-Gaussian parameter is uniformly bounded above by a deterministic constant \( \sigma \). In other words, we assume a uniform bound on the conditional $2$-Orlicz norm of these coefficients; see \eqref{Orlicz} below for a precise definition. At present, it is unclear whether the results extend to the sub-exponential setting. However, in the homogeneous setting, such as in Proposition~\ref{pre-planned 1} or Theorem~\ref{conf}, an extension is possible when the coefficient corresponding to \( \eta_t^2 \) is sub-exponential with the proof being similar to that of Theorem \ref{theorem 1A} below. 

 Unfortunately, these conditions are not satisfied by Oja's algorithm, since the corresponding recursion (see~\eqref{oja recursive} with \( U_t \equiv Q_t \)) takes the form
\[
|U_t| \lesssim \| \bm{X}_t \|^2 \eta_t \sqrt{L^{Oja}_{t-1}} +  \| \bm{X}_t \|^2 \cdot  (\text{eigen gap}) \cdot \eta_t \lb L^{Oja}_{t-1} \rb^2  + \eta_t^2 \| \bm{X}_t \|^4.
\]
The three coefficients in the recursive bound above are \( \left\{ \| \bm{X}_t \|^2, \| \bm{X}_t \|^2, \| \bm{X}_t \|^4 \right\} \). The last coefficient, \( \| \bm{X}_t \|^4 \), follows a sub-Weibull distribution, even when the data are Gaussian.

Condition~\eqref{bounded Ai} asserts that the time-dependent coefficients in~\eqref{cond-error2} are almost surely bounded by a deterministic constant \( \sigma \). As noted in Section~\ref{sec: main results}, this is a natural extension of~\eqref{as-bound}, and it is satisfied in most examples considered in the present paper. In fact, it is often the case that the coefficients \( \mc{C}^{(t)}_2 \) and \( \left\{ \mc{A}^{(t)}_i : 1 \leq i \leq m \right\} \) are functionals of the collection \( \left\{ \mc{C}^{(t)}_3, \mc{B}^{(t)}_{1}, \mc{B}^{(t)}_{2}, \dots, \mc{B}^{(t)}_{m} \right\} \), and therefore remain bounded.

Under Assumption~\ref{SE Bi}, an analogue of Theorem~\ref{master} holds, which we state below:

\begin{manualtheorem}{1A} \label{theorem 1A}
Assume \eqref{cond-error2}, \eqref{as-bound2}, and suppose that Assumption~\ref{SE Bi} holds for some \( \sigma> 0 \). Assume additionally that
\[
    \min_{1 \leq i \leq m} \left\{ (a_i + b_i) \wedge (c_i + d_i) \right\} > 1.
\]
Then, for any \( \delta \in (0, e^{-2}) \), there exist constants \( M, r > 0 \) independent of $t$, and a step-size sequence \( \{ \eta_t \}_{t \geq 1} \), such that
\[
    \mathbb{P}\left( \forall t \geq 0:\ L_t \leq M \cdot \frac{\log(\delta^{-1}) + \log\log(t+10)}{t+10} \right) \geq 1 - 2\delta,
\]
provided that
\[
    \mathbb{P}(L_0 \leq r) \geq 1 - \delta.
\]
Moreover, the step sizes satisfy \( \eta_t \asymp 1/t \) as \( t \to \infty \).
\end{manualtheorem}

Simnilar to Theorem \ref{master}, for small $\delta$, $M$ scales like $O(1/\log \lb \delta^{-1} \rb)$ and $r$ scales like $O\lb 1/\mbox{polylog} \lb \delta^{-1} \rb  \rb$, with respect to $\sigma$ and the constants $\la a_i, b_i, c_i, d_i; 1 \leq i \leq m \ra$. The proof of Theorem \ref{theorem 1A} is similar to that of Theorem \ref{master}, albeit some changes in the maximal inequality.

\subsection{Extension of Proposition \ref{pre-planned 1} and Theorem \ref{conf}}

\subsection{Some technical results}
We collect some technical results and facts that are needed in this subsection. For $p \geq 1$, define the $p$-Orlicz norm 
 \begin{align} \label{Orlicz}
     \| X \|_{\Psi_p} := \inf \la t>0: \exp \lb \frac{|X|^p}{t^p} \rb \leq 2  \ra.
 \end{align}

 It is well-known that $\|.\|_{\Psi_p}$ is a proper norm for all $p\geq 1$. The following lemma is adapted from \cite{Vershynin} to have explicit constants.
 \begin{lemma} \label{centering Orlicz}
    For $p \in \la 1;2 \ra$, we have
     \[
     \sup_{X} \frac{\| X- \mb{E}X \|_{\Psi_p}}{ \| X \|_{\Psi_p}} \leq 2
     \]
     where the supremum is taken over all random variables $X$ that is not equal to zero with probability one.
 \end{lemma}

 \noindent \textbf{Proof of Lemma \ref{centering Orlicz}}. Consider the case $p=1$. Take $\lambda>0$ and note that
 \begin{align*}
     \mb{E} \left[ \exp \lb \lambda | X - \mb{E}X | \rb  \right] &\leq \mb{E} \lb \exp \lb \lambda|X| + \lambda|\mb{E}X| \rb \rb \\
     & \leq  \exp \lb \lambda |EX| \rb \cdot \mb{E} \lb \exp \lb \lambda|X| \rb \rb \\
     &\leq \left[ \mb{E} \lb e^{\lambda|X|} \rb \right]^2 \leq \mb{E} \lb e^{2\lambda|X|} \rb.
 \end{align*}
where the last line follows from Jensen's inequality. 

Take $\ve>0$ and choose $\lambda= 1/ 2 \lb \|X\|_{\Psi_1}+\ve  \rb$, we get
\begin{align*}
    \mb{E} \left[ \exp \lb \frac{|X-\mb{E}X|}{ 2 \lb \|X\|_{\Psi_1}+\ve  \rb} \rb  \right] \leq \mb{E} \left[ \exp \lb \frac{|X|}{\|X\|_{\Psi_1}+\ve } \rb \right] \leq 2.
\end{align*}
Thus,
\[
\| X - \mb{E}X \|_{\Psi_1} \leq 2 \lb \|X\|_{\Psi_1}+\ve  \rb.
\]
By letting $\ve \to 0$, we conclude the proof for the case $p=1$.

Let us turn to the case $p=2$. By using similar argument to the case $p=1$, we have
\begin{align*}
    \mb{E} \left[ \exp \lb \lambda | X - \mb{E}X |^2 \rb  \right] &\leq \mb{E} \left[ \exp \lb 2\lambda X^2 + 2\lambda \mb{E} (X^2) \rb \right] \\
    &\leq \mb{E}\left[ e^{4\lambda X^2} \right]
\end{align*}
for all $\lambda>0$. By taking $\lambda=1/4\lb \| X \|^2_{\Psi_2} + \ve \rb$ for $\ve>0$, we obtain
\begin{align*}
    \mb{E} \left[ \exp \lb \frac{|X-\mb{E}X|^2}{ 4 \lb \|X\|^2_{\Psi_2}+\ve  \rb} \rb  \right] \leq \mb{E} \left[ \exp \lb \frac{|X|}{\|X\|^2_{\Psi_2}+\ve } \rb \right] \leq 2.
\end{align*}
Thus,
\[
\| X - \mb{E}X \|_{\Psi_2} \leq 2\sqrt{ \|X\|^2_{\Psi_2} +\ve }.
\]
The proof is completed by taking $\ve \to 0$. $\hfill$ $\square$

\begin{lemma} \label{Orlicz equivalence}
    Let $X$ be a random variable such that $\| X \|_{\Psi_p}<\infty$ for $p \in \la  1 ;2     \ra$. Then,
    \begin{itemize}
        \item If $p=1$, then 
        \[
    \mb{E} \left[ \exp \lb \lambda(X-\mb{E}X) \rb   \right] \leq \exp \lb 10\lambda^2 \| X \|^2_{\Psi_1} \rb
    \]
    for all $|\lambda|\leq1/ \lb 4 \| X \|_{\Psi_1} \rb$. 
    \item If $p=2$, then 
        \[
    \mb{E} \left[ \exp \lb \lambda(X-\mb{E}X) \rb   \right] \leq \exp \lb 6\lambda^2 \| X \|^2_{\Psi_2} \rb
    \]
    for all $\lambda \in \mb{R}$. 

    \end{itemize}

\end{lemma}

\noindent \textbf{Proof of Lemma \ref{Orlicz equivalence}}. Put $Y=X-\mb{E}X$. Then, by  Lemma \ref{centering Orlicz}, we have $\| Y \|_{\Psi_p} \leq 2 \| X \|_{\Psi_p}$. Take $\ve>0$ and  define the rescaled version of $Y$
\[
Y_1:= \frac{Y}{2\| X \|_{\Psi_p}+\ve}.
\]
It is obvious that $\| Y_1 \|_{\Psi_p}<1$ and thus, the elementary inequality $e^x \leq 1 + x + (x^2/2)e^{|x|}$ gives
\begin{align} \label{Taylor ineq}
    \mb{E} \left[  \exp\lb  \lambda Y_1 \rb \right] &\leq 1 + \mb{E} \left[  \frac{\lambda^2 Y_1^2}{2} \cdot e^{|\lambda Y_1|}  \right].
\end{align}
Consider the case $p=1$. Observe that $x^2 \leq (16/e^2)e^{|x|/2}$ so \eqref{Taylor ineq} implies
\begin{align*}
        \mb{E} \left[  \exp\lb  \lambda Y_1 \rb \right] &\leq 1 + (8/e^2)\lambda^2 \mb{E} \left[ e^{|Y_1|/2+\lambda|Y_1|} \right]. 
\end{align*}
For all $|\lambda|\leq 1/2$, we have
\begin{align*}
    \mb{E} \left[  \exp\lb  \lambda Y_1 \rb \right] &\leq 1 + (8/e^2)\lambda^2 \mb{E} \lb  e^{|Y_1|}  \rb \leq 1+ (16/e^2)\lambda^2 \leq \exp \lb  (16/e^2)\lambda^2  \rb.
\end{align*}
Consequently,
\begin{align*}
        \mb{E} \left[ \exp \lb \lambda(X-\mb{E}X) \rb   \right] &\leq \exp \lb \frac{16}{e^2}\lb 2 \| X\|_{\Psi_1} +\ve  \rb^2\lambda^2  \rb  
\end{align*}
 for all $\lambda$ satisfies
 \begin{align*}
    |\lambda| \leq  \frac{1}{4 \| X\|_{\Psi_1}+2\ve}. 
 \end{align*}
By letting $\ve \to 0$, we obtain
\begin{align*}
    \mb{E} \left[ \exp \lb \lambda(X-\mb{E}X) \rb   \right] \leq \exp \lb 9\| X\|^2_{\Psi_1}\lambda^2  \rb
\end{align*}
for all $|\lambda| \leq 1/\lb 4 \| X\|_{\Psi_1}\rb$.

Next, consider the case $p=2$. By using the bound $x \leq e^x$ and Cauchy-Schwartz inequality in \eqref{Taylor ineq}, we obtain
\begin{align*}
            \mb{E} \left[  \exp\lb  \lambda Y_1 \rb \right] &\leq 1  + \frac{\lambda^2}{2} e^{\lambda^2/2} \cdot \mb{E} \lb e^{Y_1^2} \rb \\
            &\leq 1 +\lambda^2 e^{\lambda^2/2} \\
            &\leq (1+\lambda^2)e^{\lambda^2/2} \leq e^{(3/2)\lambda^2},
\end{align*}
for all $\lambda \in \mb{R}$. 

Consequently, by letting $\ve \to 0$,
\begin{align*}
    \mb{E} \left[ \exp \lb \lambda(X-\mb{E}X) \rb   \right] \leq \exp \lb 6\| X\|^2_{\Psi_1}\lambda^2  \rb
\end{align*}
for all $\lambda \in \mb{R}$. The proof is completed. $\hfill$ $\square$

 We now modify Lemma \ref{freedman} to account for the random coefficients setting.
\begin{manuallemma}{3A} \label{freedman2}
        Let $M_t$ be an adapted process with the corresponding filtration $\la \mathcal{F}_t; t \geq 0 \ra$ with $M_0=0$. Suppose $X_t$ and $Y_t$ are positive random variables such that
    \begin{align*}
     |M_t - M_{t-1}| \leq X_t, \\
        \left|  \mb{E} \lb M_t - M_{t-1} \Big| \mathcal{F}_{t-1} \rb   \right| \leq Y_t. 
    \end{align*}
Assume additionally that
 \begin{align*}
    \mb{E} \left[ \exp \lb \frac{|X_t|}{c_t} \rb \Big|\mc{F}_{t-1} \right] \leq 2
 \end{align*}
almost surely for all $t \in [1,T]$. Then, 
    $$\mbox{P} \lb \exists t \in [1,T]: M_t \geq  V^{*}_T(\delta) + \sum_{i=1}^{t} Y_i     \rb \leq \delta $$
    for all $T \geq 1 $ and $\delta>0$, where 
    \begin{align*}
        V^{*}_T(\delta) &:= \inf_{|\lambda|<\lambda_0} \la \frac{\log \lb \delta^{-1} \rb}{\lambda} + 40\lambda\sum_{k=1}^{T} c_k^2 \ra; \\
        \lambda_0 &:= \frac{1}{8 \times \max_{1 \leq i \leq T} c_i}.
    \end{align*}
\end{manuallemma} 

\noindent \textbf{Proof of Lemma \ref{freedman2}.} As in the proof of Lemma \ref{freedman}, define $d_0:=0$ and $d_k:= M_{k}-M_{k-1} $. It is easy to check that 
$$M_k= \sum_{i=1}^{k} \left[ d_i - \mb{E} \lb d_i | \mathcal{F}_{i-1} \rb \right] + \sum_{i=1}^{k}  \mb{E} \lb d_i | \mathcal{F}_{i-1} \rb. $$ 
Let $h_i:= d_i - \mb{E} \lb d_i | \mathcal{F}_{i-1} \rb$. It is easy to check that $\la h_i; 1\leq i \leq T \ra$ forms a martingale difference sequence and
\begin{align*}
    \Big| \mb{E} \lb d_i | \mathcal{F}_{i-1} \rb \Big| &\leq Y_i; \\
     |d_i| &\leq X_i 
\end{align*}
almost surely. 

 Conditional on $\mc{F}_{t-1}$, observe that
\begin{align*}
    \| h_i \|_{\Psi_1}  \leq 2\| X_i \|_{\Psi_1} \ \leq 2c_i
\end{align*}
almost surely.

Thus, by Lemma \ref{Orlicz equivalence} and the fact that $\mb{E} \lb h_i  | \mc{F}_{i-1} \rb=0$, we have
\begin{align*}
\mb{E} \left[  e^{\lambda h_i} \Big|  \mc{F}_{i-1}\right] \leq e^{40\lambda^2 c_i^2}
\end{align*}
for all $|\lambda|\leq 1/(8c_i)$.

Therefore,
\begin{align*}
    & \mb{P} \lb \exists t \in [1,T]: M_t \geq x + \sum_{i=1}^{t} Y_i  \rb
    \leq  \mb{P} \lb \exists t \in [1,T]: \sum_{k=1}^{t} h_k \geq x \rb
\end{align*}
for all $x>0$. 

To bound the last probability, note that the process 
\[
 S_t:= \exp \lb \lambda \sum_{k=1}^t h_k - 40\lambda^2 \sum_{k=1}^{t} c_k^2 \rb
\]
forms a non-negative supermartingale, which leads to
\begin{align*}
     \mb{P} \lb \exists t \in [1,T]: \sum_{k=1}^{t} h_k \geq x \rb  &\leq \mb{P} \lb \exists t \in [1,T]: \lambda\sum_{k=1}^{t} h_k - \lambda^2 \sum_{k=1}^{T} c_k^2 \geq \lambda x - \lambda^2 \sum_{k=1}^{T} c_k^2 \rb \\
     &\leq \exp \lb -\lambda x + 40\lambda^2 \sum_{k=1}^{T} c_k^2 \rb
\end{align*}
where the last inequality follows from Ville's inequality. 

Thus, 
\begin{align*}
     \mb{P} \lb \exists t \in [1,T]: \sum_{k=1}^{t} h_k \geq \frac{x}{\lambda}+ 40\lambda \sum_{k=1}^{T} c_k^2 \rb  &\leq e^{-x},
\end{align*}
for all $|\lambda|<1/ \lb 8 \times \max \la c_1;c_2;\dots;c_T \ra \rb$. 

The proof is completed by setting $x= \log \lb \delta^{-1} \rb$. $\hfill$ $\square$

\begin{remark} \label{remark subGauss}
    If one assumes $X_t$ is conditionally sub-Gaussian in Lemma \ref{freedman2}, the variance term $V^{*}_T(\delta)$ is the same as in Lemma \ref{freedman}, up to a constant factor, that is
    \begin{align*}
        V^{SG}_T(\delta) \asymp \sqrt{ \log \lb  \delta^{-1} \rb \sum_{k=1}^{T} c_k^2}.
    \end{align*}
     This corresponds to the sub-Gaussian settings in Theorem \ref{theorem 1A}.
\end{remark}

\subsection{Proof of Theorem \ref{theorem 1A}}
It suffices to show a variant of Lemma \ref{maximal short time} under conditions \eqref{cond-error2} and \eqref{as-bound2}. Once such a result is proven, the stitching argument can be kept unchanged. An analog of Lemma \ref{maximal short time} under \eqref{cond-error2} and \eqref{as-bound2} is
\begin{manuallemma}{4A} \label{maximal short time 2}
    Let $X_0, X_1,\dots,X_n$ be a sequence of random variables that satisfy \eqref{quantitative R-S lemma} with constant step sizes $\eta$. Suppose the noise process in \eqref{quantitative R-S lemma} satisfies \eqref{cond-error2} and \eqref{as-bound2} with the corresponding parameter $\sigma$. Then, there exists a number $D_{\sigma, C_1,m}$ that depends only on $\sigma, C_1$ and $m$ such that with
         \begin{align*} 
     D^{*}_{\delta}(r) :=  \frac{D_{\sigma, C_1,m}}{(1-C_1 \eta)^n} \cdot \lb \sqrt{\log \lb \frac{1}{\delta} \rb}  \cdot \lb \sqrt{\eta r} + \sum_{i=1}^{m} \eta^{c_i} r^{d_i} \rb + \eta  + \sum_{i=1}^{m} \eta^{a_i} r^{b_i} \rb,
     \end{align*}
     we have 
    $$\mb{P} \lb \exists t \in [1,n] : L_t  \geq \lb 1- C_1\eta \rb^{t} \lb \ve_0 + D^{*}_{\delta}(r) \rb, A  \rb \leq\delta $$
    whenever $\ve_0 + D^{*}_{\delta}(r) \leq r $ and $A \subset \la L_{0} \leq \ve_0 \ra$ such that $\mb{P} \lb A \rb > 0$.
\end{manuallemma}

\noindent \textbf{Proof of Lemma \ref{maximal short time 2}}. The only step we need to adapt is the concentration inequality for the stopped process $\la  U^{*}_{\tau \wedge t}; 1 \leq t \leq n \ra$. Recall the notation from the proof of Lemma \ref{maximal short time}. In this case, we still have
Thus,
\begin{align*} 
     \left| U^{*}_{ \tau \wedge t} - U^{*}_{ \tau \wedge (t-1)}  \right|  & \leq u_t  ;  \\
    \left| \mb{E} \lb  U^{*}_{\tau \wedge t} - U^{*}_{\tau \wedge (t-1)} \Big| \mathcal{F}_{t-1} \rb \right| & \leq v_t
\end{align*}
where 
\begin{align*}
    u_t &:= (1- C_1\eta)^{-t} \cdot \left[   \mb{E} \lb \mc{C}^{(t)}_2 | \mc{F}_{t-1} \rb\eta^2 + \sum_{k=1}^{m}  \mb{E} \lb \mc{A}^{(t)}_{i} | \mc{F}_{t-1} \rb\cdot  \eta^{1+a_k} r^{b_k} \right]; \\
    v_t &:= (1- C_1\eta)^{-t} \cdot \lb \mc{C}^{(t)}_3 \eta\sqrt{r} + \sum_{k=1}^{m} \mc{B}^{(t)}_{k} \eta^{1/2+c_k} r^{d_k} \rb.
\end{align*}
By \eqref{as-bound2}, we have
\begin{align*}
    |u_t| \leq  (1- C_1\eta)^{-t} \sigma \cdot  \lb \eta^2+ \sum_{k=1}^{m}  \eta^{1+a_k} r^{b_k} \rb.
\end{align*}
Furthermore, conditionally on $\mc{F}_{t-1}$,
\begin{align*}
\| v_t \|_{\Psi_2} \leq (1- C_1\eta)^{-t} \sigma \lb  \eta\sqrt{r} + \sum_{k=1}^{m}  \eta^{1/2+c_k} r^{d_k}  \rb
\end{align*}
almost surely. Put
\[
c_t:= (1- C_1\eta)^{-t} \sigma \lb  \eta\sqrt{r} + \sum_{k=1}^{m}  \eta^{1/2+c_k} r^{d_k}  \rb.
\]
Therefore, by using Remark \ref{remark subGauss}, we have 
\begin{align*}
    \mb{P} \lb \exists t \in [1,n]: U^{*}_{\tau \wedge t} \geq D^{*}_{\delta}(r) \Big| A \rb 
    \leq & \frac{\delta}{\mb{P} \lb A \rb}
\end{align*}
where
\begin{align*}
 D^{*}_{\delta}(r)  &\stackrel{\sigma, C_1,m}{\asymp} \lb 1 - C_1 \eta \rb^{-n} \lb  \sqrt{\log \lb  \delta^{-1}  \rb \cdot 
 \lb  \sum_{t=1}^{n} c_t^2 \rb} +\eta^2 + \sum_{k=1}^{m} \eta^{1+a_k} r^{b_k} \rb \\
 &\stackrel{\sigma, C_1,m}{\asymp} \lb 1 - C_1 \eta \rb^{-n} \lb \sqrt{\log \lb \frac{1}{\delta} \rb}  \cdot \lb \sqrt{\eta r} + \sum_{i=1}^{m} \eta^{c_i} r^{d_i} \rb + \eta  + \sum_{i=1}^{m} \eta^{a_i} r^{b_i} \rb.
\end{align*}
The rest of the proof can be kept unchanged. $\hfill$ $\square$

Note that the form of $D^{*}_{\delta}(r)$ in Lemma \ref{maximal short time 2} is the same as $D_\delta(r)$ in Lemma \ref{maximal short time}, up to a constant factor. Therefore, the argument in Section \ref{stitching} carries over and this concludes  Theorem \ref{theorem 1A}. $\hfill$ $\square$

    \section{Technical proofs}
  \subsection*{Proof of Lemma \ref{simplified lemma}} Rewrite \eqref{simplified Robbins-Siegmund} in the form of \eqref{almost-supermartingale}, it is easy to check that $a_t=0$, $b_t=\beta_t$ and $c_t=-\eta_t X_t$. Since $\sum_t \beta_t <\infty$ almost surely, the Robbins-Siegmund's lemma can be applied to obtain
    \begin{align*}
        X_t \stackrel{\text{a.s}}{\to} X \  \text{and} \ \sum_{t=1}^{\infty} \eta_t X_t < \infty \ \text{a.s.}
    \end{align*}
    It suffices to show that $X \equiv 0$ a.s. To see this, suppose there exists $\delta > 0 $ such that $\mb{P} \lb X > \delta \rb > \delta$. Then, one can find a set $A$ with probability at least $\delta/2$ such that 
    \[
    \liminf_{t \to \infty} X_t > \delta/2.
    \]
Therefore, the series $ \sum_{t=1}^{\infty} \eta_t X_t$ diverges on $A$, which contradicts the condition $ \sum_{t=1}^{\infty} \eta_t X_t < \infty$ almost surely. The proof is completed. $\hfill$ $\square$


\subsection*{Proof of Proposition \ref{oja-err bound}} In the proof below, $S_t$ stands for $L^{Oja}_t$. Let us start with the error bound for \eqref{krasulina}. For simplicity, we will simply write $\bm{v}_t$ instead of $\bm{\hat{v}}_t$. Define 
\begin{align*}
    y_t &:=  \bm{X}_t^\top \bm{v}_{t-1}, \\
    \bm{z}_t &:= y_t \left[ \bm{X}_t - \frac{y_t}{\| \bm{v}_{t-1}\|^2}\cdot \bm{v}_{t-1}  \right].
\end{align*}
Note that in the update rule of \eqref{krasulina}, $\bm{z}_t$ is orthogonal to $\bm{v}_{t-1}$. To see this, write
\begin{align*}
    & \Big\langle \bm{v}_{t-1}, \eta_t y_t \left[ \bm{X}_t - \frac{y_t}{\| \bm{v}_{t-1}\|^2}\cdot \bm{v}_{t-1}  \right]   \Big\rangle \\
    = & \eta_t y_t \Big\langle \bm{v}_{t-1}, \bm{X}_t - \frac{y_t}{\| \bm{v}_{t-1}\|^2}\cdot \bm{v}_{t-1}   \Big\rangle \\
    = & \eta_t y_t \Big[ \langle \bm{v}_{t-1}, \bm{X}_t \rangle - y_t \Big] = 0.
\end{align*}
From the orthogonality between  $\bm{z}_t$ and $\bm{v}_{t-1}$, one has 
\begin{align*}
    \| \bm{v}_{t} \|^2 &= \| \bm{v}_{t-1} \|^2 + \eta_t^2 y_t^2 \Big\| \bm{X}_t - \frac{y_t}{\| \bm{v}_{t-1}\|^2}\cdot \bm{v}_{t-1} \Big\|^2  \\
    &= \| \bm{v}_{t-1} \|^2 \lb 1 + \eta_t^2 \cdot \frac{ \langle \bm{v}_{t-1}, \bm{X}_t \rangle^2}{\| \bm{v}_{t-1} \|^2} \cdot. \Big\| \bm{X}_t - \frac{y_t}{\| \bm{v}_{t-1}\|^2}\cdot \bm{v}_{t-1} \Big\|^2  \rb. 
\end{align*}
By Cauchy-Schwartz's inequality, it is easy to check that 
\begin{align*}
    \frac{ \langle \bm{v}_{t-1}, \bm{X}_t \rangle^2}{\| \bm{v}_{t-1} \|^2} &\leq \|  \bm{X}_t \|^2 \leq  B^2, \\
     \Big\| \bm{X}_t - \frac{y_t}{ \| \bm{v}_{t-1}\|^2}\cdot \bm{v}_{t-1} \Big\|^2     &\leq 2\|  \bm{X}_t \|^2 +  \frac{2y_t^2}{\| \bm{v}_{t-1}\|^2} \\
     &\leq 2B^2 +   2 \frac{ \langle \bm{v}_{t-1}, \bm{X}_t \rangle^2}{\| \bm{v}_{t-1} \|^2} \leq 4B^2.
\end{align*}
Thus,
\[
    \| \bm{v}_{t-1} \|^2 \leq     \| \bm{v}_{t} \|^2 \leq     \| \bm{v}_{t} \|^2 \lb 1 + 4B^4 \eta_t^2 \rb
\]
and 
\[
\| \bm{z}_t \|^2 \leq 4B^4 \| \bm{v}_{t-1} \|^2
\]
for all $t\geq 1$.

Recall that due to the observation in \eqref{diagonal}, we can assume the principal eigenvector is $\bm{e}_1$. Write 
\begin{align*}
    S_t 
      &= \frac{ \| \bm{v}_t \|^2 - v_{t,1}^2}{ \| \bm{v}_t \|^2} \\
      &\leq \frac{\|\bm{v}_t \|^2 - v_{t,1}^2}{\| \bm{v}_{t-1} \|^2} \\
      &\leq \frac{ \|\bm{v}_{t-1} \|^2 + \eta_t^2 \| \bm{z}_t \|^2 - \lb  v_{t-1,1}^2 + 2\eta_t v_{t-1,1}z_{t,1} + \eta_t^2 z_{t,1}^2 \rb}{\| \bm{v}_{t-1} \|^2} \\
      &= \frac{\|\bm{v}_{t-1} \|^2 - v_{t-1,1}^2}{\| \bm{v}_{t-1} \|^2} + \eta_t^2 \cdot \frac{ \| \bm{z}_t \|^2 - z_{t,1}^2  }{\| \bm{v}_{t-1} \|^2} - \frac{2\eta_t v_{t-1,1}z_{t,1}}{\| \bm{v}_{t-1} \|^2} \\
      &\leq S_{t-1} + 4B^4\eta_t^2 - \frac{2\eta_t v_{t-1,1}z_{t,1}}{\| \bm{v}_{t-1} \|^2}.
\end{align*}
Now define
\begin{align*}
    Q_t:= - \frac{2\eta_t v_{t-1,1}z_{t,1}}{\| \bm{v}_{t-1} \|^2} + \mb{E} \lb  \frac{2\eta_t v_{t-1,1}z_{t,1}}{\| \bm{v}_{t-1} \|^2} \Big| \mathcal{F}_{t-1} \rb.
\end{align*}
By using the orthogonality  between  $\bm{z}_t$ and $\bm{v}_{t-1}$, we have
\begin{align*}
    |v_{t-1,1}z_{t,1}| &=  \Big| \sum_{k=2}^{p}     v_{t-1,k}z_{t,k} \Big| \\
    &\leq \| \bm{z}_t \| \cdot \sqrt{\| \bm{v}_{t-1}\|^2 - v_{t-1,1}^2} \\
    &= \| \bm{z}_t \| \cdot \| \bm{v}_{t-1} \| \cdot \sqrt{S_t}.
\end{align*}
Thus,
\begin{align*}
    |Q_t| &\leq 2\eta_t \sqrt{S_{t-1}} \cdot \frac{\| \bm{z}_t \|}{\| \bm{v}_{t-1} \|} \\
    &\leq 2\eta_t \sqrt{S_{t-1}} \cdot \frac{ |\langle \bm{v}_{t-1}, \bm{X}_t \rangle|}{\| \bm{v}_{t-1} \|} \cdot \Big\| \bm{X}_t - \frac{y_t}{\| \bm{v}_{t-1}\|^2}\cdot \bm{v}_{t-1}  \Big\| \\
    &\leq 2\eta_t \sqrt{S_{t-1}} \cdot B \cdot 2B = 8B^2 \eta_t \sqrt{S_{t-1}}.
\end{align*}
Moreover,
\begin{align*}
    \mb{E} \lb  \frac{2\eta_t v_{t-1,1}z_{t,1}}{ \| \bm{v}_{t-1} \|^2} \Big| \mathcal{F}_{t-1} \rb &= \frac{2\eta_t v_{t-1,1}}{\| \bm{v}_{t-1} \|^2} \cdot \mb{E} \lb z_{t,1} \Big| \mathcal{F}_{t-1}  \rb \\
    &= \frac{2\eta_t v_{t-1,1}}{\| \bm{v}_{t-1} \|^2} \cdot \mb{E} \left[  y_t \lb X_{t,1} - \frac{y_t v_{t-1,1}}{\| \bm{v}_{t-1}\|^2} \rb \Big| \mathcal{F}_{t-1}  \right] \\
    &= \frac{2\eta_t v_{t-1,1}}{\| \bm{v}_{t-1} \|^2} \cdot  \Big[  \lambda_1 v_{t-1,1} - \frac{ v_{t-1,1}}{\| \bm{v}_{t-1}\|^2} \mb{E} \lb y_t^2 \Big| \mathcal{F}_{t-1} \rb \Big] \\
    &= \frac{2\eta_t \lambda_1 v_{t-1,1}^2}{\| \bm{v}_{t-1} \|^2} - \frac{2\eta_t v_{t-1,1}^2}{\| \bm{v}_{t-1} \|^4} \cdot \sum_{i=1}^{p} \lambda_i v_{t-1,i}^2 \\
    &= 2 \eta_t \lambda_1 (1-S_{t-1}) - 2\eta_t (1-S_{t-1}) \left[ \lambda_1 \frac{v_{t-1,1}^2}{\| \bm{v}_{t-1} \|^2} + \sum_{i=2}^{p} \frac{\lambda_i v^2_{t-1,i}}{\| \bm{v}_{t-1} \|^2} \right] \\
    & \geq 2\eta_t \lambda_1 (1-S_{t-1}) - 2\eta_t \lambda_1 (1-S_{t-1})^2  -2\eta_t \lambda_2 S_{t-1}(1-S_{t-1}) \\
    &= 2\eta_t \lambda_1 (1-S_{t-1})S_{t-1} -2\eta_t \lambda_2 S_{t-1}(1-S_{t-1})\\
    &= 2\eta_t(\lambda_1 - \lambda_2)S_{t-1}(1-S_{t-1}).
\end{align*}
Consequently,
\begin{align*}
S_t &\leq S_{t-1} + 4B^4 \eta_t^2  + Q_t - 2\eta_t(\lambda_1 - \lambda_2)S_{t-1}(1-S_{t-1}) \\
&=S_{t-1} \left[ 1 - 2\eta_t(\lambda_1 - \lambda_2) \right] + 2\eta_t(\lambda_1 - \lambda_2)S_{t-1}^2 + Q_t + 4B^4 \eta^2_t.
\end{align*}
This completes the proof of the error bound for \eqref{krasulina}. Let us consider \eqref{Oja}. For this update rule, it is easy to check that
\begin{align*}
    \| \bm{v}_t \|^2 &= \| \bm{v}_{t-1} \|^2 + 2\eta_t y_t^2 + \eta_t^2 \| \bm{X}_t \bm{X}_t^{\top} \bm{v}_{t-1} \|^2 \\
    &= \| \bm{v}_{t-1} \|^2 \lb 1 + \frac{2\eta y_t^2}{\| \bm{v}_{t-1} \|^2} + \eta_t^2 \frac{\| \bm{X}_t \bm{X}_t^{\top} \bm{v}_{t-1} \|^2}{\| \bm{v}_{t-1} \|^2} \rb \\
    &\leq \bm{v}_{t-1} \|^2. \lb 1 + \frac{2\eta y_t^2}{\| \bm{v}_{t-1} \|^2} + \eta_t^2 B^4 \rb.
\end{align*}
Write
\begin{align*}
    S_t &= 1 - \frac{  v_{t,1}^2}{ \| \bm{v}_t \|^2} \\
    & \leq 1 - \lb 1 + \frac{2\eta_t y_t^2}{ \| \bm{v}_{t-1} \|^2} + \eta_t^2 B^4 \rb^{-1} \cdot \frac{ v_{t-1,1}^2 + 2\eta_t y_t v_{t-1,1} X_{t,1} }{\| \bm{v}_{t-1} \|^2} \\
    &\leq 1 - \lb 1 - \frac{2\eta_t y_t^2}{ \| \bm{v}_{t-1} \|^2} - \eta_t^2 B^4 \rb \cdot \frac{ v_{t-1,1}^2 + 2\eta_t y_t v_{t-1,1} X_{t,1} }{ \| \bm{v}_{t-1} \|^2} \\
    &= 1 - \frac{ v_{t-1,1}^2 + 2\eta_t y_t v_{t-1,1} X_{t,1} }{\| \bm{v}_{t-1} \|^2} + 2 \frac{\eta_t y_t^2 v_{t-1,1}^2}{\| \bm{v}_{t-1} \|^4} \\
    &+ 4 \eta_t^2 \cdot \frac{y_t^4}{\| \bm{v}_{t-1} \|^4} + B^4 \eta_t^2 \frac{ v_{t-1,1}^2 + 2\eta_t y_t^2 }{ \| \bm{v}_{t-1} \|^2} \\
    &\leq S_{t-1} -2\eta_t \frac{v_{t-1,1}}{\| \bm{v}_{t-1} \|^2} \cdot \underbrace{y_t \lb X_{t,1} - y_t \frac{v_{t-1,1}}{\| \bm{v}_{t-1} \|^2} \rb}_{z_{t,1}} + \eta_t^2 \lb 4B^4 + B^4(1+2B^2) \rb.
\end{align*}
The middle term in the display above is exactly the same as in $Q_t$ for \eqref{krasulina}. One can use the same bound to proceed. The proof is completed. $\hfill$ $\square$

\subsection*{Proof of Proposition \ref{SC recursive}} 
In the proof below, for simplicity, we use $X_t$ to denote $L^{SGD}_t$. Since $\bm{x}_* \in \mc{X}$, it is a fixed point of $\Pi_{\mc{X}}$. Thus, we can write
\begin{align*}
    X_{t} &= \| \bm{x}_{t} - \bm{x}_* \|^2 \\
    &= \|  \Pi_{\mc{X}} \lb \bm{x}_{t-1} - \eta_{t}  g \lb \bm{x}_{t-1},\bm{\xi}_{t} \rb \rb - \Pi_{\mc{X}} \lb \bm{x}_* \rb  \|^2 \\
    & \leq \|   \bm{x}_{t-1} - \eta_{t} g \lb \bm{x}_{t-1},\bm{\xi}_{t} \rb - \bm{x}_*   \|^2 \\
    &= \|    \bm{x}_{t-1} -    \bm{x}_* \|^2 -2 \eta_{t} \langle g \lb \bm{x}_{t-1},\bm{\xi}_{t} \rb,  \bm{x}_{t-1} -    \bm{x}_*\rangle \\
    &+ \eta_t^2  \|  g \lb \bm{x}_{t-1}, \bm{\xi}_{t} \rb \|^2 \\
    &\leq X_{t-1} -2 \eta_{t} \langle g \lb \bm{x}_{t-1},\bm{\xi}_{t} \rb,  \bm{x}_{t-1} -    \bm{x}_*\rangle + B^2\eta_t^2.
\end{align*}
For the middle term in the display above, write 
\begin{align*}
   & 2\eta_t \langle g \lb \bm{x}_{t-1}, \bm{\xi}_{t} \rb ,  \bm{x}_{t-1} -    \bm{x}_*\rangle\\
   = & 2 \eta_t \langle \nabla F \lb \bm{x}_{t-1} \rb ,  \bm{x}_{t-1} -    \bm{x}_*\rangle + 2\eta_t \underbrace{\langle  g \lb \bm{x}_{t-1},\bm{\xi}_{t} \rb - \nabla F \lb \bm{x}_{t-1}\rb,  \bm{x}_{t-1} -    \bm{x}_*\rangle}_{Y_t} \\
   \geq & 2\eta_t \lb  F \lb \bm{x}_{t-1} \rb -  F \lb \bm{x}_{*}\rb + \frac{\lambda}{2} \| \bm{x}_{t-1} - \bm{x}_* \|^2 \rb  + 2\eta_t Y_t    \\
   \geq &   2\lambda \eta_t \| \bm{x}_{t-1} - \bm{x}_* \|^2 + 2\eta_t Y_t
\end{align*}
where we have used \eqref{strong convex} in the last line. Thus,
\begin{align*}
    X_t \leq \lb 1 - 2\lambda \eta_t \rb X_{t-1} +2\eta_t Y_t + B^2 \eta_t^2
\end{align*}
which gives \eqref{sgd-recursive1}. To finish the proof, note that
\begin{align*}
    \mb{E} \lb Y_t | \mc{F}_{t-1} \rb &= 0, \\
    |Y_t| &\leq \| g - \nabla F  \|_{\infty} \cdot \| \bm{x}_{t-1} - \bm{x}_*  \| \\
    &\leq B \sqrt{X_{t-1}}.
\end{align*}
The proof is completed. $\hfill$ $\square$

\subsection*{Proof of Proposition \ref{P-L recursive}}

    By using the \eqref{smoothness} and \eqref{PL condition}, we have 
\begin{align*}
    F( \bm{x}_{t}) -F( \bm{x}^{*}) &\leq F( \bm{x}_{t-1}) - F(\bm{x}^{*}) - \eta_t \langle \nabla F( \bm{x}_{t-1}), g \lb \bm{x}_{t-1}, \bm{\xi}_t \rb \rangle + \frac{\mu \eta_t^2}{2} \| g \lb \bm{x}_{t-1}, \bm{\xi}_t \rb\|^2 \\
    & \leq F(\bm{x}_{t-1}) - F(\bm{x}^{*}) - \eta_t \| \nabla F(\bm{x}_{t-1}) \|^2 + \eta_t \langle \nabla F(\bm{x}_{t-1}), \nabla F(\bm{x}_{t-1})- g \lb \bm{x}_{t-1}, \bm{\xi}_t \rb \rangle \\
    & +2\mu B^2 \eta_t^2 \\
    & \leq  F(\bm{x}_{t-1}) - F(\bm{x}^{*}) - \tau\eta_t \lb F(\bm{x}_{t-1}) - F(\bm{x}^{*})  \rb + \eta_t Y_t  +2\mu B^2 \eta_t^2 \\
    & \leq (1 -\tau \eta_t) \lb  F(\bm{x}_{t-1}) - F(\bm{x}^{*}) \rb + \eta_t Y_t+  + 2\mu B^2 \eta_t^2,
\end{align*}
where
\begin{align*}
    Y_t:= \langle \nabla F(\bm{x}_{t-1}), \nabla F(\bm{x}_{t-1})- g \lb \bm{x}_{t-1}, \bm{\xi}_t \rb \rangle.
\end{align*}
It is easy to check that 
\begin{align*}
    \mb{E} \lb Y_t | \mc{F}_{t-1} \rb &= 0; \\
    |Y_t| &\leq B \| \nabla F \lb \bm{x}_{t-1} \rb  \| \\
    &\leq B \sqrt{\mu} \cdot  \sqrt{F(\bm{x}_{t-1}) - F(\bm{x}^{*})}
\end{align*}
where the last line follows from the smoothness of $F$. The proof is completed. $\hfill$ $\square$

\subsection*{Proof of Proposition \ref{robbins-monroe constraint}}  
In the proof below, we use $S_t$ to denotes $L^{RM}_t$ for simplicity. Write
\begin{align*}
    S_{t} &= \lb X_{t} - X_{t-1} + X_{t-1} - \theta \rb^2 \\
    &= (X_{t-1} - \theta)^2 + 2(X_{t-1}-\theta)(X_{t}-X_{t-1}) + (X_{t}-X_{t-1})^2 \\
    &=S_{t-1} -2\eta_t(X_{t-1}-\theta)Y(X_{t-1}) + \eta_t^2 Y(X_{t-1})^2 \\
    &= S_{t-1} - 2\eta_t(X_{t-1}-\theta)M(X_{t-1}) - 2\eta_t(X_{t-1}-\theta) \lb Y(X_{t-1}) - M(X_{t-1})  \rb \\
    &+ \eta_t^2 Y(X_{t-1})^2 \\
    &= S_{t-1} - 2\eta_t \cdot S_{t-1} \cdot \frac{M(X_{t-1}) -M(\theta) }{X_{t-1}-\theta} -  2\eta_t(X_{t-1}-\theta) \xi_{X_{t-1}} 
    + \eta_t^2 Y(X_{t-1})^2.
\end{align*}
Note that
\begin{align*}
    \frac{M(X_{t-1}) -M(\theta)}{X_{t-1} - \theta} = M' \lb (1-u) X_{t-1} + u\theta \rb \geq R
\end{align*}
almost surely, where $u \in (0,1)$ is random and depends on $X_{t-1}$. 

Thus,
\[
S_t \leq \lb 1 - 2R\eta_t \rb S_{t-1} - 2\eta_t(X_{t-1}-\theta) \xi_{X_{t-1}} 
    + \eta_t^2 Y(X_{t-1})^2.
\]
Moreover,
\begin{align*}
    Y(X_{t-1})^2 &= \lb M(X_{t-1}) + \xi_{X_{t-1}} \rb^2 \\
    &\leq 2M\lb X_{t-1} \rb^2 + 2R_1^2 \\
    &\leq 2P^2 \lb \sqrt{S_{t-1}} \rb + 2 R_1^2 
\end{align*}
where $R_1$ is given in \eqref{R1}. The proof is completed.
$\hfill$ $\square$

    \begin{lemma} \label{decreasing}
        Let $\delta \in (0,e^{-2})$. Then, the function
        \[
        f(x):= \frac{\log\log(x) + \log(\delta^{-1})}{x}
        \]
        is decreasing on $[e,\infty]$.
    \end{lemma}

    \noindent \textbf{Proof of Lemma \ref{decreasing}.} It is easy to check that the derivative of $f$ is
    \[
    f^{'}(x):= \frac{(\log x)^{-1} - \log \log x - \log(\delta^{-1})}{x^2}.
    \]
    On $[e,\infty)$, we have $1/\log(x) <2 = \log(e^2) \leq \log(\delta^{-1})$ and $\log \log(x)\geq 0$. Thus, $f$ is decreasing. $\hfill$ $\square$

\begin{lemma} \label{tk technical}
    Let  $L$ be an integer greater than $2$ and suppose $\alpha= \max \la L-1, 33 \ra$. Defined by $t_{k}= \alpha^{k+1}-\alpha$. Then,
    \[
    31.5 \times \cdot \frac{ t_{k-1}+L-1}{\lb \alpha^2 -1 \rb \lb t_{k-1} + \alpha^2 \rb} \leq 1
    \]
    for all $k\geq 1$.

\end{lemma}

  \noindent \textbf{Proof of Lemma \ref{tk technical}.} we need to show that
    \begin{align*}
        31.5 \times \lb t_{k-1}+L-1 \rb &\leq \lb \alpha -1 \rb \lb t_{k-1} + \alpha \rb \\
        \iff t_{k-1} \lb 32.5 -\alpha \rb &\leq \alpha \lb \alpha -1 \rb - 31.5\times(L-1).
    \end{align*}
    The statement above is true for $k=1$ since $\alpha \geq \max \la L-1; 32 \ra$ and $t_0=0$. 
    
    For $k\geq 2$, note that the left-hand side in the last display is negative (since $\alpha \geq 33 $) while the right-hand side is positive. The proof is completed. $\hfill$ $\square$

    \begin{lemma} \label{19}
        Let $\Omega= \la (x,y,z) \in \mb{R}^{+}\times \mb{R}^{+} \times \mb{R}^{+} : x+y+z=1 \ra$. Then, we have
        \[
        \min_{\substack{x, y, z > 0 \\ x + y + z = 1}}  \la \max\la \frac{1}{x}; \frac{8}{y}; \frac{16}{z^2} \ra \ra \leq 31.5.
        \]
    \end{lemma}

    \noindent \textbf{Proof of Lemma \ref{19}.} It is easy to see that 
    \[
    \min_{\substack{x, y, z > 0 \\ x + y + z = 1}}  \la \max\la \frac{1}{x}; \frac{8}{y}; \frac{16}{z^2} \ra \ra \leq A
    \]
    where $A$ satisfies
    \begin{align*}
    \begin{cases}
        \frac{1}{x} = \frac{8}{y} =\frac{16}{z^2} = A; \\
        \frac{1}{A} + \frac{8}{A} + \frac{4}{\sqrt{A}} = 1.
    \end{cases}    
    \end{align*}
The second equality is a quadratic equation in $A$, and is equivalent to
\[
A -4\sqrt{A} -9=0.
\]
Solving this equation gives
\[
A= \lb 2 + \sqrt{13} \rb^2 < 31.5.
\]
 $\hfill$ $\square$

    \end{appendix}

\bibliographystyle{plain}
\bibliography{paper-ref}

\end{document}